\documentclass[10pt,reqno]{amsart}
\usepackage{lipsum}
\usepackage{amsmath}
\usepackage{mathtools}
\usepackage{enumerate}
\usepackage{amssymb}
\usepackage{tensor}
\usepackage{etoolbox}

\usepackage{tikz}
\usetikzlibrary{shapes.geometric}
\usetikzlibrary{cd}
\usetikzlibrary{calc}
\tikzcdset{every label/.append style = {font = \small}}

\usepackage{xr-hyper}
\usepackage{hyperref}

\newcommand{\A}{\mathcal{A}}
\newcommand{\B}{\mathcal{B}}
\newcommand{\C}{\mathcal{C}}
\newcommand{\D}{\mathcal{D}}
\newcommand{\E}{\mathcal{E}}

\newcommand{\I}{\mathcal{I}}

\newcommand{\K}{\mathcal{K}}
\newcommand{\M}{\mathcal{M}}
\newcommand{\N}{\mathcal{N}}
\newcommand{\V}{\mathcal{V}}
\newcommand{\X}{\mathcal{X}}
\newcommand{\Tc}{\mathcal{T}}
\newcommand{\Sc}{\mathcal{S}}
\newcommand{\bS}{\normalfont \textbf{S}}
\newcommand{\bR}{\normalfont \textbf{R}}
\newcommand{\bT}{\normalfont \textbf{T}}
\newcommand{\Q}{\mathcal{Q}}
\newcommand{\bX}{\normalfont \textbf{X}}
\newcommand{\bY}{\normalfont \textbf{Y}}
\newcommand{\bZ}{\normalfont \textbf{Z}}

\newcommand{\e}{\epsilon}
\newcommand{\g}{\mathfrak{g}}

\newcommand{\DD}{\overline{\mathcal{D}}}
\newcommand{\MM}{\overline{\mathcal{M}}}

\newcommand{\RC}{\mathcal{RC}}

\newcommand{\RTC}{\mathcal{RTC}}

\newcommand{\TC}{\mathcal{TC}}
\newcommand{\TM}{\mathcal{TM}}
\newcommand{\ModCat}{\B}

\DeclareMathOperator{\CC}{\C \boxtimes \overline{\C}}
\DeclareMathOperator{\ann}{an}
\DeclareMathOperator{\cre}{cr}

\DeclareMathOperator{\Irr}{Irr}
\DeclareMathOperator{\Blank}{--}
\DeclareMathOperator{\End}{End}
\DeclareMathOperator{\Mod}{Mod}
\DeclareMathOperator{\Hom}{Hom}
\DeclareMathOperator{\Rep}{Rep}
\DeclareMathOperator{\Vect}{\normalfont\underline{Vect}}
\DeclareMathOperator{\fld}{\mathbb{K}}

\DeclareMathOperator{\id}{id}

\DeclareMathOperator{\tid}{\normalfont \textbf{1}}
\DeclareMathOperator{\tidtc}{\tid_{\TC}}
\DeclareMathOperator{\vprod}{\cdot}
\DeclareMathOperator{\Tr}{{\normalfont Tr}}
\DeclareMathOperator{\PSL}{{\normalfont PSL}}

\newcommand{\tl}[1][\beta]{\text{\normalfont TL}}
\newcommand{\Rtl}[1][\beta]{\mathcal{R}\text{\normalfont TL}}
\newcommand{\tlred}[1][\beta]{\text{\normalfont TL}^\text{red}}
\newcommand{\Rtlred}[1][\beta]{\mathcal{R}\text{\normalfont TL}^\text{red}}

\newcommand{\Sharp}{^\sharp}

\newcommand{\defeq}{\vcentcolon=}
\newcommand{\LMod}[1]{\Mod \! \mbox{-}#1}
\newcommand{\BMod}[1]{#1,#1\mbox{-} \text{\normalfont Bimod}}
\newcommand{\Dm}[2]{ \tensor[_#1]{\DD}{_#2}}
\newcommand{\Bm}[2]{ \tensor[_#1]{B}{_#2}}

\newcommand{\Bmt}[3]{ \tensor[_#1]{{B^{\otimes #3}}}{_#2}}
\newcommand{\MX}[3]{ \tensor[_#1]{\M(#3)}{_#2}}
\newcommand{\MMX}[3]{ \tensor[_#1]{\overline{\M}(#3)}{_#2}}

\newcommand{\Ostrik}[2]{\Hom_{\M}^{\sigma}(#1,#2)}
\newcommand{\ad}{^{\text{ad}}}

\theoremstyle{plain}
\newtheorem{THM}{Theorem}[section]
\newtheorem{PROP}[THM]{Proposition}
\newtheorem{LEMMA}[THM]{Lemma}

\newtheorem{COR}[THM]{Corollary}
\theoremstyle{definition}
\newtheorem{REM}[THM]{Remark}
\newtheorem{EX}[THM]{Example}
\newtheorem{DEF}[THM]{Definition}
\newtheorem*{DEF*}{Definition}

\title{Extending the Trace of a Pivotal Monoidal Functor}

\usepackage{amsaddr}
\author{Leonard Hardiman}
\address{Institut Camille Jordan\\
Université Claude Bernard Lyon 1\\
43 Boulevard du 11 novembre 1918\\
69622 Villeurbanne Cedex France\\
ORCiD: 0000-0003-1986-6704}
\email{hardiman@math.univ-lyon1.fr}

\begin{document}

    \begin{abstract}

        We consider a pivotal monoidal functor whose domain is a modular tensor category (MTC). We show that the trace of such a functor naturally extends to a representation of the corresponding tube category. As irreducible representations of the tube category are indexed by pairs of simple objects in the underlying MTC, the simple multiplicities of this representation form a candidate modular invariant matrix. In general, this matrix will not be modular invariant, however it will always commute with the T-matrix. Furthermore, under certain additional conditions on the original functor, it is shown that the corresponding representation of the tube category is a haploid, symmetric, commutative Frobenius algebra. Such algebras are known to be connected to modular invariants, in particular a result of Kong and Runkel implies that the matrix of simple multiplicities commutes with the S-matrix if and only if the dimension of the algebra is equal to the dimension of the underlying MTC. Finally, we apply these techniques to certain pivotal monoidal functors arising from module categories over the Temperley-Lieb category and the associated MTC. This provides a novel explanation of the A-D-E pattern appearing in the classification of $ A_1^{(1)} $ modular invariants.

    \end{abstract}

\maketitle

\section{Introduction} \label{sec:introduction}

The principal aim of this article is to describe a procedure which associates an integer non-negative matrix to a pivotal monoidal functor $ \M \colon \C \to \D $. This procedure is motived by the particular case when $ \C $ is the category of modules over a vertex operator algebra and $ \M $ is a module category over $ \C $. In this case $ \M $ may be thought of as describing a boundary conformal field theory and the resulting matrix as the modular invariant of the corresponding closed conformal field theory. This introduction therefore provides a brief review of the relevant mathematical physics followed by a more detailed description of the contents of this paper.

An important property of a (closed) conformal field theory (CFT) is that it has two chiral halves: a holomorphic (or ``left-moving") half and an anti-holomorphic (or ``right-moving") half. In other words, the state space $ H $ of the theory decomposes into the direct sum
   \begin{align} \label{eq:state_space_decomp}
   H = \bigoplus\limits_{IJ} Z_{IJ} \ H_I \otimes \overline{H}_J
   \end{align}
where the $ Z_{IJ} $ are multiplicity spaces and the $ H_I $ range over the irreducible modules of a vertex operator algebra (VOA) $ \V $ (we assume that our CFT is non-heterotic, i.e.\ that $ H_I $ and $ H_J $ are modules over the same VOA). The physical ``uniqueness of the vacuum" assumption imposes that $ Z_{\tid,\tid} = \mathbb{C} $ where $ \tid $ is such that $ H_{\tid} = \V $. The CFT is called \emph{rational} if $ \V $ admits only finitely many irreducible modules; we assume that this is the case from now on. The decomposition of $ H $ given by~\eqref{eq:state_space_decomp} implies that
   \begin{align} \label{eq:partition_function}
   Z(\tau) = \sum\limits_{IJ} \dim Z_{IJ} \ \chi_{I}(\tau) \chi_J(\tau)^*
   \end{align}
where $ Z $ is the \emph{partition function} of the theory, i.e.\ $ Z(\tau) $ is the value of the theory on the torus corresponding to $ \tau \in \mathbb{H} $ and $ \chi_{I} $ is the character of the irreducible VOA module $ H_I $. As conformal structures on the torus are parametrized by $ \mathbb{H} / \PSL_2(\mathbb{Z}) $, we require that $ Z $ be invariant under the action of $ \PSL_2(\mathbb{Z}) $ on $ \mathbb{H} $.

The category of modules over a VOA has an extremely rich structure: it forms a \emph{modular tensor category} (MTC)~\cite{MR2140309}. MTCs possess many nice properties (they are semisimple, rigid, braided...) and in particular they come equipped with a representation of $ \PSL_2(\mathbb{Z}) $ given by their \emph{modular data}. Let $ \I $ be (an indexing set for) a complete set of irreducible objects in an MTC. The modular data is composed of two $ \I \times \I $-matrices known as the S-matrix and the T-matrix; they are denoted by $ \Sc $ and $ \Tc $ respectively. Using the graphical calculus of MTCs, the entries of these matrices are given as follows,
   \begin{align} \label{eq:S_and_T_matrices}
   \Tc_{IJ} \defeq \delta_{I,J}
       \begin{array}{c}
           \begin{tikzpicture}[scale = 0.25]
           \draw[thick] (-1,1) -- (-1,3);
           \draw[thick] (-1,-1) -- (-1,-3);
           \draw[thick] (1,1) to[out = 90, in = 90] (3,1);
           \draw[thick] (1,-1) to[out = -90, in = -90] (3,-1);
           \draw[thick] (3,1) -- (3,-1);
           \draw[thick] (1,1) to[out=-90,in=90] (-1,-1);
           \draw[line width = 0.3cm, white] (-1, 1) to[out=-90,in=90] (1,-1);
           \draw[thick] (-1, 1) to[out=-90,in=90] (1,-1);
           \node at (-1, 4.3) {$ I $};
           \end{tikzpicture}
       \end{array}
   \quad
   \Sc_{IJ} \defeq
       \begin{array}{c}
           \begin{tikzpicture}[scale=0.5]
           \draw [thick](-1.5,0) arc (180:0:1);
           \draw[white,line width = 0.3cm] (0.5,0) circle (1);
           \draw[thick] (0.5,0) circle (1);
           \node at (-2,1) {$ I $};
           \node at (2,1) {$ J $};
           \draw [white,line width = 0.3cm](.5,0) arc (0:-180:1);
           \draw [thick](.5,0) arc (0:-180:1);
           \end{tikzpicture}
       \end{array}.
   \end{align}
The condition that the partition function of a CFT is invariant under the action of $ \PSL_2(\mathbb{Z}) $ on $ \mathbb{H} $ may be rephrased as requiring that the $ \I \times \I $-matrix with entries $ \dim Z_{IJ} $ commutes with the modular data of the category of modules over the relevant VOA. This motivates the following definition.

   \begin{DEF} \label{def:modular_invariant}

   For a modular tensor category with tensor identity $ \tid $ and complete set of simples $ \I $, a \emph{modular invariant} is a non-negative integer $ \I \times \I $-matrix that commutes with the modular data and whose $ (\tid,\tid) $-entry is $ 1 $.

   \end{DEF}

A popular strategy when attempting to classify CFTs is to fix a VOA $ \V $ and search for all compatible partition functions. From the above discussion we see that this is related to finding the modular invariants associated to the MTC of modules over $ \V $. An example in which this has been successfully carried out is provided by the VOA constructed from the affine Lie algebra $ A_1^{(1)} $ together with a positive integer $ k $, via the Sugawara construction~\cite{PhysRev.170.1659}. The category of modules in this case is the category of integrable highest weight modules of $ A_1^{(1)} $ at level $ k $, denoted $ \Rep_k A_1^{(1)} $. In 1986 Cappelli, Itzykson and Zuber classified all possible modular invariants in this context and, to their surprise, the classification followed an A-D-E pattern~\cite{MR918402}. The appearance of this pattern intrigued many researchers in the field and was the subject of much speculation~\cite{MR1767503,MR1907190,MR1936496}. The first explanation of the pattern was provided by an operator algebra technique known as \emph{$ \alpha $-induction}, due to B\"{o}ckenhauer and Evans~\cite{MR1652746}. This technique relates the A-D-E classification of Goodman-de la Harpe-Jones subfactors to Cappelli, Itzykson and Zuber's classification~\cite{MR1867545,Ocneanu99}.

When translating from the operator algebra language to the purely categorical one an inclusion of subfactors corresponds to a module category. Modules categories also have a physical interpretation. In 1989 Cardy showed that the algebraic data of an \emph{annular partition function} in a \emph{boundary} CFT (as opposed to the toroidal partition function $ Z $) is given by a module category over the corresponding MTC~\cite{MR1048596}. From a physical point of view the correspondence between module categories and modular invariants should therefore be thought of as a ``closing up" just as an annulus closes up into a torus. Mathematically we would expect this ``closing up" to correspond to taking the \emph{trace}, in some suitable sense, of the module category. A notion of trace does exist for module categories (and more generally for monoidal functors), however it simply produces a representation of the MTC. A priori it is not at all clear how to associate a non-negative integer $ \I \times \I $-matrix to this representation. This article presents a solution to this problem by extending the representation to take values on the \emph{tube category} of the underlying MTC.

For a spherical fusion category $ \C $, the tube category, denoted $ \TC $, shares the same objects as $ \C $ but has more morphisms i.e.\ $ \Hom_{\C}(X,Y) \leq \Hom_{\TC}(X,Y) $.  The intuition is that whereas morphisms in $ \C $ may be represented graphically as diagrams drawn on a bounded region of the plane, morphisms in $ \TC $ are given by diagrams drawn on a cylinder. Section~\ref{sec:the_tm_construction} describes how, for $ \M\colon \C \to \D $ a pivotal monoidal functor, the trace of $ \M $ naturally extends to a representation of $ \TC $, which we denote $ \TM $. This extension may also be understood from the perspective of work done by Bruguieres and Natale~\cite{MR2863377,MR3161401}; this is explained in greater detail in Remark~\ref{rem:bruguieres_natale}. As irreducible representations of the tube category are indexed by pairs of elements in $ \I $, decomposing $ \TM $ into irreducibles gives a non-negative integer $ \I \times \I $-matrix, $ Z(\TM) $. For $ F $ a representation of $ \TC $, $ F $ is called T-invariant (respectively S-invariant) if $ Z(F) $ commutes with $ \Tc $ (respectively $ \Sc $).

Section~\ref{sec:t_inv} gives a graphical characterisation of T-invariance when $ \C $ is an MTC. In particular, Theorem~\ref{thm:t_invariance} proves that $ F $ is T-invariant if and only if $ F(t_X) = \id_{F(X)} $ where $ t_X \in \End_{\TC}(X) $ is the \emph{twist morphism} on $ X $, see~\eqref{eq:twist}. An immediate corollary of this is that $ \TM $ is T-invariant; this corollary will later be strengthened to Theorem~\ref{thm:t-inv-premod} which only assumes that $ \C $ is pre-modular. Section~\ref{sec:s_inv} starts by showing that, in general, $ Z(\TM) $ fails to be S-invariant. Indeed, when $ \M $ is the identity functor on $ \C $, $ Z(\TM) $ is given by
   \begin{align*}
   Z(\TM)_{IJ} =
       \begin{cases}
       1 \quad \text{if $ I = J = \tid $} \\
       0 \quad \text{else}
       \end{cases}
   \end{align*}
which doesn't commute with the S-matrix in general (this will be explained in greater detail in Example~\ref{ex:tm_not_s_inv}). However, under the assumption that $ \M $ is indecomposable and takes value in a category whose idempotent completion is multifusion, Theorem~\ref{thm:tm_is_frob} proves that $ \TM $ is a haploid, symmetric, commutative, Frobenius algebra. By a \hyperref[thm:kong_runkel]{result of Kong and Runkel}~\cite[Theorem 3.4]{MR2551797} this implies that $ Z(\TM) $ commutes with the S-matrix if and only if the dimension of $ \TM $ is equal to the dimension of $ \C $. This condition on the dimension of $ \TM $ is equivalent to requiring that $$ (\Sc \ Z(\TM) \ \Sc^{-1})_{\tid,\tid} = Z(\TM)_{\tid,\tid} $$ and is therefore always a necessary condition for S-invariance.

Section~\ref{sec:mod_cat_and_alpha_ind} describes a categorical formulation of $ \alpha $-induction given by Ostrik~\cite[Section 5]{MR1976459}. Let $ \C $ be a pre-modular tensor category and let $ \M \colon \C \to \End(\B) $ be a module category. Following \cite{MR1976459}, we define a subspace $ \Ostrik{I}{J^\vee} < \Hom_{\End(\B)}(\M(I),\M(J^\vee)) $ defined by the condition that $ \beta \in \Ostrik{I}{J^\vee} $ satisfies, for all $ X $ in $ \C $,
		\begin{equation*}
			\begin{tikzcd}[row sep=2cm,column sep=2cm,inner sep=1ex]
			\M(I) \otimes \M(X) \arrow[swap]{d}[name=D]{\beta \otimes \id}   \arrow{r}{\M(\overline{\sigma}_{XI})} &  \M(X) \otimes \M(I) \arrow{d}[name=U]{\id \otimes \beta}
			\\
			\M(J^\vee) \otimes \M(X) \arrow{r}{\M(\sigma_{XJ^\vee})} & \M(X) \otimes \M(J^\vee)
			\arrow[to path={(U) node[scale=2,midway] {$\circlearrowleft$}  (D)}]{}
			\end{tikzcd}
		\end{equation*}
where $ \sigma $ and $ \overline{\sigma} $ are the braiding on $ \C $ and its opposite respectively. Ostrik's categorical formulation of $ \alpha $-induction states that, when the dimension of all the objects in $ \C $ are positive, the $ \I \times \I $-matrix with entries given by the dimension of $ \Ostrik{I}{J^\vee} $ is a modular invariant. Theorem~\ref{thm:alpha_induction_equivalence} proves that, when $ \M $ \emph{induces a pivotal structure on its image}, the $ \TM $ construction may be applied and $ Z(\TM) $ will produce the same matrix as $ \alpha $-induction. Furthermore, this application of the $ \TM $ construction to module categories leads us to Corollary~\ref{cor:module_cats_give_frob_alg} which states that, when $ \M $ is an indecomposable module category that induces a pivotal structure on its full image, $ \TM $ is a haploid, symmetric, commutative, Frobenius algebra.

Finally, Section~\ref{sec:case_study} applies the $ \TM $ construction to a class of examples arising from module categories over the Temperley-Lieb category. The Temperley-Lieb category may be thought of as a diagrammatic presentation of the previously discussed category $ \Rep_k A_1^{(1)} $. It is shown that all module categories over $ \Rep_k A_1^{(1)} $ induce a pivotal structure on their full image and so the $ \TM $ construction may be applied. This leads to a new explanation of the A-D-E pattern that appears in the Cappelli-Itzykson-Zuber classification of $ A_1^{(1)} $ modular invariants.

There are pre-existing methods for relating module categories to modular invariant\footnote[2]{As defined in~\cite[Section 6]{MR2430629}, cf.\ Remark~\ref{rem:mod_inv_algebra}.} Frobenius algebras in $ \CC $. Any module category over $ \C $ may be realised (non-uniquely) as the category of modules of an algebra in $ \C $~\cite{MR1976459}. The \emph{full centre construction}~\cite[Definition 4.9]{MR2443266} then associates a modular invariant, commutative, symmetric Frobenius algebra in $ \CC $ to a special (as defined in, for example, \cite{MR2551797}), symmetric Frobenius algebra in $ \C $~\cite[Theorem 3.18]{MR2551797}. Furthermore every modular invariant, commutative, symmetric Frobenius algebra in $ \CC $ may be realised in this way~\cite[Theorem 3.22]{MR2551797}.

The full centre construction may also be described in terms of the module category directly~\cite[Section 3.1]{MR3406516}. Schaumann has worked on characterising the condition that the module category be equivalent to the category of modules of a special symmetric Frobenius algebra purely in terms of the module category itself. In particular he has shown that it is equivalent to requiring that the module category admits a \emph{module trace}~\cite{MR3019263}. It is possible that this could be related to the condition identified in this article: that the module category induce a pivotal structure on its full image.

\vspace{0.8em}

\noindent \textbf{Conventions.} For $ V $ and $ W $ vector spaces, we write ``$ V = W $" to indicate that $ V $ and $ W $ are isomorphic under an isomorphism that should be clear from the context. Unless otherwise specified, a sum over a variable object ranges over a complete set of simple objects. Similarly, unless otherwise specified, a sum over a variable morphism ranges over a basis of the appropriate $ \Hom $-space. All categories are assumed to be enriched over the category of finite dimensional vector spaces. For a category $ \C $ we use $ \RC $ to denote the category of \emph{contravariant} functors from $ \C $ to the category of finite dimensional vector spaces. An object in $ \RC $ is called a \emph{representation} of $ \C $. For an object $ X $ in a pivotal category we use $ X^\vee $ to denote a dual object to $ X $. The corresponding annihilation and creation maps are then denoted $ \ann_X $ and $ \cre_X $ respectively. Much of the work carried out in this article will be done relative to a fixed spherical fusion category $ \C $, when taking tensor products in this category we will omit the ``$ \otimes $" symbol and write $ XY $ for $ X \otimes Y $. However, we will write the ``$\otimes$" symbol when taking a tensor product in any other category. Many of the arguments in this article exploit the graphical calculus of spherical fusion categories. For an exposition of these techniques see, for example,~\cite{hardiman19a}. In particular, the unit and associativity isomorphisms are suppressed as well as the pivotal structure. All diagrams are read top to bottom.

\vspace{0.8em}

\noindent \textbf{Acknowledgements.} The author thanks Alastair King for his guidance during the period this work was carried out. He is also grateful to Ingo Runkel for multiple helpful conversations.

\section{Preliminaries on the Tube Category} \label{sec:preliminaries}

We start by recording some results on the tube category, which will be used throughout; for more details on these results see~\cite{hardiman_king} and~\cite{hardiman19a}. Let $ \fld $ be an algebraically closed field and let $ \C $ be a spherical fusion category over $ \fld $ with complete set of simples $ \I $. The \emph{tube category} of $ \C $, denoted $ \TC $, is a category whose objects coincide with those of $ \C $ and whose $ \Hom $-spaces are given by
    \begin{align*}
    \Hom_{\TC}(X,Y) \defeq \bigoplus\limits_S \Hom_{\C}(SX,YS)
    \end{align*}
    where, as per our conventions, the direct sum ranges over $ \I $ and the monoidal product symbol is suppressed. To depict a morphism in $ \TC $ using the graphical calculus of spherical fusion categories we take $ \alpha \in \Hom_{\C}(GX,YG) $ and write
	\begin{align} \label{eq:tc_morphism}
	\alpha_G =
		\begin{array}{c}
			\begin{tikzpicture}[scale=0.25,every node/.style={inner sep=0,outer sep=-1}]
			\node (v1) at (0,4) {};
			\node (v4) at (0,-4) {};
			\node (v2) at (4,0) {};
			\node (v3) at (-4,0) {};
			\node (v5) at (-2,2) {};
			\node (v6) at (2,-2) {};
			\node (v7) at (2,2) {};
			\node (v11) at (-2,-2) {};
			\node [draw,diamond,outer sep=0,inner sep=.5,minimum size=22,fill=white] (v9) at (0,0) {$ \alpha $};
			\node at (3,3) {$ X $};
			\node at (-3,-3) {$ Y $};
			\node at (-3,3) {$ G $};
			\node at (3,-3) {$ G $};
			\draw [thick] (v9) edge (v6);
			\draw [thick] (v9) edge (v7);
			\draw [thick] (v9) edge (v11);
			\draw [thick] (v5) edge (v9);
			\draw[very thick, red]  (v1) edge (v3);
			\draw[very thick, red]  (v2) edge (v4);
			\draw[very thick]  (v1) edge (v2);
			\draw[very thick]  (v3) edge (v4);
			\end{tikzpicture}
		\end{array}
\end{align}
as shorthand for $ \bigoplus_S\sum_b (\id_Y \otimes \ b^*) \circ \alpha \circ (b \otimes \id_X) \in \Hom_{\TC}(X,Y) $, where $ \{ b \} $ is a basis of $ \Hom_{\C}(S,G) $ and $ \{ b^* \} $ is the corresponding dual basis of $ \Hom_{\C}(G,S) $ with respect to the perfect pairing given by composition into $ \End_{\C}(S) = \fld $, see \cite[Proposition 3.1]{hardiman_king}. The intuition is that whereas morphisms in $ \C $ may be represented graphically as diagrams drawn on a bounded region of the plane, morphisms in $ \TC $ are given  by diagrams drawn on a cylinder. In particular, the red lines in~\eqref{eq:tc_morphism} should be thought of as being glued; this is compatible with our notation as one may indeed show that
	\begin{align*} 
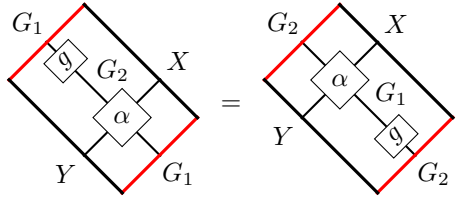
		\begin{array}{c}
			\begin{tikzpicture}[scale=0.25,every node/.style={inner sep=0,outer sep=-1}]
			\node (v1) at (-2,6) {};
			\node (v4) at (0,-4) {};
			\node (v2) at (4,0) {};
			\node (v3) at (-6,2) {};
			\node (v5) at (-4,4) {};
			\node (v6) at (2,-2) {};
			\node (v7) at (2,2) {};
			\node (v11) at (-2,-2) {};
			\node [draw,diamond,outer sep=0,inner sep=.5,minimum size=22,fill=white] (v9) at (0,0) {$ \alpha $};
			\node at (3,3) {$ X $};
			\node at (-3,-3) {$ Y $};
			\node at (-5,5) {$ G_1 $};
			\node at (3,-3) {$ G_1 $};
			\node at (-0.5,2.5) {$ G_2 $};
			\draw [thick] (v9) edge (v6);
			\draw [thick] (v9) edge (v7);
			\draw [thick] (v9) edge (v11);
			\draw [thick] (v5) edge (v9);
			\node [draw,rotate=45,outer sep=0,inner sep=2,minimum width=13,fill=white] (v8) at (-3,3) {$ g $};
			\draw[very thick, red]  (v1) edge (v3);
			\draw[very thick, red]  (v2) edge (v4);
			\draw[very thick]  (v1) edge (v2);
			\draw[very thick]  (v3) edge (v4);
			\end{tikzpicture}
		\end{array}
	=
		\begin{array}{c}
			\begin{tikzpicture}[scale=0.25,every node/.style={inner sep=0,outer sep=-1}]
			\node (v1) at (-2,6) {};
			\node (v4) at (0,-4) {};
			\node (v2) at (4,0) {};
			\node (v3) at (-6,2) {};
			\node (v5) at (-4,4) {};
			\node (v6) at (2,-2) {};
			\node (v7) at (0,4) {};
			\node (v11) at (-4,0) {};
			\node [draw,diamond,outer sep=0,inner sep=.5,minimum size=22,fill=white] (v9) at (-2,2) {$ \alpha $};
			\node at (1,5) {$ X $};
			\node at (-5,-1) {$ Y $};
			\node at (-5,5) {$ G_2 $};
			\node at (3,-3) {$ G_2 $};
			\node at (0.75,1.25) {$ G_1 $};
			\draw [thick] (v9) edge (v6);
			\draw [thick] (v9) edge (v7);
			\draw [thick] (v9) edge (v11);
			\draw [thick] (v5) edge (v9);
			\node [draw,rotate=45,outer sep=0,inner sep=2,minimum width=13,fill=white] (v8) at (1,-1) {$ g $};
			\draw[very thick, red]  (v1) edge (v3);
			\draw[very thick, red]  (v2) edge (v4);
			\draw[very thick]  (v1) edge (v2);
			\draw[very thick]  (v3) edge (v4);
			\end{tikzpicture}
		\end{array}
\end{align*}
for any $ \alpha \in \Hom_{C}(G_2X,YG_1) $ and $ g \in \Hom_{\C}(G_1,G_2) $. Composition in $ \TC $ is then defined following the intuition of vertically stacking the cylinders:
	\begin{align} \label{eq:composition_in_tc}
	\beta_H \circ \alpha_G \defeq \bigoplus\limits_{T} \sum\limits_{b}
		\begin{array}{c}
			\begin{tikzpicture}[scale=0.5,every node/.style={inner sep=0,outer sep=-1}]
			\node (v1) at (-1.25,2.75) {};
			\node (v4) at (-1.25,-5.25) {};
			\node (v2) at (2.75,-1.25) {};
			\node (v3) at (-5.25,-1.25) {};
			\node (v9) at (0.75,0.75) {};
			\node (v6) at (-3.25,-3.25) {};
			\node (v70) at (0.75,-3.25) {};
			\node (v7) at (-3.25,0.75) {};
			\node at (0.75,-1.25) {$ G $};
			\node at (-3.25,-1) {$ H $};
			\node at (1.25,1.25) {$ X $};
			\node at (-3.75,-3.75) {$ Z $};
			\node at (-1.5607,-1.0119) {$ Y $};
			\node at (1.2,-3.7) {$ T $};
			\node at (-3.75,1.25) {$ T $};
			\node (v11) at (-2.625,-1.875) {};
			\node (v12) at (-1.875,-2.625) {};
			\draw [thick] (v11) edge (v12);
			\node (v14) at (-0.625,0.125) {};
			\node (v13) at (0.125,-0.625) {};
			\draw [thick] (v13) edge (v14);
			\node (v15) at (-2.45,-0.425) {};
			\node (v16) at (-2.075,-0.05) {};
			\draw [thick] (v15) to[out=-45,in=135] (v11);
			\draw [thick] (v16) to[out=-45,in=135] (v14);
			\node (v17) at (-0.05,-2.075) {};
			\node (v18) at (-0.425,-2.45) {};
			\draw [thick] (v13) to[out=-45,in=135] (v17);
			\draw [thick] (v12) to[out=-45,in=135] (v18);
			\node [diamond,draw,outer sep=0,inner sep=0.3,minimum size=25,fill=white] (v50) at (-2.25,-2.25) {\mbox{$ \beta $}};
			\node [diamond,draw,outer sep=0,inner sep=-0.2,minimum size=25,fill=white] (v5) at (-0.25,-0.25) {\mbox{$ \alpha $}};
			\node [draw,rotate=45,outer sep=0,inner sep=2,minimum height=10,minimum width=13,fill=white] (v8) at (-2.5,0) {$ b $};
			\node [draw,rotate=45,outer sep=0,inner sep=2,minimum height=10,minimum width=13,fill=white] (v10) at (0,-2.5) {$ b^* $};
			\draw[thick]  (v50) edge (v5);
			\draw [thick] (v50) edge (v6);
			\draw[thick]  (v7) edge (v8);
			\draw[thick]  (v10) edge (v70);
			\draw[thick]  (v9) edge (v5);
			\draw[very thick, red]  (v1) edge (v3);
			\draw[very thick, red]  (v2) edge (v4);
			\draw[very thick]  (v1) edge (v2);
			\draw[very thick]  (v3) edge (v4);
			\end{tikzpicture}
		\end{array}.
    \end{align}
This intuition, together with the associativity of the tensor product, guarantee that composition in $ \TC $ is associative.

    \begin{REM}

    We note that the tensor product is merely weakly associative and yet composition in a category must be strongly associative. However, this is not an issue as the associator isomorphisms will simply modify the basis appearing in~\eqref{eq:composition_in_tc} leaving the composition unchanged.

    \end{REM}

	\begin{REM}

	If we consider the algebra $ \End_{\TC}\left(\bigoplus_{S} S \right) $ we recover Ocneanu's \emph{tube algebra}~\cite{Ocneanu1993}. As $ \bigoplus_{S} S $ is a projective generator in $ \TC $, the tube algebra is Morita equivalent to $ \TC $, i.e. its category of representations is equivalent to $ \RTC $.

	\end{REM}

    \begin{REM} \label{rem:isom_kc_endtid}

    Let $ \K(\C) $ denote the Grothendieck ring of $ \C $ and let $ \K_{\fld}(\C) $ denote $ \K(\C) \otimes_{\mathbb{Z}} \fld $. Then $ \End_{\TC}(\tid) $ and $ \K_{\fld}(\C) $ are canonically isomorphic algebras. Indeed, $ \End_{\TC}(\tid) = \bigoplus_S \End(S) = \bigoplus_S \fld $ is precisely the underlying vector space of $ \K_{\fld}(\C) $. Furthermore, composition in $ \End_{\TC}(\tid) $ corresponds to the tensor product in $ \K_{\fld}(\C) $.

    \end{REM}

	\begin{REM} \label{rem:c_subcategory}

    The canonical inclusion $ \Hom_{\C}(X,Y) \hookrightarrow \Hom_{\TC}(X,Y) $ realises $ \C $ as a wide subcategory of $ \TC $.

    \end{REM}

Remark~\ref{rem:c_subcategory} suggests the following natural question: for a given representation of $ \C $ (i.e.\ an object in $ \RC $) what additional data could be provided to specify a unique extension to an object $ F $ in $ \RTC $? This question is answered in~\cite{hardiman19a} by considering the value of the extended functor on morphisms in $ \TC $ of the form
	\begin{align*}
	\alpha_G =
		\begin{array}{c}
		 	\begin{tikzpicture}[scale=0.2,every node/.style={inner sep=0,outer sep=-1}]
		 	\node (v1) at (0,5) {};
		 	\node (v2) at (5,0) {};
		 	\node (v3) at (-5,0) {};
		 	\node (v4) at (0,-5) {};
			\node (v5) at (3.5,1.5) {};
			\node (v6) at (2.5,-2.5) {};
			\node (v7) at (1.5,3.5) {};
			\node (v8) at (-1.5,-3.5) {};
			\node (v9) at (-2.5,2.5) {};
			\node (v10) at (-3.5,-1.5) {};
		 	\node at (2.75,4.75) {$ X $};
		 	\node at (-2.75,-4.75) {$ X $};
		 	\node at (4.75,2.75) {$ G $};
		 	\node at (-4.75,-2.75) {$ G $};
		 	\node at (3.75,-3.75) {$ G $};
		 	\node at (-3.75,3.75) {$ G $};
			\draw[thick]  (v5) to[out=-135,in=135] (v6);
			\draw[thick] (v7) to[out=-135,in=45] (v8);
			\draw[thick]  (v9) to[out=-45,in=45] (v10);
		 	\draw[very thick, red]  (v1) edge (v3);
		 	\draw[very thick, red]  (v2) edge (v4);
		 	\draw[very thick]  (v1) edge (v2);
		 	\draw[very thick]  (v3) edge (v4);
		 	\end{tikzpicture}
		\end{array}
        \quad \text{where $ \alpha =\id_{GXG} $.}
	\end{align*}
    %


    \begin{PROP}\label{prop:extending_functor}

    Let $ F $ be in $ \RC $ and let $ \kappa_{G,X} \colon F(GX) \to F(XG) $ be a collection of isomorphisms which are natural in $ G $ and $ X $ and satisfy $ \kappa_{H,XG} \circ \kappa_{G,HX} = \kappa_{GH,X} $. Then there is a unique object $ (F,\kappa) $ in $ \RTC $ which satisfies $ (F,\kappa)(X) = F(X) $ for all $ X $ in $ \C $, $ (F,\kappa)(\alpha) = F(\alpha) $ for all $ \alpha \in \Hom_{\C}(X,Y) $ and $ (F,\kappa)(\alpha_G) = \kappa_{G,X} $ where $ \alpha = \id_{GXG} $ for all $ G,X $ in $ \C $.

    \proof

    See Proposition~6.1 in \cite{hardiman19a}. \endproof

    \end{PROP}

    \begin{REM} \label{rem:equiv_zc}

    As $ \C $ is a fusion category the Yoneda embedding gives an equivalence between $ \C $ and $ \RC $. As described in Section~7 of~\cite{hardiman19a}, the data required to extend the image of $ X $ under the Yoneda embedding to $ \TC $ (as given by Proposition~\ref{prop:extending_functor}) corresponds to a half braiding on $ X $. Combining these facts yields an equivalence between $ Z(\C) $  and $ \RTC $, where $ Z(\C) $ is the \emph{Drinfeld centre} of $ \C $.

    \end{REM}

We now equip $ \C $ with a (balanced) \emph{braiding}, in other words, $ \C $ is a \emph{pre-modular tensor category}. Our main tool for studying $ \TC $ in this case will be the following endomorphisms:
	\begin{align} \label{eq:idempotents}
		\e_X^Y = \frac{1}{d(\C)} \bigoplus\limits_S d(S)
			\begin{array}{c}
				\begin{tikzpicture}[scale=0.2,every node/.style={inner sep=0,outer sep=-1}]
				\node (v1) at (0,5) {};
				\node (v4) at (0,-5) {};
				\node (v2) at (5,0) {};
				\node (v3) at (-5,0) {};
				\node (v5) at (-2.5,2.5) {};
				\node (v6) at (2.5,-2.5) {};
				\node (v7) at (1.5,3.5) {};
				\node (v11) at (3.5,1.5) {};
				\node (v12) at (-1.5,-3.5) {};
				\node (v10) at (-3.5,-1.5) {};
				\node (v13) at (2.5,-2.5) {};
				\draw [thick] (v7) edge (v10);
				\draw [line width =0.5em,white] (v5) edge (v13);
				\draw [thick] (v5) edge (v13);
				\draw [line width =0.5em,white] (v11) edge (v12);
				\draw [thick] (v11) edge (v12);
				\node at (2.75,4.75) {$ X $};
				\node at (4.75,2.75) {$ Y $};
				\node at (-3.75,3.75) {$ S $};
				\node at (3.75,-3.5) {$ S $};
				\draw[very thick, red]  (v1) edge (v3);
				\draw[very thick, red]  (v2) edge (v4);
				\draw[very thick]  (v1) edge (v2);
				\draw[very thick]  (v3) edge (v4);
				\end{tikzpicture}
			\end{array}
		\in \End_{\TC}(XY)
	\end{align}
where $ d(S) $ and $ d(\C) $ are the \emph{dimensions} of $ S $ and $ \C $ respectively (for the relevant definitions see~\cite[Definition 4.7.11, 7.21.3]{Etingof15} or \cite[Section 3]{hardiman19a}). As described at the start of~\cite[Section 8]{hardiman19a}, \cite[Proposition~7.4]{hardiman19a} implies that the canonical braided functor
    \begin{align*}
    \Phi \colon \CC \to Z(\C) \cong \RTC
    \end{align*}
satisfies
    \begin{align} \label{eq:phi_and_eXY}
    \Phi(X \boxtimes Y) = (XY,\e_X^Y)\Sharp \defeq \Hom_{\TC}(\Blank,\e_X^Y)
    \end{align}
where $ \Hom_{\TC}(Z,\e_X^Y) = \{\alpha \in \Hom_{\TC}(Z,XY) \mid \e_X^Y \circ \alpha = \alpha \} $. Combining this with the fact that $ \Phi $ is an equivalence when $ \C $ is modular (see \cite[Proposition 8.20.12]{Etingof15}) and we obtain the result that, in the modular case, the set $ \{ \e_I^J \}_{I,J \in \I} $ forms a complete set of orthogonal primitive idempotents in $ \TC $.

	\begin{REM} \label{rem:isomorphic_idempotents}

	The notation $ \e_X^Y $ is chosen (as opposed to $ \e_{XY} $) as $ (XY,\e_X^Y)\Sharp $ is isomorphic (as an object in $ \RTC $) to $ (YX,\tilde{\e}_X^Y)\Sharp $, where
		\begin{align} \label{eq:alt_idem}
		\tilde{\e}_X^Y = \frac{1}{d(\C)} \bigoplus\limits_S d(S)
			\begin{array}{c}
                \begin{tikzpicture}[scale=0.2,every node/.style={inner sep=0,outer sep=-1}]
				\node (v1) at (0,5) {};
				\node (v4) at (0,-5) {};
				\node (v2) at (5,0) {};
				\node (v3) at (-5,0) {};
				\node (v5) at (-2.5,2.5) {};
				\node (v6) at (2.5,-2.5) {};
				\node (v7) at (3.5,1.5) {};
				\node (v11) at (1.5,3.5) {};
				\node (v12) at (-3.5,-1.5) {};
				\node (v10) at (-1.5,-3.5) {};
				\node (v13) at (2.5,-2.5) {};
				\draw [thick] (v7) edge (v10);
				\draw [line width =0.5em,white] (v5) edge (v13);
				\draw [thick] (v5) edge (v13);
				\draw [line width =0.5em,white] (v11) edge (v12);
				\draw [thick] (v11) edge (v12);
				\node at (2.75,4.75) {$ Y $};
				\node at (4.75,2.75) {$ X $};
				\node at (-3.75,3.75) {$ S $};
				\node at (3.75,-3.5) {$ S $};
				\draw[very thick, red]  (v1) edge (v3);
				\draw[very thick, red]  (v2) edge (v4);
				\draw[very thick]  (v1) edge (v2);
				\draw[very thick]  (v3) edge (v4);
				\end{tikzpicture}
			\end{array}
		\in \End_{\TC}(YX).
		\end{align}
	The isomorphism is in fact given by the embedding of the braiding on $ \C $ into $ \TC $. Therefore the isomorphism class of $ \e_X^Y $ is really determined by the fact that the $ X $ strand is \emph{under}-braided and the $ Y $ strand is \emph{over}-braided. This motivates the notation.

	\end{REM}

    \begin{REM} \label{rem:simples_in_kc}

    We recall from Remark~\ref{rem:isom_kc_endtid} that $ \K_{\fld}(\C) = \End_{\TC}(\tid) $. By Corollary 8.2 in~\cite{hardiman19a} we have
     	\begin{align*}
     	\Hom_{\TC}(\tid,\e_I^J) = \Hom_{\C}(\tid,IJ) = \delta_{J,I^\vee} \fld \quad \forall I,J \in \I.
     	\end{align*}
    As $ \{ \e_I^J \}_{I,J \in \I} $ forms a complete set of primitive idempotents in $ \TC $ we may conclude that $ \K_{\fld}(\C) $ is a commutative semisimple algebra generated by a set of primitive orthogonal idempotents indexed by $ \I $.

    \end{REM}

\section{The $ \TM $ Construction} \label{sec:the_tm_construction}

We are now in a good position to define $ \TM $. This definition, together with a graphical description of certain $ \Hom $-spaces into $ \TM $, provides the content of this section. Let $ \C $ be spherical fusion category, let $ \D $ be a pivotal monoidal category and let $ \M \colon \C \to \D $ be a pivotal monoidal functor. When doing graphical calculus in $ \D $ we use blue to depict the image of objects and morphisms in $ \C $ under $ \M $. For example a morphism $ \alpha \in \Hom_{\D}(A,B) $ is depicted in the normal way,
	\begin{align*}
		\begin{array}{c}
			\begin{tikzpicture}[scale=0.15,every node/.style={inner sep=0,outer sep=-1}]
			\node [draw,outer sep=0,inner sep=2,minimum size=10] (v9) at (0,-3) {$ \alpha $};
			\node at (1.7,0.8) {$A$};
			\node at (1.7,-7) {$B$};
			\draw[thick]  (0,2) edge (v9);
			\draw[thick]  (v9) edge (0,-8);
			\end{tikzpicture}
		\end{array}
	\end{align*}
whereas, for $ \beta \in \Hom_{\C}(X,Y) $, we depict $ \M(\beta) \in \Hom_{\D}(\M(Y),\M(X)) $ as
	\begin{align*}
		\begin{array}{c}
			\begin{tikzpicture}[scale=0.15,every node/.style={inner sep=0,outer sep=-1}]
			\node [draw,outer sep=0,inner sep=2,minimum size=10,color=blue] (v9) at (0,-3) {$ \beta $};
			\node[color=blue] at (1.7,0.8) {$Y$};
			\node[color=blue] at (1.7,-7) {$X$};
			\draw[thick,color=blue]  (0,2) edge (v9);
			\draw[thick,color=blue]  (v9) edge (0,-8);
			\end{tikzpicture}
		\end{array}.
	\end{align*}
Composing $ \M $ with the (contravariant) \emph{trace functor} $ \Hom_{\D}(\Blank,\tid) $ gives the following object in $ \RC $
	\begin{align*}
	\Tr \M \colon \C &\to \Vect \\
	X &\mapsto \Hom_{\D}(\M(X),\tid).
	\end{align*}
For $ X $ and $ G $ in $ \C $ we consider the isomorphism
	\begin{align} \label{eq:define_kappa}
    \begin{split}
	\kappa_{G,X}: \Tr \M(GX) &\to \Tr \M(XG)\\
		\begin{array}{c}
            \begin{tikzpicture}[scale=0.15,every node/.style={inner sep=0,outer sep=-1},xscale=-1,yscale=-1]
			\node (v1) at (-2,-0.5) {};
			\node (v3) at (2,-0.5) {};
			\node (v2) at (-2,-3) {};
			\node (v4) at (2,-3) {};
			\draw [thick,blue] (v1) edge (v2);
			\draw [thick,blue] (v3) edge (v4);
			\node [blue] at (-2,-4.5) {$ X $};
			\node [blue] at (2,-4.5) {$ G $};
			\node at (0,4.8) {};
			\node [draw,outer sep=0,inner sep=1,minimum size=10,fill=white] at (0,0) {$ \quad \alpha \quad $};
			\end{tikzpicture}
		\end{array}
	&\mapsto
		\begin{array}{c}
            \begin{tikzpicture}[scale=0.15,every node/.style={inner sep=0,outer sep=-1},xscale=-1,yscale=-1]
			\node [draw,outer sep=0,inner sep=1,minimum size=10] at (0,0) {$ \quad \alpha \quad $};
			\node (v1) at (-2,-0.5) {};
			\node (v3) at (2,-1.5) {};
			\node (v2) at (-2,-3) {};
			\node (v4) at (6,-1.5) {};
			\node (v5) at (-6,0) {};
			\node (v6) at (6,0) {};
			\node (v7) at (-6,-3) {};
			\draw [thick,blue] (v1) edge (v2);
			\draw [thick,blue] (v3) to[out=-90,in=-90] (v4);
			\draw [thick,blue] (v6) to[out=90,in=90] (v5);
			\draw [thick,blue] (v4) edge (v6);
			\draw [thick,blue] (v7) edge (v5);
			\node [blue] at (-6,-4.5) {$ G $};
			\node [blue] at (-2,-4.5) {$ X $};
			\node (v8) at (2,-0.5) {};
			\draw [thick, blue] (v8) edge (v3);
			\node [draw,outer sep=0,inner sep=1,minimum size=10,fill=white] at (0,0) {$ \quad \alpha \quad $};
			\end{tikzpicture}
		\end{array}.
    \end{split}
	\end{align}
As, for $ f $ and $ g $ morphisms in $ \C $,
	\begin{align*}
		\begin{array}{c}
            \begin{tikzpicture}[scale=0.15,every node/.style={inner sep=0,outer sep=-1},xscale=-1,yscale=-1]
			\node (v1) at (-2,-0.5) {};
			\node (v3) at (2,-1.5) {};
			\node  [draw,outer sep=0,inner sep=1,minimum height=15,color=blue,minimum width=9] (v2) at (-2,-4) {$f$};
			\node (v4) at (6,-1.5) {};
			\node (v5) at (-6,0) {};
			\node (v6) at (6,0) {};
			\node  [draw,outer sep=0,inner sep=1,minimum width=9,minimum height=15,color=blue] (v7) at (-6,-4) {$ g $};
			\draw [thick,blue] (v1) edge (v2);
			\draw [thick,blue] (v3) to[out=-90,in=-90] (v4);
			\draw [thick,blue] (v6) to[out=90,in=90] (v5);
			\draw [thick,blue] (v4) edge (v6);
			\draw [thick,blue] (v7) edge (v5);
			\node (v8) at (-2,-7) {};
			\draw [thick, color=blue] (v2) edge (v8);
			\node (v9) at (-6,-7) {};
			\draw [thick, color=blue] (v7) edge (v9);
			\node (v10) at (2,-0.5) {};
			\draw [thick, blue] (v3) edge (v10);
			\node [draw,outer sep=0,inner sep=1,minimum size=10,fill=white] at (0,0) {$ \quad \alpha \quad $};
			\end{tikzpicture}
		\end{array}
	=
		\begin{array}{c}
            \begin{tikzpicture}[scale=0.15,every node/.style={inner sep=0,outer sep=-1},xscale=-1,yscale=-1]
			\node (v1) at (-2,-1) {};
			\node (v3) at (2,-6.5) {};
			\node  [draw,outer sep=0,inner sep=1,minimum height=15,color=blue,minimum width=9] (v2) at (-2,-4) {$f$};
			\node (v4) at (6,-6.5) {};
			\node (v5) at (-6,0) {};
			\node (v6) at (6,0) {};
			\node  [draw,outer sep=0,inner sep=1,minimum width=9,minimum height=15,color=blue] (v7) at (2,-4) {$ g $};
			\draw [thick,blue] (v1) edge (v2);
			\draw [thick,blue] (v3) to[out=-90,in=-90] (v4);
			\draw [thick,blue] (v6) to[out=90,in=90] (v5);
			\draw [thick,blue] (v4) edge (v6);
			\node (v8) at (-2,-8) {};
			\draw [thick, color=blue] (v2) edge (v8);
			\node (v9) at (-6,-8) {};
			\draw [thick, color=blue] (v5) edge (v9);
			\draw [thick, color=blue] (v3) edge (v7);
			\node (v10) at (2,-1) {};
			\draw [thick, color=blue] (v7) edge (v10);
			\node [draw,outer sep=0,inner sep=1,minimum size=10,fill=white] at (0,0) {$ \quad \alpha \quad $};
			\end{tikzpicture}
		\end{array}
	\end{align*}
we have that $ \kappa_{G,X} $ is natural in both $ G $ and $ X $. Furthermore, we have
	\begin{align*}
	\kappa_{G,HX} \circ \kappa_{H,XG} = \kappa_{G,HX} \left(
		\begin{array}{c}
			\begin{tikzpicture}[scale=0.15,every node/.style={inner sep=0,outer sep=-1},xscale=-1,yscale=-1]
			\node (v1) at (-3,-1) {};
			\node (v3) at (3,-1) {};
			\node (v2) at (-3,-3) {};
			\node (v4) at (7,-1.25) {};
			\node (v5) at (-6,0) {};
			\node (v6) at (7,0) {};
			\node (v7) at (-6,-3) {};
			\draw [thick,blue] (v1) edge (v2);
			\draw [thick,blue] (v3) to[out=-90,in=-90] (v4);
			\draw [thick,blue] (v6) to[out=90,in=90] (v5);
			\draw [thick,blue] (v4) edge (v6);
			\draw [thick,blue] (v7) edge (v5);
			\node [blue] at (-6,-4.5) {$ G $};
			\node [blue] at (-3,-4.5) {$ X $};
			\node [blue] at (0,-4.5) {$ H $};
			\node (v8) at (0,-1) {};
			\node (v9) at (0,-3) {};
			\draw [thick, color=blue] (v8) edge (v9);
			\node [draw,outer sep=0,inner sep=1,minimum size=10,fill=white] at (0,0) {$ \quad \alpha \quad $};
			\end{tikzpicture}
		\end{array} \right)
	= \hspace{-0.7em}
		\begin{array}{c}
			\begin{tikzpicture}[scale=0.15,every node/.style={inner sep=0,outer sep=-1},xscale=-1,yscale=-1]
			\node (v1) at (-3,-1) {};
			\node (v3) at (3,-1) {};
			\node (v2) at (-3,-3) {};
			\node (v4) at (7,-1.25) {};
			\node (v5) at (-6,0) {};
			\node (v6) at (7,0) {};
			\node (v7) at (-6,-3) {};
			\draw [thick,blue] (v1) edge (v2);
			\draw [thick,blue] (v3) to[out=-90,in=-90] (v4);
			\draw [thick,blue] (v6) to[out=90,in=90] (v5);
			\draw [thick,blue] (v4) edge (v6);
			\draw [thick,blue] (v7) edge (v5);
			\node [blue] at (-6,-4.5) {$ G $};
			\node [blue] at (-3,-4.5) {$ X $};
			\node [blue] at (-9,-4.5) {$ H $};
			\node (v8) at (0,-1) {};
			\node (v9) at (9.5,-1.25) {};
			\draw [thick, color=blue] (v8) to[out=-90,in=-90] (v9);
			\node [draw,outer sep=0,inner sep=1,minimum size=10,fill=white] at (0,0) {$ \quad \alpha \quad $};
			\node (v10) at (9.5,0) {};
			\draw [thick, color=blue] (v9) edge (v10);
			\node (v12) at (-9,-3) {};
			\node (v11) at (-9,0) {};
			\draw [thick, color=blue] (v10) to[out=90,in=90] (v11);
			\draw [thick, color=blue] (v11) edge (v12);
			\end{tikzpicture}
		\end{array}
	= \kappa_{GH,X}
	\end{align*}
and $ \kappa_{\tid,X} = \id_{\bar{F}(X)} $. We can therefore apply Proposition~\ref{prop:extending_functor} to extend $ \Tr\M $ to a functor on $ \TC $.

    \begin{DEF}

    Let $ \C $ be spherical fusion category, let $ \D $ be a pivotal monoidal category and let $ \M \colon \C \to \D $ be a pivotal monoidal functor. Then $$  \TM \colon \TC \to \Vect  $$ is the functor in $ \RTC $ obtained by applying Proposition~\ref{prop:extending_functor} to $ \Tr \M $ and $ \kappa_{G,X} $ as given by~\eqref{eq:define_kappa}. For a more concrete description of $ \TM $ we consider $ \alpha_G \in \Hom_{\TC}(X,Y) $. Then we have
     	\begin{align*}
     	\TM(\alpha_G) \colon \Hom_{\D}(\M(Y),\tid) &\to \Hom_{\D}(\M(X),\tid)\\
    	\beta &\mapsto
    		\begin{array}{c}
                \begin{tikzpicture}[scale=0.25,every node/.style={inner sep=0,outer sep=-1},yscale=1]
    			\draw [thick, blue] (-5,-5.5) ellipse (3.5 and 3.5);
    			\node [draw,outer sep=0,inner sep=2,minimum size=16,fill=white] (v1) at (-5,-7) {$ \beta $};
    			\node [blue,draw,diamond,outer sep=0,inner sep=2,minimum size=28,fill=white] (v2) at (-2.5,-3) {$\alpha $};
    			\draw [thick, blue] (v1) to[out=90,in=-135] (v2);
    			\node (v3) at (0,0.5) {};
    			\draw [thick, blue] (v2) to[out=45,in=-90] (v3);
    			\node [blue] at (0,-5) {$ G $};
    			\node [blue] at (-5.25,-4.25) {$ Y $};
    			\node [blue] at (-1.25,0) {$ X $};
    			\end{tikzpicture}
    		\end{array}.
     	\end{align*}
    \end{DEF}

    \begin{REM}

    It is possible to define the $ \TM $ construction for functors $ \M \colon \C \to \D $ that are not pivotal by simply adding in the image of the pivotal structure in $ \C $ to~\eqref{eq:define_kappa}. However, doing so would not be compatible with our graphical conventions: the pivotal structure in $ \C $ (which is suppressed from the graphical calculus) would be mapped to a morphism in $ \D $ which could fail to be the corresponding pivotal structure (and thus not be suppressed from the graphical calculus).

    \end{REM}

    \begin{REM}

    When $ \D $ is linear it is straightforward to check that $$  \M = \M_1 \oplus \M_2 \colon X \mapsto \M_1(X) \oplus \M_2(X)  $$ satisfies $ \TM = \TM_1 \oplus \TM_2 $.

    \end{REM}

    We now once again suppose that $ \C $ is equipped with a (balanced) braiding and is therefore a pre-modular tensor category. For an object $ F $ in $ \RTC $ we consider the $ \Hom $-space
 	\begin{align*}
 	F_X^Y \defeq \Hom_{\RTC}((XY,\e_X^Y)\Sharp, F) = \{ \alpha \in F(XY) \mid F(\e_X^Y)(\alpha) = \alpha \}.
 	\end{align*}
	\begin{PROP} \label{prop:characterising_tm_xy}

	$ \TM_X^Y $ is given by the subspace of $ \TM(XY) = \Hom_{\D}(\M(XY),\tid) $ defined by the condition that $ \alpha \in \Hom_{\D}(\M(XY),\tid) $ satisfy
		\begin{align} \label{eq:phase_through}
			\begin{array}{c}
				\begin{tikzpicture}[scale=0.15,every node/.style={inner sep=0,outer sep=-1},yscale=-1]
				\node (v1) at (-1,1) {};
				\node [blue] (v11) at (2,6) {$ Z $};
				\node [blue] (v10) at (-6,-8) {$ Z $};
				\node [blue] (v11) at (-2,-8) {$ X $};
				\node [blue] (v10) at (2,-8) {$ Y $};
				\node (v9) at (-2,-6) {};
				\node (v13) at (-5,1) {};
				\node (v2) at (2,4) {};
				\node (v3) at (2,1) {};
				\node (v5) at (2,-6) {};
				\node (v4) at (-6,-6) {};
				\draw [thick, blue] (v2) edge (v3);
				\draw [thick, blue] (v9) to[out=90, in=-90] (v13);
				\draw [line width= 0.15cm,white] (v3) to[out=-90, in=90] (v4);
				\draw [thick, blue] (v3) to[out=-90, in=90] (v4);
				\draw [line width= 0.15cm,white] (v1) to[out=-90, in=90] (v5);
				\draw [thick, blue] (v1) to[out=-90, in=90] (v5);
				\node [draw,outer sep=0,inner sep=2,minimum width=22,minimum height=12,fill=white] (v10) at (-3,2) {$ \alpha $};
				\end{tikzpicture}
			\end{array}
		=
			\begin{array}{c}
				\begin{tikzpicture}[scale=0.15,every node/.style={inner sep=0,outer sep=-1},yscale=-1]
				\node (v1) at (2,1) {};
				\node [blue] (v11) at (-6,6) {$ Z $};
				\node [blue] (v10) at (-6,-8) {$ Z $};
				\node [blue] (v11) at (-2,-8) {$ X $};
				\node [blue] (v10) at (2,-8) {$ Y $};
				\node (v9) at (-2,-6) {};
				\node (v13) at (-2,1) {};
				\node (v2) at (-6,4) {};
				\node (v3) at (-6,1) {};
				\node (v5) at (2,-6) {};
				\node (v4) at (-6,-6) {};
				\draw [thick, blue] (v2) edge (v3);
				\draw [thick, blue] (v9) to[out=90, in=-90] (v13);
				\draw [line width= 0.15cm,white] (v3) to[out=-90, in=90] (v4);
				\draw [thick, blue] (v3) to[out=-90, in=90] (v4);
				\draw [line width= 0.15cm,white] (v1) to[out=-90, in=90] (v5);
				\draw [thick, blue] (v1) to[out=-90, in=90] (v5);
				\node [draw,outer sep=0,inner sep=2,minimum width=22,minimum height=12,fill=white] (v10) at (0,2) {$ \alpha $};
				\end{tikzpicture}
			\end{array}.
		\end{align}
	for all $ Z $ in $ \C $.
	\proof

	Evaluating $ \TM $ on $ \e_X^Y $ gives the map
		\begin{align*}
		\TM (\e_X^Y) \colon \Hom_{\D}(\M(XY),\tid) &\to \Hom_{\D}(\M(XY),\tid) \\
		\alpha &\mapsto \frac{1}{d(\C)} \sum\limits_S d(S)
			\begin{array}{c}
				\begin{tikzpicture}[scale=0.15,every node/.style={inner sep=0,outer sep=-1},yscale=-1]
				\node (v1) at (-2,-0.5) {};
				\node (v3) at (2,-0.5) {};
				\node (v2) at (-2,-6) {};
				\node (v4) at (2,-6) {};
				\node [blue] at (-2,-7.5) {$ X $};
				\node [blue] at (2,-7.5) {$ Y $};
				\node [blue] at (5.5,0) {$ S $};
				\draw [thick,blue] (v1) edge (v2);
				\draw [line width= 0.15cm,white] (0,0) ellipse (4 and 4);
				\draw [thick, blue] (0,0) ellipse (4 and 4);
				\draw [line width= 0.15cm,white] (v3) edge (v4);
				\draw [thick, blue] (v3) edge (v4);
				\node [draw,outer sep=0,inner sep=1,minimum height=12,minimum width=24,fill=white] (v5) at (0,0) {$  \alpha $};
				\end{tikzpicture}
			\end{array}.
		\end{align*}
	Therefore, if $ \alpha $ satisfies~\eqref{eq:phase_through}, we have
		\begin{align*}
		\TM (\e_X^Y)(\alpha) = \frac{1}{d(\C)} \sum\limits_S d(S)
			\begin{array}{c}
				\begin{tikzpicture}[scale=0.15,every node/.style={inner sep=0,outer sep=-1},yscale=-1]
				\node (v1) at (-2,-0.5) {};
				\node (v3) at (2,-0.5) {};
				\node (v2) at (-4.5,-6) {};
				\node (v4) at (-0.5,-6) {};
				\node [blue] at (-4.5,-7.5) {$ X $};
				\node [blue] at (-0.5,-7.5) {$ Y $};
				\node [blue] at (-3.5,3.5) {$ S $};
				\draw [thick,blue] (v1) to[out=-90,in=90] (v2);
				\draw [thick, blue] (-8,0) ellipse (4 and 4);
				\draw [thick, blue] (v3) to[out=-90,in=90] (v4);
				\node [draw,outer sep=0,inner sep=1,minimum height=12,minimum width=24,fill=white] (v5) at (0,0) {$  \alpha $};
				\end{tikzpicture}
			\end{array}
		 = \alpha.
		\end{align*}
	Furthermore, for $ \alpha \in \TM_X^Y $, we have
		\begin{align*}
			\begin{array}{c}
				\begin{tikzpicture}[scale=0.15,every node/.style={inner sep=0,outer sep=-1},yscale=-1]
				\node (v1) at (-1,3) {};
				\node [blue] at (3,11) {$ Z $};
				\node [blue] at (-7,-8) {$ Z $};
				\node [blue] at (-2,-8) {$ X $};
				\node [blue] at (2,-8) {$ Y $};
				\node (v9) at (-2,-6) {};
				\node (v13) at (-5,3) {};
				\node (v2) at (3,9) {};
				\node (v3) at (3,1) {};
				\node (v5) at (2,-6) {};
				\node (v4) at (-7,-6) {};
				\draw [thick, blue] (v2) edge (v3);
				\draw [thick, blue] (v9) to[out=90, in=-90] (v13);
				\draw [line width= 0.15cm,white] (v3) to[out=-90, in=90] (v4);
				\draw [thick, blue] (v3) to[out=-90, in=90] (v4);
				\draw [line width= 0.15cm,white] (v1) to[out=-90, in=90] (v5);
				\draw [thick, blue] (v1) to[out=-90, in=90] (v5);
				\node [draw,outer sep=0,inner sep=2,minimum width=22,minimum height=12,fill=white] (v100) at (-3,4) {$ \alpha $};
				\end{tikzpicture}
			\end{array}
		&= \frac{1}{d(\C)} \sum\limits_S d(S)
			\begin{array}{c}
				\begin{tikzpicture}[scale=0.15,every node/.style={inner sep=0,outer sep=-1},yscale=-1]
				\node (v1) at (-1,3) {};
				\node [blue] at (3,11) {$ Z $};
				\node [blue] at (-7,-8) {$ Z $};
				\node [blue] at (-2,-8) {$ X $};
				\node [blue] at (2,-8) {$ Y $};
				\node [blue] at (-9,4) {$ S $};
				\node (v9) at (-2,-6) {};
				\node (v13) at (-5,3) {};
				\node (v2) at (3,9) {};
				\node (v3) at (3,1) {};
				\node (v5) at (2,-6) {};
				\node (v4) at (-7,-6) {};
				\draw [thick, blue] (v2) edge (v3);
				\draw [thick, blue] (v9) to[out=90, in=-90] (v13);
				\draw [line width= 0.15cm,white] (v3) to[out=-90, in=90] (v4);
				\draw [thick, blue] (v3) to[out=-90, in=90] (v4);
				\draw [line width= 0.15cm,white] (-3,4) ellipse (4 and 4);
				\draw [thick, blue] (-3,4) ellipse (4 and 4);
				\draw [line width= 0.15cm,white] (v1) to[out=-90, in=90] (v5);
				\draw [thick, blue] (v1) to[out=-90, in=90] (v5);
				\node [draw,outer sep=0,inner sep=2,minimum width=22,minimum height=12,fill=white] (v100) at (-3,4) {$ \alpha $};
				\end{tikzpicture}
			\end{array}\\
		&= \frac{1}{d(\C)} \sum\limits_{S,T,b} d(S)
			\begin{array}{c}
				\begin{tikzpicture}[scale=0.15,every node/.style={inner sep=0,outer sep=-1},yscale=-1]
				\node (v1) at (-1,4.5) {};
				\node [blue] at (4.5,14) {$ Z $};
				\node [blue] at (-7,-8) {$ Z $};
				\node [blue] at (-2,-8) {$ X $};
				\node [blue] at (2,-8) {$ Y $};
				\node [blue] at (-9,5.5) {$ S $};
				\node [blue] at (5,5.5) {$ T $};
				\node (v9) at (-2,-6) {};
				\node (v13) at (-5,4.5) {};
				\node (v2) at (3.5,12) {};
				\node (v3) at (3.5,1.5) {};
				\node (v5) at (2,-6) {};
				\node (v4) at (-7,-6) {};
				\node (v6) at (-7,1.5) {};
				\node (v7) at (1.5,1.5) {};
				\draw [thick, blue] (v9) to[out=90, in=-90] (v13);
				\draw [line width= 0.15cm,white] (v3) to[out=-90, in=90] (v4);
				\draw [thick, blue] (v3) to[out=-90, in=90] (v4);
				\draw [line width= 0.15cm,white] (v6) to[out=-90,in=-90] (v7);
				\draw [thick, blue] (v6) to[out=-90,in=-90] (v7);
				\draw [line width= 0.15cm,white] (v1) to[out=-90, in=90] (v5);
				\draw [thick, blue] (v1) to[out=-90, in=90] (v5);
				\node [draw,outer sep=0,inner sep=2,minimum width=22,minimum height=12,fill=white] at (-3,5.5) {$ \alpha $};
				\node (v11) at (3.5,9.5) {};
				\node (v14) at (1.5,9.5) {};
				\draw [thick, blue] (v2) edge (v11);
				\node (v12) at (-7,9.5) {};
				\draw [thick, blue] (v12) to[out=90,in=90] (v14);
				\draw [thick, blue] (v12) edge (v6);
				\node [blue,draw,outer sep=0,inner sep=2,minimum width=14,minimum height=12,fill=white] (v8) at (2.5,2.5) {$ b $};
				\node [blue,draw,outer sep=0,inner sep=2,minimum width=14,minimum height=12,fill=white] (v10) at (2.5,8.5) {$ b^* $};
				\draw [thick, blue] (v8) edge (v10);
				\end{tikzpicture}
			\end{array}
        \end{align*}
        \begin{align*}
		&= \frac{1}{d(\C)} \sum\limits_{S,T,b} d(S)
			\begin{array}{c}
                \begin{tikzpicture}[scale=0.15,every node/.style={inner sep=0,outer sep=-1},yscale=-1]
				\node (v1) at (-1,4.5) {};
				\node [blue] at (-14,15) {$ Z $};
				\node [blue] at (-14,-8) {$ Z $};
				\node [blue] at (-5,-8) {$ X $};
				\node [blue] at (-1,-8) {$ Y $};
				\node [blue] at (-7,5.5) {$ S $};
				\node [blue] at (3,5.5) {$ T $};
				\node (v9) at (-5,-6) {};
				\node (v13) at (-5,4.5) {};
				\node (v2) at (-14,13) {};
				\node (v3) at (1,1.5) {};
				\node (v5) at (-1,-6) {};
				\node (v4) at (-14,-6) {};
				\node (v6) at (-8.5,3.5) {};
				\node (v7) at (-10.5,3.75) {};
				\node (v11) at (-9.5,1.5) {};
				\node (v14) at (-9.5,9.5) {};
				\node (v12) at (-8.5,7.5) {};
				\draw [thick, blue] (v9) to[out=90, in=-90] (v13);
				\draw [line width= 0.15cm,white] (v11) to[out=-90,in=-90] (v3);
				\draw [thick, blue] (v11) to[out=-90,in=-90] (v3);
				\draw [line width= 0.15cm,white] (v1) to[out=-90, in=90] (v5);
				\draw [thick, blue] (v1) to[out=-90, in=90] (v5);
				\draw [thick, blue] (v12) edge (v6);
				\node (v15) at (1,9.5) {};
				\draw [thick, blue] (v15) to[out=90, in=90] (v14);
				\draw [thick, blue] (v3) edge (v15);
				\node (v18) at (-14,3.75) {};
				\node (v17) at (-10.5,7) {};
				\node (v16) at (-14,7) {};
				\draw [thick, blue] (v16) to[out=-90, in=-90] (v17);
				\draw [thick, blue] (v7) to[out=90, in=90] (v18);
				\node (v20) at (-10.5,3.25) {};
				\node (v19) at (-10.5,7.75) {};
				\draw [thick, blue] (v19) edge (v17);
				\draw [thick, blue] (v20) edge (v7);
				\draw [thick, blue] (v18) edge (v4);
				\draw [thick, blue] (v16) edge (v2);
				\node [draw,outer sep=0,inner sep=2,minimum width=22,minimum height=12,fill=white] at (-3,5.5) {$ \alpha $};
				\node [blue,draw,outer sep=0,inner sep=2,minimum width=14,minimum height=12,fill=white,rotate=180] (v8) at (-9.5,2.5) {$ b $};
				\node [blue,draw,outer sep=0,inner sep=2,minimum width=14,minimum height=12,fill=white,rotate=180] (v10) at (-9.5,8.5) {$ b^* $};
				\end{tikzpicture}
			\end{array}\\
		&= \frac{1}{d(\C)} \sum\limits_T d(T)
			\begin{array}{c}
				\begin{tikzpicture}[scale=0.15,every node/.style={inner sep=0,outer sep=-1},yscale=-1]
				\node (v1) at (2,-2) {};
				\node [blue] (v11) at (-6,6) {$ Z $};
				\node [blue] (v10) at (-6,-8) {$ Z $};
				\node [blue] (v11) at (-2,-8) {$ X $};
				\node [blue] (v10) at (2,-8) {$ Y $};
				\node [blue] (v10) at (6,-1) {$ T $};
				\node (v9) at (-2,-6) {};
				\node (v13) at (-2,-2) {};
				\node (v2) at (-6,4) {};
				\node (v3) at (-6,1) {};
				\node (v5) at (2,-6) {};
				\node (v4) at (-6,-6) {};
				\draw [thick, blue] (v2) edge (v3);
				\draw [thick, blue] (v9) to[out=90, in=-90] (v13);
				\draw [line width= 0.15cm,white] (v3) to[out=-90, in=90] (v4);
				\draw [thick, blue] (v3) to[out=-90, in=90] (v4);
				\draw [line width= 0.15cm,white] (0,-1) ellipse (4 and 4);
				\draw [thick, blue] (0,-1) ellipse (4 and 4);
				\draw [line width= 0.15cm,white] (v1) to[out=-90, in=90] (v5);
				\draw [thick, blue] (v1) to[out=-90, in=90] (v5);
				\node [draw,outer sep=0,inner sep=2,minimum width=22,minimum height=12,fill=white] (v10) at (0,-1) {$ \alpha $};
				\end{tikzpicture}
			\end{array}
  		=
			\begin{array}{c}
				\begin{tikzpicture}[scale=0.15,every node/.style={inner sep=0,outer sep=-1},yscale=-1]
				\node (v1) at (2,-2) {};
				\node [blue] (v11) at (-6,6) {$ Z $};
				\node [blue] (v10) at (-6,-8) {$ Z $};
				\node [blue] (v11) at (-2,-8) {$ X $};
				\node [blue] (v10) at (2,-8) {$ Y $};
				\node (v9) at (-2,-6) {};
				\node (v13) at (-2,-2) {};
				\node (v2) at (-6,4) {};
				\node (v3) at (-6,1) {};
				\node (v5) at (2,-6) {};
				\node (v4) at (-6,-6) {};
				\draw [thick, blue] (v2) edge (v3);
				\draw [thick, blue] (v9) to[out=90, in=-90] (v13);
				\draw [line width= 0.15cm,white] (v3) to[out=-90, in=90] (v4);
				\draw [thick, blue] (v3) to[out=-90, in=90] (v4);
				\draw [line width= 0.15cm,white] (v1) to[out=-90, in=90] (v5);
				\draw [thick, blue] (v1) to[out=-90, in=90] (v5);
				\node [draw,outer sep=0,inner sep=2,minimum width=22,minimum height=12,fill=white] (v10) at (0,-1) {$ \alpha $};
				\end{tikzpicture}
			\end{array}
		\end{align*}
	where, to make certain string manipulations clearer, we have chosen to write $ b $ and $ b^* $ upside-down instead of writing $ b^\vee $ and $ (b^*)^\vee $ and the penultimate equality uses \cite[Lemma 3.11]{hardiman_king}. \endproof

	\end{PROP}

	\begin{REM} \label{rem:alt_phase_through}

	We recall from Remark~\ref{rem:isomorphic_idempotents} that $ (XY,\e_X^Y)\Sharp = (YX,\tilde{\e}_X^Y)\Sharp $ where $ \tilde{\e}_X^Y $ is given by~\eqref{eq:alt_idem}. Therefore $ \TM_X^Y $ may also be identified with the subspace of $ \Hom_{\D}(\M(YX),\tid) $ defined by the condition that $ \alpha \in \Hom_{\D}(\M(YX),\tid) $  satisfy
		\begin{align*}
			\begin{array}{c}
				\begin{tikzpicture}[scale=0.15,every node/.style={inner sep=0,outer sep=-1},yscale=-1]
				\node (v1) at (-5,1) {};
				\node [blue] (v11) at (2,6) {$ Z $};
				\node [blue] (v10) at (-6,-8) {$ Z $};
				\node [blue] (v11) at (2,-8) {$ X $};
				\node [blue] (v10) at (-2,-8) {$ Y $};
				\node (v9) at (2,-6) {};
				\node (v13) at (-1,1) {};
				\node (v2) at (2,4) {};
				\node (v3) at (2,1) {};
				\node (v5) at (-2,-6) {};
				\node (v4) at (-6,-6) {};
				\draw [thick, blue] (v2) edge (v3);
				\draw [thick, blue] (v9) to[out=90, in=-90] (v13);
				\draw [line width= 0.15cm,white] (v3) to[out=-90, in=90] (v4);
				\draw [thick, blue] (v3) to[out=-90, in=90] (v4);
				\draw [line width= 0.15cm,white] (v1) to[out=-90, in=90] (v5);
				\draw [thick, blue] (v1) to[out=-90, in=90] (v5);
				\node [draw,outer sep=0,inner sep=2,minimum width=22,minimum height=12,fill=white] (v10) at (-3,2) {$ \alpha $};
				\end{tikzpicture}
			\end{array}
		=
			\begin{array}{c}
				\begin{tikzpicture}[scale=0.15,every node/.style={inner sep=0,outer sep=-1},yscale=-1]
				\node (v1) at (2,1) {};
				\node [blue] (v11) at (-6,6) {$ Z $};
				\node [blue] (v10) at (-6,-8) {$ Z $};
				\node [blue] (v11) at (-2,-8) {$ Y $};
				\node [blue] (v10) at (2,-8) {$ X $};
				\node (v9) at (-2,-6) {};
				\node (v13) at (-2,1) {};
				\node (v2) at (-6,4) {};
				\node (v3) at (-6,1) {};
				\node (v5) at (2,-6) {};
				\node (v4) at (-6,-6) {};
				\draw [thick, blue] (v2) edge (v3);
				\draw [thick, blue] (v9) to[out=90, in=-90] (v13);
				\draw [line width= 0.15cm,white] (v3) to[out=-90, in=90] (v4);
				\draw [thick, blue] (v3) to[out=-90, in=90] (v4);
				\draw [line width= 0.15cm,white] (v1) to[out=-90, in=90] (v5);
				\draw [thick, blue] (v1) to[out=-90, in=90] (v5);
				\node [draw,outer sep=0,inner sep=2,minimum width=22,minimum height=12,fill=white] (v10) at (0,2) {$ \alpha $};
				\end{tikzpicture}
			\end{array}
		\end{align*}
	for all $ Z $ in $ \C $.
	\end{REM}

	\begin{DEF} \label{def:zf}

	For any $ F $ in $ \RTC $ one may consider the $ \I \times \I $ integer matrix.
		\begin{align*}
		Z(F) \defeq (\dim F_I^J)_{I,J \in \I}.
		\end{align*}
	\end{DEF}

	\begin{REM}

	We recall that, if $ \C $ is modular, the set $ \{ (IJ,\e_I^J)\Sharp \}_{I,J \in \I} $ forms a complete set of simples in $ \RTC $. Therefore an entry of $ Z(F) $ simply gives the multiplicity of the corresponding simple object in $ F $.

	\end{REM}

\section{T-Invariance} \label{sec:t_inv}

	\begin{DEF} \label{def:t-invariance}

	Let $ \C $ be a pre-modular tensor category and let $ \Tc $ be the T-matrix of $ \C $ as defined by~\eqref{eq:S_and_T_matrices}. We call an object $ F $ in $ \RTC $ \emph{T-invariant} if $ Z(F) $ commutes with $ \Tc $.

	\end{DEF}

The principal goal of this section is to give a graphical characterisation of T-invariance when $ \C $ is modular. We consider the following automorphism of $ X $ in $ \TC $,
	\begin{align} \label{eq:twist}
	t_X \defeq
		\begin{array}{c}
			\begin{tikzpicture}[scale=0.25,every node/.style={inner sep=0,outer sep=-1}]
		 	\node (v1) at (0.5,4.5) {};
		 	\node (v2) at (4.5,0.5) {};
		 	\node (v3) at (-4.5,-0.5) {};
		 	\node (v4) at (-0.5,-4.5) {};
			\node (v6) at (-2,2) {};
			\node (v7) at (2.5,2.5) {};
			\node (v8) at (-2.5,-2.5) {};
			\node (v9) at (2,-2) {};
		 	\node at (3.5,3.5) {$ X $};
		 	\node at (-3.5,-3.5) {$ X $};
		 	\node at (3,-3) {$ X^\vee $};
		 	\node at (-3,3) {$ X^\vee $};
			\draw[thick] (v6) to[out=-45,in=-135] (v7);
			\draw[thick] (v9) to[out=135,in=45] (v8);
		 	\draw[very thick, red]  (v1) edge (v3);
		 	\draw[very thick, red]  (v2) edge (v4);
		 	\draw[very thick]  (v1) edge (v2);
		 	\draw[very thick]  (v3) edge (v4);
		 	\end{tikzpicture}
		\end{array}.
	\end{align}
  	\begin{LEMMA} \label{lemma:commuting_twists}

	For all $ \beta \in \Hom_{\TC}(X,Y) $ we have,
		\begin{align*}
		\beta \circ t_X = t_Y \circ \beta.
		\end{align*}
	\proof

	Let $ \beta $ be in $ \Hom_{\TC}(X,Y) $. As any morphism in $ \TC $ may be written as a linear combination of elements of the form $ \alpha_G $ we may assume w.l.o.g.\ that $ \beta = \alpha_G $. We have
		\begin{align*}
		\alpha_G \circ t_{X} =
			\begin{array}{c}
				\begin{tikzpicture}[scale=0.25,every node/.style={inner sep=0,outer sep=-1}]
				\node (v1) at (0,5) {};
				\node (v2) at (5,0) {};
				\node (v3) at (-6,-1) {};
				\node (v4) at (-1,-6) {};
				\node (v6) at (-2,3) {};
				\node at (-5,2) {$ G $};
				\node at (2,-5) {$ G $};
				\node at (3.5,3.5) {$ X $};
				\node at (4,-3) {$ X^\vee $};
				\node at (-3,4) {$ X^\vee $};
				\node at (-4.5,-4.5) {$ Y $};
				\node (v5) at (-1,-1) {};
				\node (v10) at (3,-2) {};
				\draw [thick] (v5) to[out=45,in=135] (v10);
				\node (v14) at (1,-4) {};
				\node (v13) at (-4,1) {};
				\node (v7) at (2.5,2.5) {};
				\draw[thick]  (v6) to[out=-45,in=-135] (v7);
				\node (v8) at (-3.5,-3.5) {};
				\node [draw,diamond,outer sep=0,inner sep=.5,minimum size=20,fill=white] (v9) at (-1.5,-1.5) {$ \alpha $};
				\draw [thick] (v13) edge (v9);
				\draw [thick] (v9) edge (v14);
				\draw[very thick, red]  (v1) edge (v3);
				\draw[very thick, red]  (v2) edge (v4);
				\draw[very thick]  (v1) edge (v2);
				\draw[very thick]  (v3) edge (v4);
				\draw [thick] (v9) edge (v8);
				\end{tikzpicture}
			\end{array}
		=
			\begin{array}{c}
				\begin{tikzpicture}[scale=0.25,every node/.style={inner sep=0,outer sep=-1}]
				\node (v1) at (-2,3) {};
				\node (v2) at (3,-2) {};
				\node (v3) at (-8,-3) {};
				\node (v4) at (-3,-8) {};
				\node at (-5,2) {$ G $};
				\node at (2,-5) {$ G $};
				\node at (1.5,1.5) {$ X $};
				\node at (0,-7) {$ Y^\vee $};
				\node at (-7,0) {$ Y^\vee $};
				\node at (-6.5,-6.5) {$ Y $};
				\node (v5) at (-1,-1) {};
				\node (v14) at (-1,-6) {};
				\node (v13) at (-4,1) {};
				\node (v7) at (0.5,0.5) {};
				\node (v8) at (-5.5,-5.5) {};
				\draw [thick] (v8) to[out=45,in=135] (v14);
				\node (v6) at (1,-4) {};
				\node (v11) at (-2,-2) {};
				\node (v12) at (-6,-1) {};
				\draw [thick] (v12) to[out=-45,in=-135] (v11);
				\draw [thick] (v7) edge (v5);
				\node [draw,diamond,outer sep=0,inner sep=.5,minimum size=20,fill=white] (v9) at (-1.5,-1.5) {$ \alpha $};
				\draw [thick] (v9) edge (v13);
				\draw [thick] (v9) edge (v6);
				\draw[very thick, red]  (v1) edge (v3);
				\draw[very thick, red]  (v2) edge (v4);
				\draw[very thick]  (v1) edge (v2);
				\draw[very thick]  (v3) edge (v4);
				\end{tikzpicture}
			\end{array}
		= t_X \circ \alpha_G.
		\end{align*}
	as desired. \endproof

	\end{LEMMA}
As described in Section~\ref{sec:preliminaries}, if $ \C $ is modular then $ \e_I^J $ is a primitive idempotent. In particular we have $ \End_{\TC}(\e_I^J) = \fld $. However, by Lemma~\ref{lemma:commuting_twists}, we have
 	\begin{align*}
 	\e_I^J \circ t_{IJ} \circ \e_I^J  =\e_I^J \circ \e_I^J \circ t_{IJ} = \e_I^J \circ t_{IJ}
 	\end{align*}
so $ \e_I^J \circ t_{IJ} \in \End_{\TC}(\e_I^J) $. Therefore $ \e_I^J \circ t_{IJ} = \lambda \e_I^J  $ for some $ \lambda \in \fld $. This turns out to also be true in the case when $ \C $ is only assumed to be a pre-modular tensor category.

	\begin{PROP} \label{prop:twisting_the_idempotents}

	Let $ \C $ be a pre-modular tensor category and let $ I,J $ be in $ \I $. Then
		\begin{align*}
		\e_I^J \circ t_{IJ} = \frac{\Tc_{II}}{\Tc_{JJ}} \e_I^J.
		\end{align*}
	\proof

	We have

		\begin{align*}
		\e_I^J \circ t_{IJ} &= \bigoplus\limits_S  \sum\limits_{T,b} d(T)
			\begin{array}{c}
				\begin{tikzpicture}[scale=0.25,every node/.style={inner sep=0,outer sep=-1}]
				\node (v1) at (-2,7) {};
				\node (v2) at (7,-2) {};
				\node (v3) at (-8,1) {};
				\node (v4) at (1,-8) {};
				\node (v6) at (-2,3) {};
				\node (v7) at (2,3) {};
				\node at (5,-6) {$ S $};
				\node at (-6,5) {$ S $};
				\node at (-4.25,-0.5) {$ T $};
				\node at (4,3) {$ J $};
				\node at (3,4) {$ I $};
				\node at (-4,-5) {$ J $};
				\node at (-5,-4) {$ I $};
				\draw[thick] (v6) to[out=-45,in=-135] (v7);
				\node (v5) at (-4,-3) {};
				\node (v10) at (3,-2) {};
				\draw [thick] (v5) to[out=45,in=135] (v10);
				\node (v12) at (3,2) {};
				\node (v11) at (-3,2) {};
				\draw [thick] (v11) to[out=-45,in=-135] (v12);
				\node (v14) at (1,-4) {};
				\node (v13) at (-4,1) {};
				\draw [line width=0.15cm,white] (v13) edge (v14);
				\draw [thick] (v13) edge (v14);
				\node (v8) at (-3,-4) {};
				\node (v9) at (2,-3) {};
				\draw[line width=0.15cm,white] (v9) to[out=135,in=45] (v8);
				\draw[thick] (v9) to[out=135,in=45] (v8);
				\node [draw,rotate=45,outer sep=0,inner sep=2,minimum width=25,minimum height=11,fill=white] (v18) at (2.5,-3.5) {$ b^* $};
				\node [draw,rotate=45,outer sep=0,inner sep=2,minimum width=25,minimum height=11,fill=white] (v15) at (-3.5,2.5) {$ b $};
				\node (v16) at (-5,4) {};
				\node (v17) at (4,-5) {};
				\draw [thick] (v15) edge (v16);
				\draw [thick] (v17) edge (v18);
				\draw[very thick, red]  (v1) edge (v3);
				\draw[very thick, red]  (v2) edge (v4);
				\draw[very thick]  (v1) edge (v2);
				\draw[very thick]  (v3) edge (v4);
				\end{tikzpicture}
			\end{array}.
		\end{align*}
	Therefore the $ S $-summand of $ \e_I^J \circ t_{IJ} $ is given by
		\begin{align*}
		\sum\limits_{T,b} d(T)
			\begin{array}{c}
				\begin{tikzpicture}[scale=0.25,every node/.style={inner sep=0,outer sep=-1}]
				\node (v1) at (1.5,-3.5) {};
				\node (v5) at (4.5,-4) {};
				\node (v3) at (3,-4) {};
				\node (v4) at (7.5,-4) {};
				\node (v6) at (6,-4) {};
				\node (v2) at (1.5,-8) {};
				\draw [thick] (v5) to[out=-90,in=-90] (v6);
				\draw [line width=0.15cm,white] (v1) to[out=-90,in=-90] (v2);
				\draw [line width=0.15cm,white] (v3) to[out=-90,in=-90] (v4);
				\draw [thick] (v3) to[out=-90,in=-90] (v4);
				\node at (0.5,-4.75) { $ T $ };
				\node at (-1.5,-11) { $ I $ };
				\node at (7.5,-0.5) { $ J $ };
				\node at (0,-11) { $ J $ };
				\node at (6,-0.5) { $ I $ };
				\node (v7) at (6,-1.5) {};
				\node (v8) at (7.5,-1.5) {};
				\draw [thick] (v6) edge (v7);
				\draw [thick] (v4) edge (v8);
				\node (v9) at (3,-3.5) {};
				\node (v10) at (4.5,-3.5) {};
				\draw [thick] (v3) edge (v9);
				\draw [thick] (v5) edge (v10);
				\node (v11) at (3,-8) {};
				\node (v12) at (3,-7.5) {};
				\node (v13) at (4.5,-8) {};
				\node (v14) at (4.5,-7.5) {};
				\draw [thick] (v11) edge (v12);
				\draw [thick] (v13) edge (v14);
				\node (v15) at (0,-7.5) {};
				\node (v16) at (-1.5,-7.5) {};
				\node (v19) at (-1.5,-10) {};
				\draw [thick] (v14) to[out=90,in=90] (v16);
				\draw [thick] (v16) edge (v19);
				\draw [line width=0.15cm,white] (v1) edge (v2);
				\draw [thick] (v1) edge (v2);
				\draw [line width=0.15cm,white] (v12) to[out=90,in=90] (v15);
				\draw [thick] (v12) to[out=90,in=90] (v15);
				\node (v17) at (0,-10) {};
				\draw [thick] (v15) edge (v17);
				\node [draw,outer sep=0,inner sep=2,minimum width=25,minimum height=10,fill=white] (v18) at (3,-3) {$  b $};
				\node [draw,outer sep=0,inner sep=2,minimum width=25,minimum height=10,fill=white] (v180) at (3,-8.5) {$  b $};
				\node (v20) at (3,-1.5) {};
				\draw [thick] (v20) edge (v18);
				\node (v21) at (3,-10) {};
				\draw [thick] (v180) edge (v21);
				\node at (3,-0.5) {$ S $};
				\node at (3,-11) {$ S $};
				\end{tikzpicture}
			\end{array}
		&= \sum\limits_{T,b} d(T) \
			\begin{array}{c}
				\begin{tikzpicture}[scale=0.25,every node/.style={inner sep=0,outer sep=-1}]
				\node (v1) at (1.5,-3.5) {};
				\node (v5) at (4.5,-4) {};
				\node (v3) at (3,-4) {};
				\node (v4) at (7.5,-4) {};
				\node (v6) at (6,-4) {};
				\node (v2) at (1.5,-8) {};
				\draw [thick] (v5) to[out=-90,in=-90] (v6);
				\draw [line width=0.15cm,white] (v1) to[out=-90,in=-90] (v2);
				\draw [line width=0.15cm,white] (v3) to[out=-90,in=-90] (v4);
				\draw [thick] (v3) to[out=-90,in=-90] (v4);
				\node at (0.5,-6) { $ T $ };
				\node at (7.5,-0.5) { $ J $ };
				\node at (9,-14) { $ J $ };
				\node at (6,-0.5) { $ I $ };
				\node at (6,-14) { $ I $ };
				\node (v7) at (6,-1.5) {};
				\node (v8) at (7.5,-1.5) {};
				\draw [thick] (v6) edge (v7);
				\draw [thick] (v4) edge (v8);
				\node (v9) at (3,-3.5) {};
				\node (v10) at (4.5,-3.5) {};
				\draw [thick] (v3) edge (v9);
				\draw [thick] (v5) edge (v10);
				\node (v11) at (3,-8) {};
				\node (v12) at (3,-7.5) {};
				\node (v13) at (4.5,-8) {};
				\node (v14) at (4.5,-7.5) {};
				\draw [thick] (v11) edge (v12);
				\draw [thick] (v13) edge (v14);
				\node (v15) at (9,-7.5) {};
				\node (v150) at (10,-9.5) {};
				\node (v1500) at (9,-8) {};
				\node (v16) at (6,-7.5) {};
				\node (v19) at (6,-8) {};
				\draw [thick] (v14) to[out=90,in=90] (v16);
				\draw [thick] (v16) edge (v19);
				\draw [line width=0.15cm,white] (v1) edge (v2);
				\draw [thick] (v1) edge (v2);
				\draw [line width=0.15cm,white] (v12) to[out=90,in=90] (v15);
				\draw [thick] (v12) to[out=90,in=90] (v15);
				\node (v17) at (11,-9.5) {};
				\node [draw,outer sep=0,inner sep=2,minimum width=25,minimum height=10,fill=white] (v18) at (3,-3) {$  b $};
				\node [draw,outer sep=0,inner sep=2,minimum width=25,minimum height=10,fill=white] (v180) at (3,-8.5) {$  b $};
				\node (v20) at (11,-8) {};
				\node (v21) at (10,-8) {};
				\draw [thick] (v17) edge (v20);
				\node (v22) at (9,-9.5) {};
				\draw [thick] (v1500) to[out=-90,in=90] (v150);
				\draw [line width=0.15cm,white] (v21) to[out=-90,in=90] (v22);
				\draw [thick] (v21) to[out=-90,in=90] (v22);
				\draw [thick] (v20) to[out=90,in=90] (v21);
				\draw [thick] (v150) to[out=-90,in=-90] (v17);
				\node (v24) at (7,-8) {};
				\node (v23) at (7,-9.5) {};
				\node (v25) at (6,-9.5) {};
				\draw [thick] (v24) to[out=-90,in=90] (v25);
				\draw [line width=0.15cm,white] (v19) to[out=-90,in=90] (v23);
				\draw [thick] (v19) to[out=-90,in=90] (v23);
				\node (v26) at (8,-8) {};
				\node (v27) at (8,-9.5) {};
				\draw [thick] (v24) to[out=90,in=90] (v26);
				\draw [thick] (v23) to[out=-90,in=-90] (v27);
				\draw [thick] (v26) edge (v27);
				\node (v28) at (6,-13) {};
				\node (v29) at (9,-13) {};
				\draw [thick] (v25) edge (v28);
				\draw [thick] (v1500) edge (v15);
				\node (v20) at (3,-1.5) {};
				\draw [thick] (v20) edge (v18);
				\node (v21) at (4,-10.5) {};
				\draw [thick] (v180) to[out=-90,in=135] (v21);
				\node at (3,-0.5) {$ S $};
				\node at (10.5,-14) {$ S $};
				\node (v30) at (10.5,-13) {};
				\draw [line width=0.15cm,white] (v21) to[out=-45,in=90] (v30);
				\draw [thick] (v21) to[out=-45,in=90] (v30);
				\draw [line width=0.15cm,white] (v22) edge (v29);
				\draw [thick] (v22) edge (v29);
				\end{tikzpicture}
			\end{array}\\
		&= \frac{\Tc_{II}}{\Tc_{JJ}} \sum\limits_{T,b} d(T) \
			\begin{array}{c}
				\begin{tikzpicture}[scale=0.25,every node/.style={inner sep=0,outer sep=-1}]
				\node (v1) at (1.5,-3.5) {};
				\node (v5) at (4.5,-4) {};
				\node (v3) at (3,-4) {};
				\node (v4) at (7.5,-4) {};
				\node (v6) at (6,-4) {};
				\node (v2) at (1.5,-8) {};
				\draw [thick] (v5) to[out=-90,in=-90] (v6);
				\draw [line width=0.15cm,white] (v1) to[out=-90,in=-90] (v2);
				\draw [line width=0.15cm,white] (v3) to[out=-90,in=-90] (v4);
				\draw [thick] (v3) to[out=-90,in=-90] (v4);
				\node at (2.5,-6) { $ T $ };
				\node at (7.5,-0.5) { $ J $ };
				\node at (4,-14) { $ J $ };
				\node at (6,-0.5) { $ I $ };
				\node at (2.5,-14) { $ I $ };
				\node (v7) at (6,-1.5) {};
				\node (v8) at (7.5,-1.5) {};
				\draw [thick] (v6) edge (v7);
				\draw [thick] (v4) edge (v8);
				\node (v9) at (3,-3.5) {};
				\node (v10) at (4.5,-3.5) {};
				\draw [thick] (v3) edge (v9);
				\draw [thick] (v5) edge (v10);
				\node (v11) at (3,-8) {};
				\node (v12) at (3,-7.5) {};
				\node (v13) at (4.5,-8) {};
				\node (v14) at (4.5,-7.5) {};
				\draw [thick] (v11) edge (v12);
				\draw [thick] (v13) edge (v14);
				\draw [thick] (v1) edge (v2);
				\node (v17) at (6,-9.5) {};
				\node (v19) at (7.5,-9.5) {};
				\node (v15) at (6,-7.5) {};
				\node (v16) at (7.5,-7.5) {};
				\draw [thick] (v14) to[out=90,in=90] (v15);
				\draw [thick] (v12) to[out=90,in=90] (v16);
				\draw [thick] (v15) edge (v17);
				\draw [thick] (v16) edge (v19);
				\node [draw,outer sep=0,inner sep=2,minimum width=25,minimum height=10,fill=white] (v18) at (3,-3) {$  b $};
				\node [draw,outer sep=0,inner sep=2,minimum width=25,minimum height=10,fill=white] (v180) at (3,-8.5) {$  b $};
				\node (v20) at (3,-1.5) {};
				\draw [thick] (v20) edge (v18);
				\node (v21) at (3,-9.5) {};
				\draw [thick] (v180) edge (v21);
				\node (v22) at (2.5,-13) {};
				\node (v23) at (4,-13) {};
				\node (v24) at (7,-13) {};
				\draw [thick] (v17) to[out=-90,in=90] (v22);
				\draw [line width=0.15cm,white] (v21) to[out=-90,in=90] (v24);
				\draw [thick] (v21) to[out=-90,in=90] (v24);
				\draw [line width=0.15cm,white] (v19) to[out=-90,in=90] (v23);
				\draw [thick] (v19) to[out=-90,in=90] (v23);
				\node at (7,-14) {$ S $};
				\node at (3,-0.5) {$ S $};
				\end{tikzpicture}
			\end{array}
        \end{align*}
        \begin{align*}
		&= \frac{\Tc_{II}}{\Tc_{JJ}}  d(S)
			\begin{array}{c}
				\begin{tikzpicture}[scale=0.25,every node/.style={inner sep=0,outer sep=-1}]
				\node at (5,-14) { $ J $ };
				\node at (7,-8) { $ J $ };
				\node at (3,-14) { $ I $ };
				\node at (5,-8) { $ I $ };
				\node at (7,-14) {$ S $};
				\node at (3,-8) {$ S $};
				\node (v17) at (5,-9) {};
				\node (v19) at (7,-9) {};
				\node (v21) at (3,-9) {};
				\node (v22) at (3,-13) {};
				\node (v23) at (5,-13) {};
				\node (v24) at (7,-13) {};
				\draw [thick] (v17) to[out=-90,in=90] (v22);
				\draw [line width=0.15cm,white] (v21) to[out=-90,in=90] (v24);
				\draw [thick] (v21) to[out=-90,in=90] (v24);
				\draw [line width=0.15cm,white] (v19) to[out=-90,in=90] (v23);
				\draw [thick] (v19) to[out=-90,in=90] (v23);
				\end{tikzpicture}
			\end{array}
		\end{align*}
	where the final equality is due to an application of Lemma~3.11 in~\cite{hardiman_king} for $ X = IJ $. As this is exactly the $ S $-summand of $ \frac{\Tc_{II}}{\Tc_{JJ}} \e_I^J $ we are done. \endproof

	\end{PROP}

We may now prove the main result of this section.

	\begin{THM} \label{thm:t_invariance}

	Let $ \C $ be a modular tensor category and let $ F $ be an object in $ \RTC $. $ F $ is T-invariant if and only if $ F(t_X) = \id_{F(X)} $ for all $ X $ in $ \C $.

	\proof

	As $ \C $ is modular the $ (IJ,\e_I^J)\Sharp $ form a complete set of simple objects in $ \TC $. We can therefore decompose $ F $ as
		\begin{align*}
		F = \bigoplus\limits_{IJ} F_I^J \vprod (IJ,\e_I^J)\Sharp.
		\end{align*}
	Evaluating this on $ t_X $ gives
		\begin{align*}
		F(t_X) =  \bigoplus\limits_{IJ} \id_{F_I^J} \otimes (IJ,\e_I^J)\Sharp(t_X).
		\end{align*}
	By Lemma~\ref{lemma:commuting_twists} and Proposition~\ref{prop:twisting_the_idempotents} we have, for $ \alpha \in \Hom_{\TC}(X,\e_I^J) $,
		\begin{align*}
		(IJ,\e_I^J)\Sharp(t_X)(\alpha) = \e_I^J \circ \alpha \circ t_X = \e_I^J \circ t_{IJ} \circ \alpha = \frac{\Tc_{II}}{\Tc_{JJ}} \e_I^J \circ \alpha = \frac{\Tc_{II}}{\Tc_{JJ}} \alpha.
		\end{align*}
	Therefore
		\begin{align*}
		\bigoplus\limits_{IJ} \id_{F_I^J} \otimes (IJ,\e_I^J)\Sharp(t_X) &= \bigoplus\limits_{IJ} \frac{\Tc_{II}}{\Tc_{JJ}} \id_{F_I^J \otimes (IJ,\e_I^J)\Sharp(X)}.
		\end{align*}
	This is equal to $ \id_{F(X)} $ if and only if $ F_I^J \neq 0 $ implies $ \frac{\Tc_{II}}{\Tc_{JJ}} = 1 $. As $ \Tc $ is diagonal that is precisely the condition that $ Z(F) $ commutes with $ \Tc $. \endproof

	\end{THM}

	\begin{COR} \label{cor:t-inv}

	Let $ \C $ be a modular tensor category and let $ \M \colon \C \to \D $ be a pivotal monoidal functor. Then $ \TM $ is T-invariant.

	\proof

	For $ \alpha \in \TM(X) = \Hom_{\D}(\M(X),\tid) $ we have
		\begin{align*}
		\TM(t_X) \colon
			\begin{array}{c}
				\begin{tikzpicture}[scale=0.15,every node/.style={inner sep=0,outer sep=-1},yscale=-1]
				\node (v1) at (0,-1) {};
				\node (v2) at (0,-3) {};
				\draw [thick,blue] (v1) edge (v2);
				\node [blue] at (0,-4.5) {$ X $};
				\node at (0,4.8) {};
				\node [draw,outer sep=0,inner sep=1,minimum width=15,minimum height=11,fill=white] at (0,0) {$ \alpha $};
				\end{tikzpicture}
			\end{array}
		&\mapsto
			\begin{array}{c}
				\begin{tikzpicture}[scale=0.15,every node/.style={inner sep=0,outer sep=-1},yscale=-1]
				\node (v1) at (0,-1) {};
				\node (v2) at (0,-2) {};
				\draw [thick,blue] (v1) edge (v2);
				\node [blue] at (-3,-6) {$ X $};
				\node [draw,outer sep=0,inner sep=1,minimum width=15,minimum height=11,fill=white] at (0,0) {$ \alpha $};
				\node (v3) at (3,-2) {};
				\node (v4) at (3,1) {};
				\node (v5) at (-3,1) {};
				\node (v6) at (-3,-4) {};
				\draw [thick, blue] (v2) to[out=-90,in=-90] (v3);
				\draw [thick, blue] (v3) edge (v4);
				\draw [thick, blue] (v4) to[out=90,in=90] (v5);
				\draw [thick, blue] (v5) edge (v6);
				\end{tikzpicture}
			\end{array}
		=
			\begin{array}{c}
				\begin{tikzpicture}[scale=0.15,every node/.style={inner sep=0,outer sep=-1},yscale=-1]
				\node (v1) at (0,-1) {};
				\node (v2) at (0,-3) {};
				\draw [thick,blue] (v1) edge (v2);
				\node [blue] at (0,-4.5) {$ X $};
				\node at (0,4.8) {};
				\node [draw,outer sep=0,inner sep=1,minimum width=15,minimum height=11,fill=white] at (0,0) {$ \alpha $};
				\end{tikzpicture}
			\end{array}
		\end{align*}
	as $ \M $ is pivotal. Therefore Theorem~\ref{thm:t_invariance} applies, and $ \TM $ is T-invariant.	\endproof

	\end{COR}

\section{S-Invariance and Frobenius Algebras} \label{sec:s_inv}

	\begin{DEF}

	Let $ \C $ be a pre-modular tensor category and let $ \Sc $ be the S-matrix of $ \C $ as defined by~\eqref{eq:S_and_T_matrices}. We call an object $ F $ in $ \RTC $ \emph{S-invariant} if $ Z(F) $ commutes with $ \Sc $.

	\end{DEF}

We start with an example which illustrates that, even when $ \C $ is modular, $ \TM $ is \emph{not} necessarily S-invariant.

	\begin{EX} \label{ex:tm_not_s_inv}

	Let $ \C $ be an modular tensor category and let $ \M $ be the identity functor on $ \C $. Then
		\begin{align*}
	      \TM_I^J \leq \Hom_{\C}(IJ,\tid) = \delta_{I,J^\vee} \left\langle
            \begin{array}{c}
                \begin{tikzpicture}[scale=0.25,every node/.style={inner sep=0,outer sep=-1},yscale=-1]
                \node (v1) at (-1,0) {};
                \node (v2) at (1,0) {};
                \draw [thick] (v1) to[out=90,in=90] (v2);
                \node at (-1,-1.04) {$ J $};
                \node at (1,-1) {$ J^\vee $};
                \end{tikzpicture}
            \end{array}
         \right\rangle_{\fld}.
		\end{align*}
	We now suppose $ \delta_{I,J^\vee}  = 1 $. By Proposition~\ref{prop:characterising_tm_xy}, $ \TM_{J^\vee}^J $ is non-trivial if and only if
		\begin{align*}
			\begin{array}{c}
        		\begin{tikzpicture}[scale=0.15,every node/.style={inner sep=0,outer sep=-1},yscale=-1]
				\node (v1) at (-0.5,-0.5) {};
				\node (v11) at (2,3.5) {$ Z $};
				\node (v10) at (-6,-8) {$ Z $};
				\node (v11) at (-2,-8) {$ J^\vee $};
				\node (v10) at (2,-8) {$ J $};
				\node (v9) at (-2,-6) {};
				\node (v13) at (-4.5,-0.5) {};
				\node (v2) at (2,1.5) {};
				\node (v3) at (2,-0.5) {};
				\node (v5) at (2,-6) {};
				\node (v4) at (-6,-6) {};
				\draw [thick] (v2) edge (v3);
				\draw [thick] (v9) to[out=90, in=-90] (v13);
				\draw [line width= 0.15cm,white] (v3) to[out=-90, in=90] (v4);
				\draw [thick] (v3) to[out=-90, in=90] (v4);
				\draw [line width= 0.15cm,white] (v1) to[out=-90, in=90] (v5);
				\draw [thick] (v1) to[out=-90, in=90] (v5);
				\draw [thick] (v13) to[out=90, in=90] (v1);
				\end{tikzpicture}
			\end{array}
		=
			\begin{array}{c}
				\begin{tikzpicture}[scale=0.15,every node/.style={inner sep=0,outer sep=-1},yscale=-1]
				\node (v1) at (2,-0.5) {};
				\node (v11) at (-6,3.5) {$ Z $};
				\node (v10) at (-6,-8) {$ Z $};
				\node (v11) at (-2,-8) {$ J^\vee $};
				\node (v10) at (2,-8) {$ J $};
				\node (v9) at (-2,-6) {};
				\node (v13) at (-2,-0.5) {};
				\node (v2) at (-6,2) {};
				\node (v3) at (-6,1) {};
				\node (v5) at (2,-6) {};
				\node (v4) at (-6,-6) {};
				\draw [thick] (v2) edge (v3);
				\draw [thick] (v9) to[out=90, in=-90] (v13);
				\draw [line width= 0.15cm,white] (v3) to[out=-90, in=90] (v4);
				\draw [thick] (v3) to[out=-90, in=90] (v4);
				\draw [line width= 0.15cm,white] (v1) to[out=-90, in=90] (v5);
				\draw [thick] (v1) to[out=-90, in=90] (v5);
				\draw [thick] (v13) to[out=90, in=90] (v1);
				\end{tikzpicture}
			\end{array}
		\end{align*}
	for all $ Z $ in $ \C $. Post-composing this equality with $ \id_Z \otimes \ann_{J} $ and taking the trace implies $ \Sc_{ZJ} = d(Z) d(J) $ for all $ Z $ in $ \C $. Therefore the $ J $-th column in $ \Sc $ is proportional to the $ \tid $-th column. As $ \C $ is modular this implies $ J = \tid $. In summary, we have
		\begin{align*}
		Z(\TM)_{IJ} =
			\begin{cases}
			1 \quad \text{if $ I = J = \tid $} \\
			0 \quad \text{else}.
			\end{cases}
		\end{align*}
	Conjugating this matrix with $ \Sc $ and using the fact that $ \Sc_{I,J} = d(\C) \Sc_{I,J^\vee}^{-1} $ gives us
	 	\begin{align*}
	 	\left( \Sc \ Z(\TM) \ \Sc^{-1} \right)_{\tid \tid} = \frac{1}{d(\C)}
	 	\end{align*}
	implying that $ \Sc $-invariance will fail whenever $ d(\C) \neq 1 $.

	\end{EX}

However, when $ \C $ is a modular tensor category, we have a helpful theorem from~\cite{MR2551797}. For the prerequisite definitions on Frobenius algebras see Definition~\ref{def:frob_alg}.

	\begin{THM}[Theorem 3.4,~\cite{MR2551797}] \label{thm:kong_runkel}

	Let $ A $ be a haploid, symmetric, commutative Frobenius algebra in $ \CC $. Then the $ \I \times \I $-matrix with entries $ \hom(I \boxtimes J, A) $ commutes with the S-matrix of $ \C $ (where, as before, $ \hom $ denotes the dimension of the relevant $ \Hom $-space) if and only if
		\begin{align} \label{eq:s11_condition}
		d(A) = d(\C).
		\end{align}
	\end{THM}

	\begin{REM}

	We note that
		\begin{align*}
		d(A) &= \sum\limits_{IJ} \hom_{\CC}(I \boxtimes J, A) d(I \boxtimes J)\\
		&= \sum\limits_{IJ} \hom_{\CC}(I \boxtimes J, A) d(I)d(J) \\
		&= \sum\limits_{IJ} \Sc_{\tid,I} \hom_{\CC}(I \boxtimes J, A) \Sc_{J,\tid}.
		\end{align*}
	As $ \Sc_{I,J} = d(\C) \Sc_{I,J^\vee}^{-1} $, Condition~\eqref{eq:s11_condition} is precisely the condition that the matrix with entries $ \hom(I \boxtimes J, A) $ commutes with the S-matrix evaluated at $ (\tid,\tid) $ for an arbitrary object $ A $ in $ \CC $. Condition~\eqref{eq:s11_condition} is therefore certainly necessary, the content of the theorem is that, when $ A $ is a haploid, symmetric, commutative Frobenius algebra, it is also sufficient.

	\end{REM}

	\begin{REM} \label{rem:mod_inv_algebra}

	\cite[Theorem 3.4]{MR2551797} actually proves that when $ A $ is a haploid, symmetric, commutative Frobenius algebra~\eqref{eq:s11_condition} implies an equality which is strictly stronger than the result stated here. In particular, \cite[Theorem 3.4]{MR2551797} proves that $ A $ will be a \emph{modular invariant algebra}. This notion is defined and motivated in~\cite[Section 6]{MR2430629}.

	\end{REM}

	\begin{REM}

	As explained in the proof of~\cite[Theorem 3.4]{MR2551797}, there exists an MTC (the category of of local $ A $-modules) whose dimension is given by $ \frac{d(\C)^2}{d(A)^2} $. Combining this with the fact that any MTC over the complex numbers has dimension at least 1~\cite[Theorem 2.3.]{MR2183279} tells us that, in the case when $ \fld = \mathbb{C} $, the dimension of $ A $ cannot exceed $ d(\C) $.

	\end{REM}

We recall that, when $ \C $ is modular, $ \Phi\colon \CC \to \RTC $ is an equivalence and $ I \boxtimes J \mapsto (IJ,\e_I^J)\Sharp $. Therefore, for $ F $ in $ \RTC $, the $ \I \times \I $ matrix with entries $ \hom_{\RTC}(\Phi(I \boxtimes J), F) $ is precisely $ Z(F) $. The goal of the following section is to prove that $ \TM $ is a commutative algebra in $ \RTC $, and then, under a further condition on $ \M $, to show that it is also a haploid, symmetric, commutative Frobenius algebra.

The above stated goal assumes that $ \RTC $ is a braided pivotal monoidal category; this is indeed the case as $ \RTC = Z(\C) $ (see Remark~\ref{rem:equiv_zc}) and $ Z(\C) $ admits a canonical MTC structure~\cite[Corollary 8.20.13]{Etingof15}. To achieve this goal we are therefore going to have to work with the monoidal product, braiding and pivotal structure that $ \RTC $ inherits from $ Z(\C) $. In general this is not easy; for instance it is hard to express the tensor product of two generic objects in $ \RTC $. However, if we restrict our attention to functors coming from idempotents of the form~\eqref{eq:idempotents} these structures may be described graphically.

 	\begin{DEF} \label{def:an_associative_product}

	Let $ \otimes_{\TC} \colon \TC \times \TC \to \TC  $ be the bifunctor given by $$ X \otimes_{\TC} Y = XY $$ for $ X, Y $ in $ \TC $ and
	 	\begin{align} \label{eq:an_associative_product}
	 	f \otimes_{\TC} g = d(\C) \bigoplus\limits_{S} \frac{1}{d(S)}
			\begin{array}{c}
				\begin{tikzpicture}[scale=0.5,every node/.style={inner sep=0,outer sep=-1}]
				\node [draw,diamond,outer sep=0,inner sep=.5,minimum size=24] (v10) at (-1.5,1.5) {$ f_S $};
				\node [draw,diamond,outer sep=0,inner sep=.5,minimum size=24] (v20) at (0,0) {$ g_S $};
				\node (v6) at (-2.5,2.5) {};
				\node (v11) at (1,-1) {};
				\node (v8) at (1,1) {};
				\node (v9) at (-1,-1) {};
				\node (v7) at (-2.5,0.5) {};
				\node (v5) at (-0.5,2.5) {};
				\draw [thick] (v6) edge (v10);
				\draw [thick] (v5) edge (v10);
				\draw [thick] (v10) edge (v7);
				\draw [thick] (v10) edge (v20);
				\draw [thick] (v20) edge (v8);
				\draw [thick] (v20) edge (v9);
				\draw [thick] (v20) edge (v11);
				\node at (-3,3) {$ S $};
				\node at (1.5,-1.5) {$ S $};
				\node at (0,3) {$ W $};
				\node at (-3,0) {$ X $};
				\node at (1.5,1.5) {$ Y $};
				\node at (-1.5,-1.5) {$ Z $};
				\node (v1) at (-1.5,3.5) {};
				\node (v4) at (0,-2) {};
				\node (v2) at (2,0) {};
				\node (v3) at (-3.5,1.5) {};
				\draw[very thick, red]  (v1) edge (v3);
				\draw[very thick, red]  (v2) edge (v4);
				\draw[very thick]  (v1) edge (v2);
				\draw[very thick]  (v3) edge (v4);
				\end{tikzpicture}
			\end{array} \in \Hom_{\TC}(WY,XZ)
	 	\end{align}
	for $ f \in \Hom_{\TC}(W,X) $ and $ g \in \Hom_{\TC}(Y,Z) $.

	\end{DEF}

We note that this product does not give a monoidal product as there is no unit. The tensor identity $ \tid $ in $ \C $ fails to give a unit as the functor
	\begin{align*}
	\Blank \otimes_{\TC} \tid \colon \TC \to \TC
	\end{align*}
maps $ \alpha \in \Hom_{\TC}(X,Y) $ to $ d(\C) \alpha_{\tid} \in \Hom_{\C}(X,Y) $ and so the unit isomorphisms fail to be natural.

	\begin{REM} \label{rem:iso_idem_via_braiding}

	The scalars appearing in~\eqref{eq:an_associative_product} are chosen to guarantee that $ \otimes_{\TC} $ is well-behaved with respect to idempotents of the form $ \e_X^Y $. Indeed, we have
		\begin{align*}
		(XYAB,\e_{X}^{Y}\otimes_{\TC}\e_{A}^{B})\Sharp = (XAYB,\e_{XA}^{YB})\Sharp
		\end{align*}
	where the isomorphism is once again given by the braiding.

	\end{REM}

   \begin{PROP} \label{prop:graphical_tensor}

    Let $ \C $ be a modular tensor category. For $ X,Y $ objects in $ \C $, we have
        \begin{align*}
        (XY,\e_X^Y)\Sharp \otimes (AB,\e_A^B)\Sharp = (XAYB,\e_{XA}^{YB})\Sharp = (XYAB,\e_{X}^{Y}\otimes_{\TC}\e_{A}^{B})\Sharp.
        \end{align*}
    Furthermore, for $ \alpha \in \Hom_{\TC}(\e_X^Y,\e_A^B) $ and $ \beta \in \Hom_{\TC}(\e_{X'}^{Y'},\e_{A'}^{B'}) $, we have
        \begin{align*}
        \alpha \otimes \beta = \alpha \otimes_{\TC} \beta
        \end{align*}
    where $ \otimes_{\TC} $ is the associative product given by Definition~\ref{def:an_associative_product}.

    \proof

	By \eqref{eq:phi_and_eXY}, we have
		\begin{align*}
		(XY,\e_X^Y)\Sharp \otimes (AB,\e_A^B)\Sharp &= \Phi(X \boxtimes Y) \otimes \Phi(A \boxtimes B) \\
		&= \Phi(XA \boxtimes YB) \\
		&= (XAYB,\e_{XA}^{YB})\Sharp
		\end{align*}
	where $ \Phi $ is as defined in Section~\ref{sec:preliminaries}. As described in Remark~\ref{rem:iso_idem_via_braiding}, we then have
        \begin{align*}
        (XAYB,\e_{XA}^{YB})\Sharp = (XYAB,\e_{X}^{Y}\otimes_{\TC}\e_{A}^{B})\Sharp.
        \end{align*}
	where the natural isomorphism is given by the braiding. This proves the first half of the proposition.

	Let $ f,f',g,g' $ be in $ \Hom_{\C}(X,A) $, $\Hom_{\C}(X',A')$, $\Hom_{\C}(Y,B)$ and $\Hom_{\C}(Y',B') $ respectively and let $ \alpha $ and $ \beta $ be given by
		\begin{align} \label{eq:pullback_phi}
		\alpha = \Phi(f \boxtimes g) \quad \text{and} \quad \beta = \Phi(f' \boxtimes g').
		\end{align}
	Then
		\begin{align*}
		\alpha \otimes \beta &= \Phi((f \boxtimes g) \otimes (f \boxtimes g)) \\
		&= \Phi( (f \otimes f') \boxtimes (g \otimes g') )\\
		&= \alpha \otimes_{\TC} \beta.
		\end{align*}
	as desired.

	As $ \Phi $ is fully faithful any morphism in $ \TC $ may be written as a sum of morphisms of the form~\eqref{eq:pullback_phi}. This implies $ \alpha \otimes \beta = \alpha \otimes_{\TC} \beta $ for arbitrary $ \alpha $ and $ \beta $. \endproof

    \end{PROP}

The braiding between $ (XY,\e_X^Y)\Sharp $ and $ (AB,\e_A^B)\Sharp $ is then given by the following morphism,
	\begin{align*}
	\sigma_{X,A}^{Y,B} &= \frac{1}{d(\C)} \bigoplus d(T)
		\begin{array}{c}
			\begin{tikzpicture}[scale=0.4,every node/.style={inner sep=0,outer sep=-1}]
			\node (v1) at (-0.5,4.5) {};
			\node (v4) at (-0.5,-3.5) {};
			\node (v2) at (3.5,0.5) {};
			\node (v3) at (-4.5,0.5) {};
			\node (v28) at (-3.75,1.25) {};
			\node at (-4.25,1.75) {$ T $};
			\node at (0.5,4.5) {$ X $};
			\node at (2.5,2.5) {$ A $};
			\node at (1.5,3.5) {$ Y $};
			\node at (3.5,1.5) {$ B $};
			\node at (0.75,-3.25) {$ T $};
			\node (v5) at (0,4) {};
			\node (v11) at (1,3) {};
			\node (v7) at (2,2) {};
			\node (v9) at (3,1) {};
			\node (v8) at (-4,0) {};
			\node (v10) at (-3,-1) {};
			\node (v12) at (-1,-3) {};
			\node (v6) at (-2,-2) {};
			\node (v13) at (0.25,-2.75) {};
			\draw [thick] (v7) to[out=-135,in=45] (v8);
			\draw [line width= 0.2cm,white] (v5) to[out=-135,in=45] (v6);
			\draw [thick] (v5) to[out=-135,in=45] (v6);
			\draw [line width= 0.2cm,white] (v28) edge (v13);
			\draw [thick] (v28) edge (v13);
			\draw [line width= 0.2cm,white] (v11) to[out=-135,in=45] (v12);
			\draw [thick] (v11) to[out=-135,in=45] (v12);
			\draw [line width= 0.2cm,white] (v9) to[out=-135,in=45] (v10);
			\draw [thick] (v9) to[out=-135,in=45] (v10);
			\draw[very thick, red]  (v1) edge (v3);
			\draw[very thick, red]  (v2) edge (v4);
			\draw[very thick]  (v1) edge (v2);
			\draw[very thick]  (v3) edge (v4);
			\end{tikzpicture}
		\end{array}\\
	\in  &\Hom_{\TC}(\e_{X}^{Y}\otimes_{\TC}\e_{A}^{B}, \e_{A}^{B}\otimes_{\TC}\e_{X}^{Y})\\
	= &\Hom_{\RTC}((XYAB,\e_{X}^{Y}\otimes_{\TC}\e_{A}^{B})\Sharp , (ABXY,\e_{A}^{B}\otimes_{\TC}\e_{X}^{Y})\Sharp)
	\end{align*}
and the creation and annihilation morphisms for $ (XY,\e_X^Y)\Sharp $ and $ (Y^\vee X^\vee,\e_{X^\vee}^{Y^\vee})\Sharp $ are given by
	\begin{align*}
	\frac{1}{d(\C)} \bigoplus d(T)
		\begin{array}{c}
			\begin{tikzpicture}[scale=0.4,every node/.style={inner sep=0,outer sep=-1},xscale=-1,yscale=-1]
			\node (v1) at (-0.5,4.5) {};
			\node (v4) at (1.5,-1.5) {};
			\node (v2) at (3.5,0.5) {};
			\node (v3) at (-2.5,2.5) {};
			\node (v28) at (-1.75,3.25) {};
			\node at (-2.25,3.75) {$ T $};
			\node at (0.625,4.5) {$ X^\vee $};
			\node at (2.625,2.625) {$ Y $};
			\node at (1.625,3.5) {$ Y^\vee $};
			\node at (3.625,1.625) {$ X $};
			\node at (2.75,-1.25) {$ T $};
			\node (v5) at (0,4) {};
			\node (v11) at (1,3) {};
			\node (v7) at (2,2) {};
			\node (v9) at (3,1) {};
			\node (v8) at (-2,2) {};
			\node (v10) at (-1,1) {};
			\node (v12) at (1,-1) {};
			\node (v6) at (0,0) {};
			\node (v13) at (2.25,-0.75) {};
			\draw [thick] (v11) to[out=-135,in=-135] (v7);
			\draw [thick] (v9) to[out=-135,in=-135] (v5);
			\draw [thick] (v28) edge (v13);
			\draw[very thick, red]  (v1) edge (v3);
			\draw[very thick, red]  (v2) edge (v4);
			\draw[very thick]  (v1) edge (v2);
			\draw[very thick]  (v3) edge (v4);
			\end{tikzpicture}
		\end{array}
	\quad \text{and}\quad  \frac{1}{d(\C)} \bigoplus d(T)
		\begin{array}{c}
			\begin{tikzpicture}[scale=0.4,every node/.style={inner sep=0,outer sep=-1}]
			\node (v1) at (-0.5,4.5) {};
			\node (v4) at (1.5,-1.5) {};
			\node (v2) at (3.5,0.5) {};
			\node (v3) at (-2.5,2.5) {};
			\node (v28) at (-1.75,3.25) {};
			\node at (-2.25,3.75) {$ T $};
			\node at (0.5,4.5) {$ X $};
			\node at (2.75,2.625) {$ Y^\vee $};
			\node at (1.5,3.5) {$ Y $};
			\node at (3.75,1.625) {$ X^\vee $};
			\node at (2.75,-1.25) {$ T $};
			\node (v5) at (0,4) {};
			\node (v11) at (1,3) {};
			\node (v7) at (2,2) {};
			\node (v9) at (3,1) {};
			\node (v8) at (-2,2) {};
			\node (v10) at (-1,1) {};
			\node (v12) at (1,-1) {};
			\node (v6) at (0,0) {};
			\node (v13) at (2.25,-0.75) {};
			\draw [thick] (v11) to[out=-135,in=-135] (v7);
			\draw [thick] (v9) to[out=-135,in=-135] (v5);
			\draw [thick] (v28) edge (v13);
			\draw[very thick, red]  (v1) edge (v3);
			\draw[very thick, red]  (v2) edge (v4);
			\draw[very thick]  (v1) edge (v2);
			\draw[very thick]  (v3) edge (v4);
			\end{tikzpicture}
		\end{array}
	\end{align*}
respectively. Note the tensor identity in $ \RTC $ is $ (\tid, \e_{\tid}^{\tid})\Sharp $ and \emph{not} $ \Hom_{\TC}(\Blank, \tid) $.

As $ \TM $ is not of the form $ (XY,\e_X^Y)\Sharp $ equipping it with the structure of a Frobenius algebra directly is difficult. However, as we are assuming that $ \C $ is modular then we can decompose $ \TM $ as follows:
	\begin{align*}
	\TM = \bigoplus_{I,J} \TM_I^J \vprod \e_I^J.
	\end{align*}
We may then define the Frobenius structure in terms of this decomposition. This is the approach adopted by the following section.

\section{$ \TM $ as a Frobenius Algebra} \label{sec:tm_as_a_frob_alg}

As before let $ \C $ be an MTC and let $ \M $ be a pivotal monoidal functor from $ \C $ to $ \D $ where $ \D $ is a pivotal monoidal category. Before preceding with the description of a Frobenius algebra structure on $ \TM $ we give a brief outline of our strategy. Our first step is to equip $ \TM $ with the structure of an algebra. We do this by specifying a map
   \begin{align*}
   \nabla_{\bX}^{\bY,\bZ} \colon \Hom_{\RTC}(\bX,\bY\bZ) &\to \Hom(\TM_{\bY} \otimes \TM_{\bZ}, \TM_{\bX})
   \end{align*}
for all $ \bX,\bY,\bZ $ in $ \RTC $ of the form $ (AB,\e_A^B)\Sharp $.  As $ \C $ is modular $ \{ (IJ,\e_I^J)\Sharp \}_{I,J \in \I} $ forms a complete set of simples and $ \nabla_{\bX}^{\bY,\bZ} $ determines a map $ \nabla \colon \TM \otimes \TM \to \TM $ as described in Remark~\ref{rem:decompose_morphisms}. We then verify that this does indeed give an commutative algebra structure via Lemma~\ref{lem:build_an_algebra}.

An important property of Frobenius algebras is that they naturally carry a self-dual structure. Indeed, it is simple to check that, for a Frobenius algebra $ A $, the maps
    \begin{align} \label{eq:frobenius_dual_maps}
        \begin{array}{c}
            \begin{tikzpicture}[scale=0.3, yscale = -1]
            \coordinate (v1) at (-0.5,0.3) {};
            \coordinate (v2) at (0.5,0.3) {};
            \coordinate (v3) at (0,-0.45) {};
            \draw[thick] (-1,2) to[out=-90] ($ (v1) + (0.11,-0.1) $);
            \draw[thick] (1,2) to[out=-90,in=45] ($ (v2) + (-0.11,-0.1) $);
            \draw[thick] (0,0) -- (0,-2);
            \fill [blue] (v1) -- (v2) -- (v3) -- cycle;
			\draw[thick, fill=white] (0,-2) circle (0.2);
            \node at (-1,2.8) {$A$};
            \node at (1,2.8) {$A$};
            \end{tikzpicture}
        \end{array}
    \quad \text{and} \quad
        \begin{array}{c}
            \begin{tikzpicture}[scale=0.3]
            \coordinate (v1) at (-0.5,0.3) {};
            \coordinate (v2) at (0.5,0.3) {};
            \coordinate (v3) at (0,-0.45) {};
            \draw[thick] (-1,2) to[out=-90] ($ (v1) + (0.11,-0.1) $);
            \draw[thick] (1,2) to[out=-90,in=45] ($ (v2) + (-0.11,-0.1) $);
            \draw[thick] (0,0) -- (0,-2);
            \fill [blue] (v1) -- (v2) -- (v3) -- cycle;
			\draw[thick, fill=white] (0,-2) circle (0.2);
            \node at (-1,2.8) {$A$};
            \node at (1,2.8) {$A$};
            \end{tikzpicture}
        \end{array}
    \end{align}
are self-dualizing structure maps on $ A $ (where we are using the notation of Definition~\ref{def:frob_alg}). We therefore proceed by identifying self-dualizing structure maps on $ \TM $. These self-dualizing structure maps may also be given in terms of simple multiplicity spaces. Under certain additional conditions on $ \M $, we describe a collection of perfect pairings $ \TM_X^Y \otimes \TM_{X^\vee}^{Y^\vee} \to \fld $ which may then be used to construct self-dualizing structure maps via Lemma~\ref{lem:dualizing_maps}.

Once the self-dualizing structure maps are established we note that if $ \TM $ \emph{were} a Frobenius algebra the Frobenius condition~\eqref{eq:frob_condition} tells us that the coproduct could be written as both sides of the following condition
    \begin{align} \label{eq:unique_coproduct}
        \begin{array}{c}
            \begin{tikzpicture}[scale=0.3]
            \coordinate (v1) at (-0.5,0.3) {};
            \coordinate (v2) at (0.5,0.3) {};
            \coordinate (v3) at (0,-0.45) {};
            \draw[thick] (-1,2) to[out=-90] ($ (v1) + (0.11,-0.1) $);
            \draw[thick] (1,1) to[out=-90,in=45] ($ (v2) + (-0.11,-0.1) $);
            \draw[thick] (1,1) to[out=90,in=90] (3,1);
            \draw[thick] (0,0) -- (0,-2);
             \draw[thick] (3,1) -- (3,-2);
            \fill [blue] (v1) -- (v2) -- (v3) -- cycle;
            \node at (-1,2.8) {$A$};
            \node at (0,-3.2) {$A$};
            \node at (3,-3.2) {$A$};
            \end{tikzpicture}
        \end{array}
    \quad = \quad
        \begin{array}{c}
            \begin{tikzpicture}[scale=0.3,xscale=-1]
            \coordinate (v1) at (-0.5,0.3) {};
            \coordinate (v2) at (0.5,0.3) {};
            \coordinate (v3) at (0,-0.45) {};
            \draw[thick] (-1,2) to[out=-90] ($ (v1) + (0.11,-0.1) $);
            \draw[thick] (1,1) to[out=-90,in=45] ($ (v2) + (-0.11,-0.1) $);
            \draw[thick] (1,1) to[out=90,in=90] (3,1);
            \draw[thick] (0,0) -- (0,-2);
             \draw[thick] (3,1) -- (3,-2);
            \fill [blue] (v1) -- (v2) -- (v3) -- cycle;
            \node at (-1,2.8) {$A$};
            \node at (0,-3.2) {$A$};
            \node at (3,-3.2) {$A$};
            \end{tikzpicture}
        \end{array}.
    \end{align}
So both of these morphisms being equal is certainly a necessary condition. In fact, it is also sufficient (see, for example, Proposition 2.1 in~\cite{MR2075605}). We therefore verify Condition~\eqref{eq:frob_condition} for $ \TM $ via Lemma~\ref{lem:checking_balanced_condition}. This concludes the outline of the strategy.

	\begin{DEF}

	Let $ \bX, \bY $ and $ \bZ $ be given by $ (AB,\e_A^B)\Sharp, (CD,\e_C^D)\Sharp $ and $ (EF,\e_E^F)\Sharp $ respectively. Let $ \alpha $ be in $ \Hom_{\RTC}(\bX,\bY\bZ) = \Hom_{\TC}(\e_A^B,\e_C^D \otimes_{\TC}\e_E^F) $. Recall from Proposition~\ref{prop:characterising_tm_xy} that $ \TM_A^B $ is identified with the subspace of $ \Hom_D(\M(AB),\tid) $ characterised by~\eqref{eq:phase_through}. We consider the map
		\begin{align*}
		\Hom_D(\M(CD),\tid) \otimes \Hom_D(\M(EF),\tid) &\to \Hom_D(\M(AB),\tid)\\
		f \otimes g &\mapsto \TM(\alpha)(f \otimes_\D g).
		\end{align*}
	We note that the image of this map is in $ \TM_{\bX} = \TM_A^B $ as
		\begin{align*}
		\TM(\e_A^B) \circ \TM(\alpha)(f \otimes_\D g) = \TM(\alpha \circ \e_A^B)(f \otimes_\D g) = \TM(\alpha)(f \otimes_\D g).
		\end{align*}
	Therefore restricting this map to the subspace $ \TM_{\bY} \otimes \TM_{\bZ} $ gives a map
		\begin{align*}
		\nabla_{\bX}^{\bY,\bZ}(\alpha) \colon \TM_{\bY} \otimes \TM_{\bZ} \to \TM_{\bX}.
		\end{align*}
	Let $ \nabla: \TM \otimes \TM \to \TM $ be the map constructed from $ \nabla_{\bX}^{\bY,\bZ}(\alpha) $ as described in Remark~\ref{rem:decompose_morphisms}.

	\end{DEF}

    \begin{PROP} \label{prop:algebra_structure}

    The morphisms $ \nabla $ and
        \begin{align*}
        u \defeq \id_{\tid_{\D}} \in \TM_{\tid}^{\tid} = \Hom_{\RTC}(\tidtc, \TM)
        \end{align*}
    form a product/unit pair that make $ \TM $ a \emph{commutative} algebra.

    \proof

	Let $ \bX, \bY $ and $ \bZ $ be given by $ (AB,\e_A^B)\Sharp, (CD,\e_C^D)\Sharp $ and $ (EF,\e_E^F)\Sharp $ respectively. To prove the desired result we shall apply Lemma~\ref{lem:build_an_algebra} by showing that~\eqref{eq:condition_assoc},\eqref{eq:condition_unit} and~\eqref{eq:condition_commutative}  are satisfied. We first note that \eqref{eq:condition_unit} reduces to a triviality in this case.

	To verify~\eqref{eq:condition_assoc} we let $ f,g $ and $ h $ be in $ \TM_A^B,\TM_C^D $ and $ \TM_E^F $ respectively and compute,
		\begin{align*}
		&\nabla^{\bX\bY,\bZ}_{\bX\bY\bZ} (\e_A^B \otimes_{\TC} \e_C^D \otimes_{\TC} \e_E^F) \big( \nabla^{\bX,\bY}_{\bX\bY} (\e_A^B \otimes_{\TC} \e_C^D) (f \otimes g) \otimes h \big) \\
		& = \frac{1}{d(\C)^2} \sum\limits_{S,T} d(S) d(T)
			\begin{array}{c}
				\begin{tikzpicture}[scale=0.15,every node/.style={inner sep=0,outer sep=-1},yscale=-1]
				\node (v1) at (-2,-0.5) {};
				\node (v3) at (2,-0.5) {};
				\node (v2) at (-2,-9) {};
				\node (v4) at (2,-9) {};
				\node (v7) at (10,-9) {};
				\node (v13) at (14,-9) {};
				\node (v9) at (18,-9) {};
				\node (v10) at (6,-0.5) {};
				\node (v6) at (10,-0.5) {};
				\node (v12) at (14,-0.5) {};
				\node (v8) at (18,-0.5) {};
				\node (v11) at (6,-9) {};
				\node [blue] at (-2,-11) {$ A $};
				\node [blue] at (18,-11) {$ F $};
				\node [blue] at (14,-11) {$ E $};
				\node [blue] at (10,-11) {$ D $};
				\node [blue] at (6,-11) {$ C $};
				\node [blue] at (2,-11) {$ B $};
				\node [blue] at (12,4) {$ S $};
				\node [blue] at (20,6) {$ T $};
				\draw [thick,blue] (v1) edge (v2);
				\draw [thick, blue] (v10) edge (v11);
				\draw [thick, blue] (v12) edge (v13);
				\draw [line width= 0.2cm,white] (4,0) ellipse (8 and 4);
				\draw [thick, blue] (4,0) ellipse (8 and 4);
				\draw [line width= 0.2cm,white] (8,0) ellipse (14 and 7);
				\draw [thick, blue] (8,0) ellipse (14 and 7);
				\draw [line width= 0.2cm,white] (v3) edge (v4);
				\draw [line width= 0.2cm,white] (v6) edge (v7);
				\draw [line width= 0.2cm,white] (v8) edge (v9);
				\draw [thick, blue] (v3) edge (v4);
				\draw [thick, blue] (v6) edge (v7);
				\draw [thick, blue] (v8) edge (v9);
				\node [draw,outer sep=0,inner sep=1,minimum height=12,minimum width=24,fill=white] at (0,0) {$  f $};
				\node [draw,outer sep=0,inner sep=1,minimum height=12,minimum width=24,fill=white] at (8,0) {$  g $};
				\node [draw,outer sep=0,inner sep=1,minimum height=12,minimum width=24,fill=white] at (16,0) {$  h $};
				\end{tikzpicture}
			\end{array}\\
		& = f \otimes_D g \otimes_D h
	\end{align*}
	\begin{align*}
		& = \frac{1}{d(\C)^2} \sum\limits_{S,T} d(S) d(T)
			\begin{array}{c}
				\begin{tikzpicture}[scale=0.15,every node/.style={inner sep=0,outer sep=-1},yscale=-1]
				\node (v1) at (-2,-0.5) {};
				\node (v3) at (2,-0.5) {};
				\node (v2) at (-2,-9) {};
				\node (v4) at (2,-9) {};
				\node (v7) at (10,-9) {};
				\node (v13) at (14,-9) {};
				\node (v9) at (18,-9) {};
				\node (v10) at (6,-0.5) {};
				\node (v6) at (10,-0.5) {};
				\node (v12) at (14,-0.5) {};
				\node (v8) at (18,-0.5) {};
				\node (v11) at (6,-9) {};
				\node [blue] at (-2,-11) {$ A $};
				\node [blue] at (18,-11) {$ F $};
				\node [blue] at (14,-11) {$ E $};
				\node [blue] at (10,-11) {$ D $};
				\node [blue] at (6,-11) {$ C $};
				\node [blue] at (2,-11) {$ B $};
				\node [blue] at (4,4) {$ S $};
				\node [blue] at (-4,6) {$ T $};
				\draw [thick,blue] (v1) edge (v2);
				\draw [thick, blue] (v10) edge (v11);
				\draw [thick, blue] (v12) edge (v13);
				\draw [line width= 0.2cm,white] (12,0) ellipse (8 and 4);
				\draw [thick, blue] (12,0) ellipse (8 and 4);
				\draw [line width= 0.2cm,white] (8,0) ellipse (14 and 7);
				\draw [thick, blue] (8,0) ellipse (14 and 7);
				\draw [line width= 0.2cm,white] (v3) edge (v4);
				\draw [line width= 0.2cm,white] (v6) edge (v7);
				\draw [line width= 0.2cm,white] (v8) edge (v9);
				\draw [thick, blue] (v3) edge (v4);
				\draw [thick, blue] (v6) edge (v7);
				\draw [thick, blue] (v8) edge (v9);
				\node [draw,outer sep=0,inner sep=1,minimum height=12,minimum width=24,fill=white] at (0,0) {$  f $};
				\node [draw,outer sep=0,inner sep=1,minimum height=12,minimum width=24,fill=white] at (8,0) {$  g $};
				\node [draw,outer sep=0,inner sep=1,minimum height=12,minimum width=24,fill=white] at (16,0) {$  h $};
				\end{tikzpicture}
			\end{array}\\
		& = \nabla^{\bX,\bY\bZ}_{\bX\bY\bZ} (\e_A^B \otimes_\TC \e_C^D\otimes_{\TC}\e_E^F) \big(  f \otimes \nabla^{\bY,\bZ}_{\bY\bZ}(\e_C^D\otimes_{\TC}\e_E^F)(g \otimes h) \big)
		\end{align*}
	where we are simply using multiple instances of Proposition~\ref{prop:characterising_tm_xy}. Finally, once again by Proposition~\ref{prop:characterising_tm_xy}, we have
		\begin{align*}
		\nabla^{\bZ,\bY}_{\bY\bZ} \big( \hspace{-.3em}
			\begin{array}{c}
				\begin{tikzpicture}[scale=0.25,every node/.style={inner sep=0,outer sep=-1}]
				\node (v4) at (-0.5,-1) {};
				\node (v6) at (1,-1) {};
				\node (v7) at (-0.5,0.5) {};
				\node (v5) at (1,0.5) {};
				\draw [thick] (v4) to[out=90,in=-90] (v5);
				\draw [line width= 0.15cm,white] (v6) to[out=90,in=-90] (v7);
				\draw [thick] (v6) to[out=90,in=-90] (v7);
				\node at (0,1.5) {};
				\end{tikzpicture}
			\end{array}
		\hspace{-.3em} \big) (h \otimes g) &= \frac{1}{d(\C)} \sum\limits_S d(S)
			\begin{array}{c}
				\begin{tikzpicture}[scale=0.15,every node/.style={inner sep=0,outer sep=-1},yscale=-1]
				\node (v1) at (-2,-0.5) {};
				\node (v3) at (2,-0.5) {};
				\node (v2) at (6,-10) {};
				\node (v4) at (10,-10) {};
				\node (v7) at (2,-10) {};
				\node (v10) at (6,-0.5) {};
				\node (v6) at (10,-0.5) {};
				\node (v11) at (-2,-10) {};
				\node [blue] at (6,-12) {$ E $};
				\node [blue] at (2,-12) {$ D $};
				\node [blue] at (-2,-12) {$ C $};
				\node [blue] at (10,-12) {$ F $};
				\node [blue] at (16,-2) {$ S $};
				\draw [thick, blue] (v10) to[out=-90,in=90] (v11);
				\draw [line width= 0.2cm,white] (v1) to[out=-90,in=90] (v2);
				\draw [thick,blue] (v1) to[out=-90,in=90] (v2);
				\draw [line width= 0.2cm,white] (4,-2) ellipse (10 and 6);
				\draw [thick, blue] (4,-2) ellipse (10 and 6);
				\draw [line width= 0.2cm,white] (v3) to[out=-90,in=90] (v4);
				\draw [thick, blue] (v3) to[out=-90,in=90] (v4);
				\draw [line width= 0.2cm,white] (v6) to[out=-90,in=90] (v7);
				\draw [thick, blue] (v6) to[out=-90,in=90] (v7);
				\node [draw,outer sep=0,inner sep=1,minimum height=12,minimum width=24,fill=white] at (0,0) {$  h $};
				\node [draw,outer sep=0,inner sep=1,minimum height=12,minimum width=24,fill=white] at (8,0) {$  g $};
				\end{tikzpicture}
			\end{array}\\
		&= \frac{1}{d(\C)} \sum\limits_S d(S)
			\begin{array}{c}
				\begin{tikzpicture}[scale=0.15,every node/.style={inner sep=0,outer sep=-1},yscale=-1]
				\node (v1) at (6,-2.5) {};
				\node (v3) at (10,-2.5) {};
				\node (v2) at (6,-10) {};
				\node (v4) at (10,-10) {};
				\node (v7) at (2,-10) {};
				\node (v10) at (-2,-2.5) {};
				\node (v6) at (2,-2.5) {};
				\node (v11) at (-2,-10) {};
				\node [blue] at (6,-12) {$ E $};
				\node [blue] at (2,-12) {$ D $};
				\node [blue] at (-2,-12) {$ C $};
				\node [blue] at (10,-12) {$ F $};
				\node [blue] at (16,-2) {$ S $};
				\draw [thick, blue] (v10) to[out=-90,in=90] (v11);
				\draw [line width= 0.2cm,white] (v1) to[out=-90,in=90] (v2);
				\draw [thick,blue] (v1) to[out=-90,in=90] (v2);
				\draw [line width= 0.2cm,white] (4,-2) ellipse (10 and 6);
				\draw [thick, blue] (4,-2) ellipse (10 and 6);
				\draw [line width= 0.2cm,white] (v3) to[out=-90,in=90] (v4);
				\draw [thick, blue] (v3) to[out=-90,in=90] (v4);
				\draw [line width= 0.2cm,white] (v6) to[out=-90,in=90] (v7);
				\draw [thick, blue] (v6) to[out=-90,in=90] (v7);
				\node [draw,outer sep=0,inner sep=1,minimum height=12,minimum width=24,fill=white] at (8,-2) {$  h $};
				\node [draw,outer sep=0,inner sep=1,minimum height=12,minimum width=24,fill=white] at (0,-2) {$  g $};
				\end{tikzpicture}
			\end{array}\\
		&= \nabla^{\bY,\bZ}_{\bY\bZ}(\e_C^D \otimes_{\TC} \e_E^F)(g \otimes h)
		\end{align*}
	which proves~\eqref{eq:condition_commutative}.    \endproof

    \end{PROP}

The next step is to equip $ \TM $ with self-dualizing structure maps. For this to work we need to make some additional assumptions on $ \M\colon \C \to \D $. Firstly we assume that the idempotent completion of $ \D $, denoted $ \DD $, is a multifusion category. Secondly we assume that $ \MM \colon \C \to \DD $, obtained by composing $ \M $ with this embedding, is \emph{indecomposable}. In other words, that there do not exist functors $ \MM_1 \colon \C \to \DD_1 $ and $ \MM_2 \colon \C \to \DD_2 $ such that $ \MM = \MM_1 \oplus \MM_2 $, where $ \DD_i \leq \DD $.

As described in \cite[Section 4.3]{Etingof15}, $ \DD $ decomposes into $ \bigoplus_{i,j \in I} \Dm{i}{j} $ where $ I $ is an indexing set for the primitive idempotents in $ \End_{\DD}(\tid) $. Therefore the condition that $ \MM $ is indecomposable is equivalent to requiring that there exists no subset $ K \subset I $ such that $ \MMX{i}{j}{X} = \MMX{j}{i}{X} = 0 $ for all $ X $ in $ \C $, $ i \in K $ and $ j \in I \setminus K $.

	\begin{PROP} \label{prop:indecomposable_implies_haploid}

	$ \MM $ is indecomposable if and only if $ \TM_{\tid}^{\tid} = \fld $. Furthermore, in this case, any non-zero $ \alpha \in \TM_X^Y \leq \Hom_{\D}(\M(XY),\tid) $ has a left-inverse in $ \Hom_{\D}(\tid,\M(XY)) $ for all $ X,Y $ in $ \C $.

	\proof

	By Proposition~\ref{prop:characterising_tm_xy}, $ \TM_{\tid}^{\tid} $ is given by the subspace of $ \End_{\D}(\tid) $ such that
		\begin{align*}
		\alpha \otimes \id_{\M(Z)} = \id_{\M(Z)} \otimes \ \alpha \quad \forall Z \text{ in }  \C.
		\end{align*}
	Embedding this equality into $ \DD $ and decomposing gives
		\begin{align*}
		\alpha_i \id_{\MMX{i}{j}{Z}}= \alpha_j \id_{\MMX{i}{j}{Z}} \quad \forall Z \text{ in }  \C.
		\end{align*}
	This implies $ \alpha_i = \alpha_j $ for all $ i,j \in I $ if and only if $ \MM $ is indecomposable. This proves the first claim.

	To prove the second claim we recall the characterisation of $ \TM_X^Y $ provided by Proposition~\ref{prop:characterising_tm_xy}, i.e.\ the subspace of $ \Hom_{\D}(\M(XY),\tid) $ such that
		\begin{align*}
		(\alpha  \otimes \id_{\M(Z)}) \circ \phi = \id_{\M(Z)} \otimes \ \alpha \quad \forall Z \text{ in }  \C.
		\end{align*}
	where $ \phi $ is a certain isomorphism. Embedding this equality into $ \DD $ and decomposing gives
		\begin{align*}
		(\alpha_i \otimes \id_{\MMX{i}{j}{Z}}) \circ \tensor[_i]{\phi}{_j} = \id_{\MMX{i}{j}{Z}} \otimes \ \alpha_j \quad \forall Z \text{ in }  \C.
		\end{align*}
	Therefore, if $ \MM $ is indecomposable, $ \alpha_i = 0 $ for any $ i \in I $ implies $ \alpha = 0 $. This proves the second claim.	\endproof

	\end{PROP}

We are now ready to equip $ \TM $ with some self-dualizing structure maps. To accomplish this we shall use Lemma~\ref{lem:dualizing_maps}. We therefore first establish the following perfect pairing.

    \begin{LEMMA} \label{lem:perfect_pairing}

	Let $ X $ and $ Y $ be in $ \C $. As usual $ \TM_X^Y $ is identified with a subspace of $ \Hom_{\D}(\M(XY),\tid) $, however, as described in Remark~\ref{rem:alt_phase_through} we identify $ \TM_{X^\vee}^{Y^\vee} $ with a subspace of $ \Hom_{\D}(\M(Y^\vee X^\vee),\tid) $. The map
        \begin{align*}
        \langle \Blank, \Blank \rangle \colon \TM_X^Y \otimes \TM_{X^\vee}^{Y^\vee} &\to \TM_{\tid}^{\tid} = \fld \\
        f \otimes g &\mapsto
        	\begin{array}{c}
				\begin{tikzpicture}[scale=0.25,every node/.style={inner sep=0,outer sep=-1},yscale=-1]
				\node (v1) at (0.5,-0.25) {};
				\node (v7) at (-1.5,1.5) {};
				\node (v9) at (2.5,1.5) {};
				\draw [thick,blue] (v7) to[out=-90,in=180] (v1);
				\draw [thick,blue] (v9) to[out=-90,in=0] (v1);
				\node (v2) at (-0.75,1.5) {};
				\node (v4) at (1.75,1.5) {};
				\node (v3) at (0.5,0.375) {};
				\draw [thick, blue] (v2) to[out=-90,in=180] (v3);
				\draw [thick, blue] (v3) to[in=-90,out=0] (v4);
				\node [draw,outer sep=0,inner sep=1,minimum width = 13, minimum height = 13, fill=white] at (-1,2) {$ f $};
				\node [draw,outer sep=0,inner sep=1,minimum width = 13, minimum height = 13,fill=white] at (2,2) {$ g $};
				\end{tikzpicture}
        	\end{array}
        \end{align*}
    is a perfect pairing.

    \proof

	Given a non-zero $ f \in \TM_X^Y $, by Proposition~\ref{prop:indecomposable_implies_haploid} there exists $ g \in \Hom_D(\M(Y^\vee X^\vee),\tid) $ such that
		\begin{align*}
			\begin{array}{c}
				\begin{tikzpicture}[scale=0.25,every node/.style={inner sep=0,outer sep=-1},yscale=-1]
				\node (v1) at (0.5,-0.25) {};
				\node (v7) at (-1.5,1.5) {};
				\node (v9) at (2.5,1.5) {};
				\draw [thick,blue] (v7) to[out=-90,in=180] (v1);
				\draw [thick,blue] (v9) to[out=-90,in=0] (v1);
				\node (v2) at (-0.75,1.5) {};
				\node (v4) at (1.75,1.5) {};
				\node (v3) at (0.5,0.375) {};
				\draw [thick, blue] (v2) to[out=-90,in=180] (v3);
				\draw [line width= 0.15cm,white] (v3) to[in=-90,out=0] (v4);
				\draw [thick, blue] (v3) to[in=-90,out=0] (v4);
				\node [draw,outer sep=0,inner sep=1,minimum width = 13, minimum height = 13, fill=white] at (-1,2) {$ f $};
				\node [draw,outer sep=0,inner sep=1,minimum width = 13, minimum height = 13,fill=white] at (2,2) {$ g $};
				\end{tikzpicture}
			\end{array}
		= \id_{\tid}.
		\end{align*}
	We therefore have
	 	\begin{align} \label{eq:perfect_pairing_equation}
		\id_{\tid} = \frac{1}{d(\C)} \sum\limits_S d(S)
			\begin{array}{c}
				\begin{tikzpicture}[scale=0.25,every node/.style={inner sep=0,outer sep=-1},yscale=-1]
				\node (v1) at (0.5,-0.25) {};
				\node (v7) at (-1.5,1.5) {};
				\node (v9) at (2.5,1.5) {};
				\draw [thick,blue] (v7) to[out=-90,in=180] (v1);
				\draw [thick,blue] (v9) to[out=-90,in=0] (v1);
				\node (v2) at (-0.75,1.5) {};
				\node (v4) at (1.75,1.5) {};
				\node (v3) at (0.5,0.375) {};
				\draw [thick, blue] (v2) to[out=-90,in=180] (v3);
				\draw [thick, blue] (v3) to[in=-90,out=0] (v4);
				\node [draw,outer sep=0,inner sep=1,minimum width = 13, minimum height = 13, fill=white] at (-1,2) {$ f $};
				\node [draw,outer sep=0,inner sep=1,minimum width = 13, minimum height = 13,fill=white] at (2,2) {$ g $};
				\draw [thick, blue] (0.5,2) ellipse (3 and 3);
				\node[blue] at (4.5,2) {$ S $};
				\end{tikzpicture}
			\end{array}
		= \frac{1}{d(\C)} \sum\limits_S d(S)
			\begin{array}{c}
				\begin{tikzpicture}[scale=0.25,every node/.style={inner sep=0,outer sep=-1},yscale=-1]
				\node (v1) at (0.5,-0.25) {};
				\node (v7) at (-1.5,1.5) {};
				\node (v9) at (2.5,1.5) {};
				\draw [thick,blue] (v7) to[out=-90,in=180] (v1);
				\draw [thick,blue] (v9) to[out=-90,in=0] (v1);
				\node (v2) at (-0.75,1.5) {};
				\node (v4) at (1.75,1.5) {};
				\node (v3) at (0.5,0.375) {};
				\draw [thick, blue] (v2) to[out=-90,in=180] (v3);
				\draw [line width= 0.15cm,white] (2,2) ellipse (1.6 and 1.6);
				\draw [thick, blue] (2,2) ellipse (1.6 and 1.6);
				\draw [line width= 0.15cm,white] (v3) to[in=-90,out=0] (v4);
				\draw [thick, blue] (v3) to[in=-90,out=0] (v4);
				\node [draw,outer sep=0,inner sep=1,minimum width = 13, minimum height = 13, fill=white] at (-1,2) {$ f $};
				\node [draw,outer sep=0,inner sep=1,minimum width = 13, minimum height = 13,fill=white] at (2,2) {$ g $};
				\node[blue] at (4.5,2) {$ S $};
				\end{tikzpicture}
			\end{array}
	 	\end{align}
	by Proposition~\ref{prop:characterising_tm_xy}. We now consider $ \tilde{g} = \TM(\tilde{\e}_{X^\vee}^{Y^\vee})(g) \in \TM_{X^\vee}^{Y^\vee} $ (where $ \tilde{\e}_{X^\vee}^{Y^\vee} $ is given by~\eqref{eq:alt_idem}). Then the right-hand side of~\eqref{eq:perfect_pairing_equation} is $ \langle f, \tilde{g} \rangle $ and so we are done.	\endproof

	\end{LEMMA}

	\begin{REM} \label{rem:symmetric_pairing}

	We note that this perfect pairing is \emph{symmetric} with respect to the pivotal structure, i.e.\ $ \langle f , g \rangle = \langle g , f \rangle $ where $ f \in \TM_X^Y = \TM_{X^{\vee\vee}}^{Y^{\vee\vee}} $.

	\end{REM}

    \begin{PROP} \label{prop:tm_is_frob}

    We consider $ \TM $ equipped with the algebra structure from Proposition \ref{prop:algebra_structure}. We also equip $ \TM $ with the self-dualizing maps given by Lemma~\ref{lem:perfect_pairing} and Lemma~\ref{lem:dualizing_maps} with $ c = d $ (the dimension map for $ \RTC $). Then $ \TM $ satisfies \eqref{eq:frob_condition}, i.e.\ is a Frobenius algebra.

    \proof

    Let $ \bR, \bS $ and $ \bT $ be given by $ (IJ,\e_I^J)\Sharp, (KL,\e_K^L)\Sharp $ and $ (MN,\e_M^N)\Sharp $ respectively where $ I,J,K,L,M,N \in \I $. Let $ f,g $ and $ h $ be in $ \TM_I^J,\TM_K^L $ and $ \TM_M^N $ respectively and let $ \beta $ be in $ \Hom_{\RTC}(\bS\bT,\bR) = \Hom_{\TC}(\e_K^L \otimes_{\TC} \e_M^N, \e_I^J) $. We have
        \begin{align*}
		&d(\bS)\left( g^* \circ \nabla_{\bS}^{\bR,\bT^\vee}\left(
				\begin{array}{c}
						\begin{tikzpicture}[scale=0.25,every node/.style={inner sep=0,outer sep=-1}]
						\draw[thick] (0,0) -- (0,-2);
						\draw[thick] (-0.6,0) -- (-0.6,2);
						\draw[thick] (2,0.7)-- (2,-2);
						\draw[thick] (0.6,0.7) to[out=90,in=90] (2,0.7);
						\node [draw,outer sep=0,inner sep=2,minimum width=20,fill=white] at (0,0) {$\scriptstyle \beta $};
						\node at (-0.6,-1.65) {$\scriptscriptstyle \bR $};
						\node at (1.3,-1.55) {$\scriptscriptstyle \bT^\vee $};
						\node at (-1.1,1.6) {$\scriptscriptstyle \bS $};
						\end{tikzpicture}
				\end{array}
		\right) \right)(f \otimes (h^*)^\vee) \\
		= \ &d(\bS) \ g^* \left(
			\begin{array}{c}
                \begin{tikzpicture}[scale=0.25,every node/.style={inner sep=0,outer sep=-1}]
                \node (v4) at (-1.5,9.5) {};
                \node (v6) at (-0.75,8.75) {};
                \node (v7) at (0,8) {};
                \node (v9) at (0.75,7.25) {};
                \node (v3) at (-0.25,11.5) {};
                \node (v5) at (1,11.5) {};
                \draw [thick, blue] (v3) to[out=-90,in=45] (v4);
                \draw [thick, blue] (v5) to[out=-90,in=45] (v6);
                \node (v10) at (2.25,6.75) {};
                \node (v8) at (3.75,6.5) {};
                \draw [thick, blue] (v7) to[out=45,in=90] (v8);
                \draw [thick, blue] (v9) to[out=45,in=90] (v10);
                \node (v11) at (-0.25,3.75) {};
                \node (v12) at (1.75,3.75) {};
                \draw [thick, blue] (v10) to[out=-90,in=90] (v11);
                \node (v13) at (-2.5,7.25) {};
                \node (v15) at (-1.5,6.25) {};
                \node (v14) at (-5,3.75) {};
                \node (v16) at (-3.5,3.75) {};
                \draw [thick, blue] (v13) to[in=90,out=-135] (v14);
                \draw [thick, blue] (v15) to[in=90,out=-135] (v16);
                \node (v17) at (0.25,6.25) {};
                \node (v20) at (-2.25,8.75) {};
                \node (v18) at (4.25,2.25) {};
                \node (v19) at (-7,2.25) {};
                \draw [thick, blue] (v18) to[out=-90,in=-90] (v19);
                \draw [thick, blue] (v20) to[out=135,in=90] (v19);
                \draw [line width= 0.2cm,white] (v17) to[out=-45,in=90] (v18);
                \draw [thick, blue] (v17) to[out=-45,in=90] (v18);
                \draw [line width= 0.2cm,white] (v8) to[out=-90,in=90] (v12);
                \draw [thick, blue] (v8) to[out=-90,in=90] (v12);
                \node [draw,outer sep=0,inner sep=2,minimum size=16,fill=white] (v1) at (-4.25,2.75) {$ f $};
                \node [draw,outer sep=0,inner sep=2,minimum size=16,fill=white] (v2) at (0.75,2.75) {$ (h^*)^\vee $};
                \draw [blue,rotate=-45,fill=white] (-8.298,6.1266) rectangle (-3.798,3.1266);
                \node [blue] at (2,11.25) {$ L $};
                \node [blue] at (-1.5,11.25) {$ K $};
                \node [blue] at (4.75,0.25) {$ G $};
                \node [blue] at (-1,7.5) {$ \beta $};
                \end{tikzpicture}
			\end{array} \right)
		=
			\begin{array}{c}
                \begin{tikzpicture}[scale=0.25,every node/.style={inner sep=0,outer sep=-1}]
				\node (v4) at (-1.5,9.5) {};
				\node (v6) at (-0.5,9) {};
				\node (v7) at (0,8) {};
				\node (v9) at (0.75,7.25) {};
				\node (v3) at (-0.75,11) {};
				\node (v5) at (0.5,11.25) {};
				\draw [thick, blue] (v3) to[out=-90,in=45] (v4);
				\draw [thick, blue] (v5) to[out=-90,in=45] (v6);
				\node (v10) at (3,6.75) {};
				\node (v8) at (4.5,6.75) {};
				\draw [thick, blue] (v7) to[out=45,in=90] (v8);
				\draw [thick, blue] (v9) to[out=45,in=90] (v10);
				\node (v13) at (-2.5,7.25) {};
				\node (v15) at (-1.5,6.25) {};
				\node (v14) at (-4.5,4.75) {};
				\node (v16) at (-3,4.75) {};
				\draw [thick, blue] (v13) to[in=90,out=-135] (v14);
				\draw [thick, blue] (v15) to[in=90,out=-135] (v16);
				\node (v17) at (0.25,6.25) {};
				\node (v20) at (-2.25,8.75) {};
				\node (v18) at (-1.75,1.25) {};
				\node (v19) at (-7,6.5) {};
				\draw [thick, blue] (v18) to[out=-135,in=-135] (v19);
				\draw [thick, blue] (v20) to[out=135,in=45] (v19);
				\draw [thick, blue] (v17) to[out=-45,in=45] (v18);
				\node (v11) at (3,4.75) {};
				\node (v12) at (4.5,4.75) {};
				\draw [thick, blue] (v10) edge (v11);
				\draw [thick, blue] (v8) edge (v12);
				\node (v21) at (-8.25,11) {};
				\node (v22) at (-6.5,11) {};
				\draw [thick, blue] (v3) to[out=90,in=90] (v22);
				\draw [thick, blue] (v5) to[out=90,in=90] (v21);
				\node [draw,outer sep=0,inner sep=2,minimum size=16,fill=white] (v1) at (-3.75,3.75) {$ f $};
				\node [draw,outer sep=0,inner sep=2,minimum size=16,fill=white] (v2) at (3.75,3.75) {$ (h^*)^\vee $};
				\draw [blue,rotate=-45,fill=white] (-8.298,6.1266) rectangle (-3.798,3.1266);
				\node [blue] at (1.75,11.25) {$ L $};
				\node [blue] at (-2.25,11.25) {$ K $};
				\node [blue] at (0.5,1.75) {$ G $};
				\node [blue] at (-1,7.5) {$ \beta $};
				\node [draw,outer sep=0,inner sep=2,minimum size=16,fill=white] (v2) at (-7.25,10) {$ (g^*)^\vee $};
				\end{tikzpicture}
			\end{array}\\
		= \ &d(\bT) \ h^* \left(
			\begin{array}{c}
                \begin{tikzpicture}[scale=0.25,every node/.style={inner sep=0,outer sep=-1}]
				\node (v4) at (-1.5,9.5) {};
				\node (v6) at (-0.5,9) {};
				\node (v7) at (0,8) {};
				\node (v9) at (0.75,7.25) {};
				\node (v3) at (-0.75,10.5) {};
				\node (v5) at (0.5,10.75) {};
				\draw [thick, blue] (v3) to[out=-90,in=45] (v4);
				\draw [thick, blue] (v5) to[out=-90,in=45] (v6);
				\node (v10) at (4.5,13.25) {};
				\node (v8) at (3.25,13.25) {};
				\draw [thick, blue] (v7) to[out=45,in=-90] (v8);
				\draw [thick, blue] (v9) to[out=45,in=-90] (v10);
				\node (v13) at (-2.5,7.25) {};
				\node (v15) at (-1.5,6.25) {};
				\node (v14) at (-3.25,5) {};
				\node (v16) at (-1.75,5) {};
				\draw [thick, blue] (v13) to[in=90,out=-135] (v14);
				\draw [thick, blue] (v15) to[in=90,out=-135] (v16);
				\node (v17) at (0.25,6.25) {};
				\node (v20) at (-2.25,8.75) {};
				\node (v18) at (-8,1.25) {};
				\node (v19) at (-8,8.5) {};
				\draw [thick, blue] (v18) to[out=135,in=-135] (v19);
				\draw [thick, blue] (v17) to[out=-45,in=-45] (v18);
				\node (v21) at (-7,10.75) {};
				\node (v22) at (-5.5,10.5) {};
				\draw [thick, blue] (v3) to[out=90,in=90] (v22);
				\draw [thick, blue] (v5) to[out=90,in=90] (v21);
				\node (v11) at (-7,4.75) {};
				\node (v12) at (-5.5,4.75) {};
				\draw [thick, blue] (v22) edge (v12);
				\draw [line width= 0.2cm,white] (v20) to[out=135,in=45] (v19);
				\draw [thick, blue] (v20) to[out=135,in=45] (v19);
				\draw [line width= 0.2cm,white] (v21) edge (v11);
				\draw [thick, blue] (v21) edge (v11);
				\node [draw,outer sep=0,inner sep=2,minimum size=16,fill=white] (v1) at (-2.5,4) {$ f $};
				\draw [blue,rotate=-45,fill=white] (-8.298,6.1266) rectangle (-3.798,3.1266);
				\node [blue] at (1.75,12.75) {$ M $};
				\node [blue] at (5.75,12.75) {$ N $};
				\node [blue] at (0.5,1.75) {$ G $};
				\node [blue] at (-1,7.5) {$ \beta $};
				\node [draw,outer sep=0,inner sep=2,minimum size=16,fill=white] (v2) at (-6.25,4) {$ (g^*)^\vee $};
				\end{tikzpicture}
			\end{array} \right)\\
		= \ &d(\bT) \left( h^* \circ \nabla_{\bT}^{\bS^\vee,\bR}\left(
		 	\begin{array}{c}
				\begin{tikzpicture}[scale=0.25,every node/.style={inner sep=0,outer sep=-1}, xscale=-1]
				\draw[thick] (0,0) -- (0,-2);
				\draw[thick] (-0.6,0) -- (-0.6,2);
				\draw[thick] (2,0.7)-- (2,-2);
				\draw[thick] (0.6,0.7) to[out=90,in=90] (2,0.7);
				\node [draw,outer sep=0,inner sep=2,minimum width=20,fill=white] at (0,0) {$\scriptstyle \beta $};
				\node at (-0.5,-1.65) {$\scriptscriptstyle \bR $};
				\node at (1.2,-1.5) {$\scriptscriptstyle \bS^\vee $};
				\node at (-1.1,1.6) {$\scriptscriptstyle \bT $};
				\end{tikzpicture}
			\end{array}
		\right) \right)((g^*)^\vee \otimes f)
        \end{align*}
	where we have used Proposition~\ref{prop:characterising_tm_xy} and Lemma~\ref{lem:dualizing_maps} multiple times. \endproof
\end{PROP}
	\begin{THM} \label{thm:tm_is_frob}

	Let $ \C $ be an MTC and let $ \M $ be a pivotal tensor functor from $ \C $ to $ \D $ such that $ \MM $ is indecomposable. Then $ \TM $ is a haploid, symmetric, commutative Frobenius algebra.

	\proof

	This result follows from Proposition~\ref{prop:indecomposable_implies_haploid}, Proposition~\ref{prop:algebra_structure}, Proposition~\ref{prop:tm_is_frob} and Remark~\ref{rem:symmetric_pairing}.	\endproof

	\end{THM}

    \begin{REM} \label{rem:bruguieres_natale}

    The work of Buguieres and Natale~\cite{MR2863377,MR3161401} allows for an alternative perspective $ \TM $ and the algebraic structures it admits. They show that if a monoidal functor $ \M \colon \C \to \D $ admits a right adjoint $ \M\ad \colon \D \to \C  $ then $ \M\ad(\tid_{\D}) $ admits a canonical half-braiding and product yielding a commutative algebra in $ Z(\C) $~\cite[Proposition 5.1]{MR3161401}. As we have supposed $ \C $ to be modular and therefore semisimple, the functor
        \begin{align*}
        \M\ad \colon \D &\to \C \\
        A &\mapsto \bigoplus_{S} \Hom_D(\M(S),A) \cdot S
        \end{align*}
    is right adjoint to $ \M $ and the resulting object $ \M\ad(\tid_{\D}) $ in $ Z(\C) = \RTC $ coincides with $ \TM $.
    \end{REM}

Combining Theorem~\ref{thm:tm_is_frob} with the aforementioned \hyperref[thm:kong_runkel]{result of Kong and Runkel}~\cite[Theorem 3.4]{MR2551797} we obtain the following.

    \begin{COR} \label{cor:tm_mod_inv}

    Let $ \C $ be an MTC and let $ \M $ be a pivotal tensor functor from $ \C $ to $ \D $ such that $ \MM $ is indecomposable and $ \sum_{IJ} \dim(\TM_I^J) d(I) d(J) = d(\C) $. Then $ Z(\TM) $ is a modular invariant.

    \proof

    As $ \sum_{IJ} \dim(\TM_I^J) d(I) d(J) $ is the dimension of $ \TM $ as an object in $ \RTC \cong \CC $, the claim follows immediately from Theorem~\ref{thm:tm_is_frob} and\hyperref[thm:kong_runkel]{~\cite[{Theorem 3.4}]{MR2551797}}. \endproof

    \end{COR}

\section{Module Categories and $ \alpha $-induction} \label{sec:mod_cat_and_alpha_ind}

	\begin{DEF}

	Let $ \C $ be a monoidal category. A \emph{module category} over $ \C $ is a monoidal category $ \ModCat $ together with a monoidal functor $ \M \colon \C \to \End(\ModCat) $, where $ \End(\ModCat) $ is the category of endofunctors on $ \ModCat $. For the remainder of this article all module categories are assumed to be semisimple with finitely many simple objects.

	\end{DEF}

From a physical point of view Cardy~\cite{MR1048596} showed that the algebraic data of an \emph{annular partition function} in a boundary (rational) conformal field theory is given by a module category over the corresponding MTC. The process known as $ \alpha $-induction is an operator algebra technique developed by B\"{o}ckenhauer and Evans~\cite{MR1652746} that produces a \emph{toroidal partition function} (as described in the introduction) from an annular partition function. The aim of this section is to apply the $ \TM $ construction in the case when $ \M $ is a module category. In particular, we shall show that in this case $ \TM $ may be related to $ \alpha $-induction.

Ostrik~\cite[Section 5]{MR1976459} rephrased $ \alpha $-induction using categorical language in the following way. Let $ \M \colon \C \to \D $ be a module category over a pre-modular tensor category $ \C $, where $ \D $ denotes $ \End(\B) $ and $ \B $ is a semisimple category with finitely many simple objects. For $ X,Y $ in $ \C $ we consider the subspace
 	\begin{align*}
 	\Ostrik{X}{Y} \leq \Hom_{\D}(\M(X),\M(Y))
 	\end{align*}
defined by the condition that $ \beta \in \Ostrik{X}{Y} $ satisfies, for all $ Z $ in $ \C $,
		\begin{equation}
			\begin{tikzcd}[row sep=2cm,column sep=2cm,inner sep=1ex] \label{eq:commuting_diag_alpha}
			\M(X) \otimes \M(Z) \arrow[swap]{d}[name=D]{\beta \otimes \id}   \arrow{r}{\M(\overline{\sigma}_{ZX})} &  \M(Z) \otimes \M(X) \arrow{d}[name=U]{\id \otimes \beta}
			\\
			\M(Y) \otimes \M(Z) \arrow{r}{\M(\sigma_{ZY})} & \M(Z) \otimes \M(Y)
			\arrow[to path={(U) node[scale=2,midway] {$\circlearrowleft$}  (D)}]{}
			\end{tikzcd}
		\end{equation}
where $ \sigma $ and $ \overline{\sigma} $ are the braiding on $ \C $ and its opposite respectively. The principal claim of $ \alpha $-induction is then as follows. Under the assumption that the dimensions of all the objects in $ \C $ are positive, the $ \I \times \I $-matrix whose entries are given by the dimension of $ \Ostrik{I}{J^\vee} $ commutes with the modular data of $ \C $. Furthermore if $ \M $ is irreducible then this matrix is a \emph{modular invariant} (see Definition~\ref{def:modular_invariant}).

	\begin{REM}

	The claim found in~\cite{MR1976459} is actually that the $ \I \times \I $-matrix whose entries are given by the dimension of $ \Ostrik{I}{J} $ commutes with the modular data of $ \C $. However as the modular data always commutes with the charge conjugation matrix these statements are equivalent.

	\end{REM}

In~\cite{MR1976459} Ostrik also provides the following example to prove the necessity of the condition that the objects in $ \C $ have positive dimension.

	\begin{EX} \label{ex:ostrik_alpha_fail}

	Let $ \C $ be the fusion category of representations of $ \mathbb{Z}/2\mathbb{Z} $. A complete set of simples in $ \C $ is given by $ \{\underline{0}, \underline{1} \} $, where $ \underline{0} $ is the tensor unit and $ \underline{1} \otimes \underline{1} = \underline{0} $. We may then equip $ \C $ with the pivotal structure $ \delta_{\underline{1}} = - \id_{\underline{1}} $ (so that $ d(\underline{1}) = -1 $). One may also check that setting $ \sigma_{\underline{1},\underline{1}} = \id_{\underline{0}} $ defines a (degenerate) braiding on the category and we obtain a pre-modular tensor category. We then consider the module category
       \begin{align*}
       \M\colon \C &\to \Vect \\
	   \underline{0} &\mapsto \mathbb{K} \\
       \underline{1} &\mapsto \mathbb{K}.
   \end{align*}
   As the braiding is given by the identity, we have $ \sigma = \overline{\sigma} $ and Equation~\eqref{eq:commuting_diag_alpha} reduces to a tautology. Therefore $ \Ostrik{\underline{0}}{\underline{1}} = \fld $ and the resulting dimension matrix fails to commute with the T-matrix
   $$ \Tc = \left( \begin{array}{cc}
   1 & 0 \\ 0& -1
   \end{array} \right). $$

	\end{EX}

We start by remarking that Condition~\ref{eq:commuting_diag_alpha} makes sense even when $ \D $ is an arbitrary tensor category. Therefore to connect Ostrik's formulation of $ \alpha $-induction to $ \TM $ we have the following.

 	\begin{THM} \label{thm:alpha_induction_equivalence}

	Let $ \C $ be a pre-modular tensor category and let $ \M \colon \C \to \D $ be a pivotal monoidal functor. Then $ \Ostrik{X}{Y^\vee} \cong \TM_X^Y $.

	\proof

	Graphically Condition~\eqref{eq:commuting_diag_alpha} is given by
		\begin{align} \label{eq:graphical_alpha_condition}
			\begin{array}{c}
				\begin{tikzpicture}[scale=0.15,every node/.style={inner sep=0,outer sep=-1}]
				\node [draw,outer sep=0,inner sep=2,minimum size=10] (v2) at (4,-7) {$ \beta  $};
				\node (v8) at (0.5,-0.5) {};
				\node (v3) at (4,-10) {};
				\draw [thick, blue] (v3) edge (v2);
				\node (v6) at (1.5,-1.5) {};
				\node (v4) at (2,0) {};
				\node (v5) at (0,-2) {};
				\draw [thick, blue] (v4) edge (v5);
				\node (v9) at (-2,4) {};
				\draw [thick, blue] (v6) to[out=-45,in=90] (v2);
				\node (v1) at (-2,-5.5) {};
				\node (v7) at (-2,-10) {};
				\draw [thick, blue] (v5) to[out=-135,in=90] (v1);
				\draw [thick, blue] (v7) edge (v1);
				\draw [thick, blue] (v8) to[out=135,in=-90] (v9);
				\node [blue] (v10) at (-2,6) {$ X $};
				\node [blue] (v11) at (4,6) {$ Z $};
				\node [blue] (v10) at (-2,-12) {$ Z $};
				\node [blue] (v11) at (4,-12) {$ Y^\vee $};
				\node (v12) at (4,4) {};
				\draw [thick, blue] (v4) to[out=45,in=-90] (v12);
				\end{tikzpicture}
			\end{array}
		=
			\begin{array}{c}
				\begin{tikzpicture}[scale=0.15,every node/.style={inner sep=0,outer sep=-1}]
				\node (v2) at (4,4) {};
				\node (v8) at (1,-5) {};
				\node (v3) at (4,-10) {};
				\node (v6) at (1,-5) {};
				\node (v4) at (1.5,-4.5) {};
				\node (v5) at (0.5,-5.5) {};
				\node [draw,outer sep=0,inner sep=2,minimum size=10] (v9) at (-2,1) {$ \beta $};
				\draw [thick, blue] (v6) to[out=-45,in=90] (v3);
				\node (v1) at (-2,4) {};
				\node (v7) at (-2,-10) {};
				\draw [thick, blue] (v5) to[out=-135,in=90] (v7);
				\draw [thick, blue] (v8) to[out=135,in=-90] (v9);
				\node [blue] (v10) at (-2,6) {$ X $};
				\node [blue] (v11) at (4,6) {$ Z $};
				\node [blue] (v10) at (-2,-12) {$ Z $};
				\node [blue] (v11) at (4,-12) {$ Y^\vee $};
				\node (v12) at (4,1) {};
				\draw [thick, blue] (v4) to[out=45,in=-90] (v12);
				\draw [thick, blue] (v9) edge (v1);
				\draw [thick, blue] (v2) edge (v12);
				\end{tikzpicture}
			\end{array}
		\end{align}
	for all $ Z $ in $ \C $. As $ \TM_X^Y $ is a subspace of $ \Hom_{\D}(\M(XY),\tid) = \Hom_{\D}(\M(X),\M(Y^\vee)) $ we only have to check that Condition~\eqref{eq:graphical_alpha_condition} is equivalent to Condition~\eqref{eq:phase_through}.

	Suppose $ \beta \in \Hom_{\D}(\M(X),\M(Y^\vee)) $ satisfies Condition \eqref{eq:commuting_diag_alpha}. Then we have
		\begin{align*}
			\begin{array}{c}
                \begin{tikzpicture}[scale=0.15,every node/.style={inner sep=0,outer sep=-1},yscale=-1]
				\node (v9) at (-2,3) {};
				\node (v1) at (-4.5,3) {};
				\node (v10) at (-4.5,0.5) {};
				\draw [thick, blue] (v9) to[out=90,in=90] (v1);
				\node (v2) at (0,0.5) {};
				\node (v3) at (-7,-4) {};
				\node (v4) at (-4.5,-4) {};
				\node (v5) at (-2,-4) {};
				\draw [thick, blue] (v10) to[out=-90,in=90] (v4);
				\draw [line width=0.2cm,white] (v2) to[out=-90,in=90] (v3);
				\draw [thick, blue] (v2) to[out=-90,in=90] (v3);
				\draw [line width=0.2cm,white] (v9) to[out=-90,in=90] (v5);
				\draw [thick, blue] (v9) to[out=-90,in=90] (v5);
				\node (v6) at (0,4.5) {};
				\draw [thick, blue] (v6) edge (v2);
				\node[blue] at (-7,-5.5) {$ Z $};
				\node[blue] at (-4.5,-5.5) {$ X $};
				\node[blue] at (-2,-5.5) {$ Y $};
				\node [draw,outer sep=0,inner sep=2,minimum size=10,fill=white] at (-4.5,1.5) {$ \beta $};
				\end{tikzpicture}
			\end{array}
		=
			\begin{array}{c}
                \begin{tikzpicture}[scale=0.15,every node/.style={inner sep=0,outer sep=-1},yscale=-1]
                \node (v9) at (-2,2) {};
                \node (v1) at (-4.5,2) {};
                \node (v10) at (-4.5,-3) {};
                \node (v2) at (0,3) {};
                \node (v3) at (-7,0) {};
                \node (v4) at (-4.5,-4) {};
                \node (v5) at (-2,-4) {};
                \draw [thick, blue] (v10) to[out=-90,in=90] (v4);
                \draw [line width=0.2cm,white] (v2) to[out=-90,in=90] (v3);
                \draw [thick, blue] (v2) to[out=-90,in=90] (v3);
                \draw [line width=0.2cm,white] (v9) to (v5);
                \draw [thick, blue] (v9) to (v5);
                \node (v6) at (0,4.5) {};
                \draw [thick, blue] (v6) edge (v2);
                \node[blue] at (-7,-5.5) {$ Z $};
                \node[blue] at (-4.5,-5.5) {$ X $};
                \node[blue] at (-2,-5.5) {$ Y $};
                \node (v7) at (-7,-4) {};
                \draw [thick, blue] (v7) edge (v3);
                \draw [line width=0.2cm,white] (v1) edge (v10);
                \draw [thick, blue] (v1) edge (v10);
                \draw [thick, blue] (v9) to[out=90,in=90] (v1);
                \node [draw,outer sep=0,inner sep=2,minimum size=10,fill=white] at (-4.5,-1.5) {$ \beta $};
                \end{tikzpicture}
			\end{array}
		=
			\begin{array}{c}
                \begin{tikzpicture}[scale=0.15,every node/.style={inner sep=0,outer sep=-1},yscale=-1]
				\node (v9) at (-2,2.5) {};
				\node (v1) at (-4.5,2.5) {};
				\draw [thick, blue] (v9) to[out=90,in=90] (v1);
				\node (v2) at (-7,5.5) {};
				\node (v3) at (-7,-4) {};
				\node (v4) at (-4.5,-4) {};
				\node (v5) at (-2,-4) {};
				\draw [thick, blue] (v1) to[out=-90,in=90] (v4);
				\draw [line width=0.2cm,white] (v2) to[out=-90,in=90] (v3);
				\draw [thick, blue] (v2) to[out=-90,in=90] (v3);
				\draw [line width=0.2cm,white] (v9) to[out=-90,in=90] (v5);
				\draw [thick, blue] (v9) to[out=-90,in=90] (v5);
				\node (v6) at (-7,6) {};
				\draw [thick, blue] (v6) edge (v2);
				\node[blue] at (-7,-5.5) {$ Z $};
				\node[blue] at (-4.5,-5.5) {$ X $};
				\node[blue] at (-2,-5.5) {$ Y $};
				\node [draw,outer sep=0,inner sep=2,minimum size=10,fill=white] (v9) at (-4.5,-0.5) {$ \beta $};
				\end{tikzpicture}
			\end{array}.
		\end{align*}
	Furthermore, for $ \alpha \in \TM_I^J $, we have
		\begin{align*}
			\begin{array}{c}
                \begin{tikzpicture}[scale=0.15,every node/.style={inner sep=0,outer sep=-1}]
                \node [blue] (v10) at (-4.5,4) {$ X $};
                \node [blue] (v11) at (-1.5,4) {$ Z $};
                \node [blue] (v10) at (-7,-9) {$ Z $};
                \node [blue] (v11) at (2,-9) {$ Y^\vee $};
                \node (v1) at (-0.5,-3) {};
                \node (v4) at (2,-7) {};
                \node (v2) at (-0.5,-2) {};
                \node (v3) at (2,-2) {};
                \node (v7) at (-7,-7) {};
                \node (v8) at (-7,-3) {};
                \node (v9) at (-1.5,2) {};
                \node (v5) at (-4.5,-3) {};
                \node (v6) at (-4.5,2) {};
                \draw [thick, blue] (v5) to[out=90,in=-90] (v6);
                \draw [thick, blue] (v7) edge (v8);
                \draw [line width=0.2cm,white] (v8) to[out=90,in=-90] (v9);
                \draw [thick, blue] (v8) to[out=90,in=-90] (v9);
                \draw [line width=0.2cm,white] (v1) edge (v2);
                \draw [thick, blue] (v1) edge (v2);
                \draw [line width=0.2cm,white] (v3) edge (v4);
                \draw [thick, blue] (v3) edge (v4);
                \draw [thick, blue] (v2) to[out=90,in=90] (v3);
                \node [draw,outer sep=0,inner sep=2,minimum width=24,minimum height=12,fill=white] (v10) at (-2.5,-4) {$\alpha$};
                \end{tikzpicture}
			\end{array}
		=
			\begin{array}{c}
                \begin{tikzpicture}[scale=0.15,every node/.style={inner sep=0,outer sep=-1}]
                \node [blue] (v10) at (-4.5,4) {$ X $};
                \node [blue] (v11) at (3.5,4) {$ Z $};
                \node [blue] (v10) at (-7,-9) {$ Z $};
                \node [blue] (v11) at (2,-9) {$ Y^\vee $};
                \node (v1) at (-0.5,-3) {};
                \node (v4) at (2,-7) {};
                \node (v2) at (-0.5,1) {};
                \node (v3) at (2,1) {};
                \node (v7) at (-7,-7) {};
                \node (v8) at (-7,-3) {};
                \node (v9) at (3.5,2) {};
                \node (v5) at (-4.5,-3) {};
                \node (v6) at (-4.5,2) {};
                \draw [thick, blue] (v5) to[out=90,in=-90] (v6);
                \draw [thick, blue] (v7) edge (v8);
                \draw [line width=0.2cm,white] (v8) to[out=90,in=-90] (v9);
                \draw [thick, blue] (v8) to[out=90,in=-90] (v9);
                \draw [line width=0.2cm,white] (v1) edge (v2);
                \draw [thick, blue] (v1) edge (v2);
                \draw [line width=0.2cm,white] (v3) edge (v4);
                \draw [thick, blue] (v3) edge (v4);
                \draw [thick, blue] (v2) to[out=90,in=90] (v3);
                \node [draw,outer sep=0,inner sep=2,minimum width=24,minimum height=12,fill=white] (v10) at (-2.5,-4) {$\alpha$};
                \end{tikzpicture}
			\end{array}
		=
			\begin{array}{c}
                \begin{tikzpicture}[scale=0.15,every node/.style={inner sep=0,outer sep=-1}]
                \node [blue] (v10) at (-4.5,4) {$ X $};
                \node [blue] (v11) at (4.5,4) {$ Z $};
                \node [blue] (v10) at (-2,-9) {$ Z $};
                \node [blue] (v11) at (2,-9) {$ Y^\vee $};
                \node (v1) at (-0.5,0) {};
                \node (v4) at (2,-7) {};
                \node (v2) at (-0.5,1) {};
                \node (v3) at (2,1) {};
                \node (v7) at (-2,-7) {};
                \node (v8) at (-2,-6) {};
                \node (v9) at (4.5,-2) {};
                \node (v5) at (-4.5,0) {};
                \node (v6) at (-4.5,2) {};
                \draw [thick, blue] (v5) to[out=90,in=-90] (v6);
                \draw [thick, blue] (v7) edge (v8);
                \draw [line width=0.2cm,white] (v8) to[out=90,in=-90] (v9);
                \draw [thick, blue] (v8) to[out=90,in=-90] (v9);
                \draw [line width=0.2cm,white] (v1) edge (v2);
                \draw [thick, blue] (v1) edge (v2);
                \draw [line width=0.2cm,white] (v3) edge (v4);
                \draw [thick, blue] (v3) edge (v4);
                \draw [thick, blue] (v2) to[out=90,in=90] (v3);
                \node [draw,outer sep=0,inner sep=2,minimum width=24,minimum height=12,fill=white] (v10) at (-2.5,-1) {$\alpha$};
                \node (v12) at (4.5,2) {};
                \draw [thick, blue] (v12) edge (v9);
                \end{tikzpicture}
			\end{array}
		\end{align*}
	where the final equality uses Proposition~\ref{prop:characterising_tm_xy}. This is equivalent to Condition~\eqref{eq:commuting_diag_alpha} as desired. \endproof

	\end{THM}

The alternative characterisation of $ \TM_I^J $ given by Theorem~\ref{thm:alpha_induction_equivalence} allows for the following generalization of Corollary~\ref{cor:t-inv} to the pre-modular case.

 	\begin{THM} \label{thm:t-inv-premod}

	Let $ \C $ be a pre-modular tensor category, let $ \D $ be a pivotal monoidal category and let $ \M\colon \C \to \D $ be a pivotal monoidal functor. Then $ \TM $ is T-invariant.

	\proof

	Let $ I,J \in \Irr(\C) $ be such that $ \TM_I^J \neq 0 $. Then, by Theorem~\ref{thm:alpha_induction_equivalence}, there exists a non-zero map $ \beta \in \Hom_{\D}(\M(I),\M(J^\vee)) $ that satisfies~\eqref{eq:graphical_alpha_condition}. We have
		\begin{align*}
		\Tc_{II} \beta =
			\begin{array}{c}
				\begin{tikzpicture}[scale=0.15,every node/.style={inner sep=0,outer sep=-1}]
				\node (v3) at (4,-10) {};
				\draw [thick, blue] (v3) edge (4,-4);
				\node (v6) at (4,-4) {};
				\node [blue] (v10) at (4.5,2) {$ I $};
				\draw[thick, blue] (6,-1) node (v1) {} to[out=-90,in=90] (4,-4);
				\draw[white,line width = 0.2cm] (4,-1) to[out=-90,in=90] (6,-4);
				\draw[thick, blue] (4,-1) node (v9) {} to[out=-90,in=90] (6,-4) node (v7) {};
				\node [blue] (v11) at (4,-12) {$ J^\vee $};
				\node (v5) at (8,-4) {};
				\node (v4) at (8,-1) {};
				\draw [thick, blue] (v1) to[out=90,in=90] (v4);
				\draw [thick, blue] (v4) edge (v5);
				\draw [thick, blue] (v7) to[out=-90,in=-90] (v5);
				\node (v8) at (4,0) {};
				\draw [thick, blue] (v8) edge (v9);
				\node [draw,outer sep=0,inner sep=2,minimum size=10,fill=white] at (4,-7) {$ \beta  $};
				\end{tikzpicture}
			\end{array}
		=
			\begin{array}{c}
				\begin{tikzpicture}[scale=0.15,every node/.style={inner sep=0,outer sep=-1}]
				\node (v3) at (4,-10) {};
				\draw [thick, blue] (v3) edge (4,-7);
				\node (v6) at (4,-7) {};
				\node [blue] (v10) at (4.5,2) {$ I $};
				\draw[thick, blue] (4,-2) node (v9) {} to[out=-90,in=90] (6,-7) node (v7) {};
				\draw[white,line width = 0.2cm] (6,-2) node (v1) {} to[out=-90,in=90] (4,-7);
				\draw[thick, blue] (6,-2) node (v1) {} to[out=-90,in=90] (4,-7);
				\node [blue] (v11) at (4,-12) {$ J^\vee $};
				\node (v5) at (8,-7) {};
				\node (v4) at (8,-2) {};
				\draw [thick, blue] (v1) to[out=90,in=90] (v4);
				\draw [thick, blue] (v4) edge (v5);
				\draw [thick, blue] (v7) to[out=-90,in=-90] (v5);
				\node (v8) at (4,0) {};
				\draw [thick, blue] (v8) edge (v9);
				\node [draw,outer sep=0,inner sep=2,minimum size=10,fill=white,rotate=180] at (8,-4.5) {$ \beta  $};
				\end{tikzpicture}
			\end{array}
		=
			\begin{array}{c}
				\begin{tikzpicture}[scale=0.15,every node/.style={inner sep=0,outer sep=-1}]
				\node (v3) at (4,-10) {};
				\draw [thick, blue] (v3) edge (4,-9);
				\node (v6) at (4,-9) {};
				\node [blue] (v10) at (4.5,2) {$ I $};
				\draw[thick, blue] (6,-6) node (v1) {} to[out=-90,in=90] (4,-9);
				\draw[white,line width = 0.2cm] (4,-6) to[out=-90,in=90] (6,-9);
				\draw[thick, blue] (4,-6) node (v9) {} to[out=-90,in=90] (6,-9) node (v7) {};
				\node [blue] (v11) at (4,-12) {$ J^\vee $};
				\node (v5) at (8,-9) {};
				\node (v4) at (8,-6) {};
				\draw [thick, blue] (v1) to[out=90,in=90] (v4);
				\draw [thick, blue] (v4) edge (v5);
				\draw [thick, blue] (v7) to[out=-90,in=-90] (v5);
				\node (v8) at (4,0) {};
				\draw [thick, blue] (v8) edge (v9);
				\node [draw,outer sep=0,inner sep=2,minimum size=10,fill=white] at (4,-3) {$ \beta  $};
				\end{tikzpicture}
			\end{array}
		= \Tc_{JJ} \beta
		\end{align*}
	where $ \Tc $ denotes the T-matrix and to make certain string manipulations clearer, we have chosen to write $ \beta $ upside-down instead of writing $ \beta^\vee $. Therefore $ Z(\TM)_{IJ} \neq 0 $ implies $ \Tc_{II} = \Tc_{JJ} $. As $ \Tc $ is diagonal that is precisely the condition that $ Z(F) $ commutes with $ \Tc $. \endproof

	\end{THM}

Our goal is therefore to reinterpret a module category $ \M \colon \C \to \End(\B) $ as a pivotal monoidal functor. Initially, this may seem improbable as $ \End(\B) $ admits a canonical pivotal structure and, excluding pathological examples, modules categories fail to be pivotal. However, we can study many interesting examples if we only require that $ \M $ \emph{induce a pivotal structure on its full image}. Let $ \D $ be the full image of $ \M $ in $ \End(\B) $, i.e.\ objects in $ \D $ are of the form $ \M(X) $ for $ X $ in $ \C $ and $$  \Hom_{\D}(\M(X), \M(Y)) = \Hom_{\End(\B)}(\M(X), \M(Y)).  $$ Clearly $ \D $ is a rigid monoidal category.
Furthermore, it comes with a natural candidate pivotal structure: $ \M(\delta_X) $, where $ \delta_X \colon ^\vee X \to X^\vee $ gives the pivotal structure on $ \C $. As $ \M $ is a functor, $ \M(\delta_X) $ is natural with respect to morphisms in $ \C $; however, to give a pivotal structure on $ \D $ it must be natural with respect \emph{all} morphisms in $ \D $. In other words the diagram
	\begin{equation}
		\begin{tikzcd}[row sep=2cm,column sep=2cm,inner sep=1ex] \label{eq:new_condition}
		\M(^\vee Y) \arrow[swap]{d}[name=D]{\M(\delta_Y)}   \arrow{r}{^\vee\alpha} &  \M(^\vee X) \arrow{d}[name=U]{\M(\delta_X)}
		\\
		\M(Y^\vee) \arrow{r}{\alpha^\vee} & \M(X^\vee)
		\arrow[to path={(U) node[scale=2,midway] {$\circlearrowleft$}  (D)}]{}
		\end{tikzcd}
	\end{equation}
must commute for all $ \alpha \in \Hom_{\D}(\M(X),\M(Y)) $. When this is satisfied and $ \D $ is equipped with the resulting pivotal structure, the functor $ \M \colon \C \to \D $ is automatically pivotal. We may therefore construct $ \TM $ and Theorem~\ref{thm:alpha_induction_equivalence} guarantees that $ Z(\TM) $ will give the same matrix as $ \alpha $-induction. Furthermore, the inclusion $ \D \hookrightarrow \End(\B) $ fully embeds $ \D $ into a multifusion category. We therefore obtain the following corollary of Theorem~\ref{thm:tm_is_frob} and Corollary~\ref{cor:tm_mod_inv}.

	\begin{COR} \label{cor:module_cats_give_frob_alg}

	Let $ \C $ be an MTC and let $ \M \colon \C \to \End(\B) $ be an indecomposable semisimple module category over $ \C $ with finitely many simple objects that induces a pivotal structure on its full image. Then $ \TM $ is a haploid, symmetric, commutative Frobenius algebra. In particular, $ Z(\TM) $ is a modular invariant if and only $ \sum_{IJ} \dim(\TM_I^J) d(I) d(J) = d(\C) $.

	\end{COR}

	\begin{REM}

	Ostrik's Example~\ref{ex:ostrik_alpha_fail} also shows that the condition that $ \M $ be  pivotal is necessary for the results of Section~\ref{sec:t_inv} and Section~\ref{sec:tm_as_a_frob_alg}. Indeed, one may check that his example fails to induce a pivotal structure on its full image.

	\end{REM}

\section{A Case Study: The Temperley-Lieb Category} \label{sec:case_study}

\renewcommand{\fld}{\mathbb{C}}

The goal of this section is to describe a class of interesting examples of module categories that induce a pivotal structure on their full images.

Let $ \beta $ be in $ \fld^* $. The Temperley-Lieb category $ \tl(\beta) $ is a $ \fld $-linear category whose set of objects is given by $ \{ \underline{n} \}_{n \in \mathbb{N}} $ where $ \underline{n} $ may be thought of as a collection of $ n $ dots along an interval. The space $ \Hom(\underline{m},\underline{n}) $ is the span of planar $ (m,n) $-tangles modulo the relation $ t \sqcup u - \beta t $ where t is a planar $ (m,n) $-tangle and $ u $  is the unknot. Composition is then given by tangle composition. For a more detailed description, see, for example,~\cite{Wenzl02}. $ \tl(\beta) $ is a monoidal category whose tensor product satisfies $ \underline{n} \otimes \underline{m} = \underline{n+m} $. Furthermore, $ \tl(\beta) $ is rigid and every object admits a canonical choice of self dualizing maps
	\begin{align*}
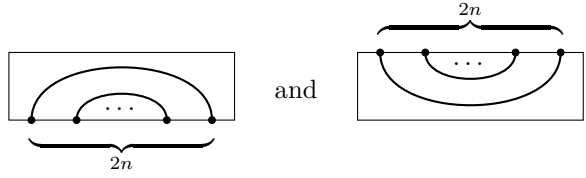
		\begin{array}{c}
			\begin{tikzpicture}[scale=0.3, yscale=-1]
			\draw[fill] (-2,2) circle (0.15);
			\draw[fill] (0,2) circle (0.15);
			\draw[fill] (6,2) circle (0.15);
			\draw[fill] (4,2) circle (0.15);
			\draw[thick] (-2,2) to[out=-90, in=-90] (6,2);
			\draw[thick] (0,2) to[out=-90, in=-90] (4,2);
			\draw (-3,2) rectangle (7,-1);
			\node at (2,1.5) {\dots};
			\node at (2,3.5) {$ \underbrace{\hspace{7em}}_{2n} $};
			\node at (2,-2.5) {\vphantom{$ \overbrace{\hspace{7em}}^{2n} $}};
			\end{tikzpicture}
		\end{array}
	\quad \text{and} \quad
		\begin{array}{c}
			\begin{tikzpicture}[scale=0.3]
			\draw[fill] (-2,2) circle (0.15);
			\draw[fill] (0,2) circle (0.15);
			\draw[fill] (6,2) circle (0.15);
			\draw[fill] (4,2) circle (0.15);
			\draw[thick] (-2,2) to[out=-90, in=-90] (6,2);
			\draw[thick] (0,2) to[out=-90, in=-90] (4,2);
			\draw (-3,2) rectangle (7,-1);
			\node at (2,1.5) {\dots};
			\node at (2,3.5) {$ \overbrace{\hspace{7em}}^{2n} $};
			\node at (2,-2.5) {\vphantom{$ \underbrace{\hspace{7em}}_{2n} $}};
			\end{tikzpicture}
		\end{array}
	\end{align*}
denoted $ \cre_n \colon \underline {0} \to \underline{n} \otimes \underline{n} $ and $ \ann_n \colon \underline{n} \otimes \underline{n} \to \underline {0} $ respectively. In the case when $ \beta = -[2]_q $ for a primitive even root of unity $ q $, $ \tl(\beta) $ admits a unique tensor ideal $ \N $~\cite{Wenzl02}. Quotienting $ \tl(\beta) $ by $ \N $ and idempotent completing the result yields a spherical fusion category we denote $ \C $. Using a skein relation we may define a non-degenerate braiding on $ \C $ giving us a modular tensor category~\cite[Theorem 7.5.3]{Turaev16}. Let $ h $ be the the smallest positive integer such that $ q^{2h} = 1 $ or equivalently $ [h]_q = 0 $; $ h $ is called the \emph{Coxeter number} of $ \C $. $ \C $ turns out to be equivalent to the category of integrable highest weight modules of $ A_1^{(1)} $ at level $ k=h-2 $, denoted $ \Rep_k A_1^{(1)} $ (for this equivalence to be pivotal one must equip $ \tl(\beta) $ with a `twisted' pivotal structure, or alternatively, consider the so-called ``disoriented" diagrammatic category presented in \cite[p. 5]{Morrison09}; for further details on this issue see \cite{Snyder09}). In particular, a complete set of simples in $ \C $ has size $ h-1 $.

Let $ \Q $ be a symmetric quiver with non-degenerate eigenvalue $ \beta $ (here non-degenerate signifies there exists an eigenvector $ (x_i) = x $  with non-zero entries), let $ A $ be the basic algebra spanned by vertices in $ \Q $ and let $ B $ be the $ A $-bimodule spanned by arrows in $ \Q $. We can construct a module category over $ \tl(\beta) $ as follows
    \begin{align} \label{eq:define_functor}
        \begin{split}
        \M \colon \tl(\beta) &\to \BMod{A} = \End(\Mod{A})\\
        \M(\underline{n}) &= B^{\otimes n} \\
        \M(\ann_1)_{ij} &= \left( \phi_{ij} \colon v \otimes w \mapsto x_j \langle v, w \rangle \right)\\
        \M(\cre_1)_{ij} &= \Big( \varphi_{ij} \colon 1 \mapsto x_i^{-1} \sum_b b \otimes b^* \Big).    \end{split}
    \end{align}
As $ \cre_1 $ and $ \ann_1 $ tensor generate $ \tl(\beta) $ this fully determines $ \M $. One advantage of considering module categories of this form is that they induce a pivotal structure on the full image (cf. Section~\ref{sec:mod_cat_and_alpha_ind}). To prove this we first consider the following lemma.

	\begin{LEMMA}

	Let $ \phi^n $ and $ \varphi^n $ denote the image of $ \M(\ann_{n}) $ and $ \M(\cre_{n}) $ respectively. For $ i,j \in Q_0 $ and $ n \in \mathbb{N}^+ $ we have the following
		\begin{align} \label{eq:first_pivotal_trick}
		\phi_{ji}^n(w \otimes v) = \frac{x_i}{x_j} \ \phi_{ij}^n(v \otimes w)
		\end{align}
	and
		\begin{align} \label{eq:second_pivotal_trick}
		\varphi_{ji}^n = \frac{x_i}{x_j} \ T_{ij}^n \circ \varphi_{ij}^n
		\end{align}
	where $ T_{ij}^n $ is the canonical isomorphism from $ \Bmt{i}{j}{n} \otimes \Bmt{j}{i}{n} $ to $ \Bmt{j}{i}{n} \otimes \Bmt{i}{j}{n} $.

	\proof

	We proceed by induction on $ n $. The base case $ n = 1 $ is clear. Assuming the hypothesis for all integers up to $ n-1 $, we take $ b \in \Bm{i}{k}, \ v \in \Bmt{k}{j}{n-1}, \ w \in \Bmt{j}{k}{n-1} $ and compute,
		\begin{align*}
		\phi_{ij}^n(b \otimes v \otimes w \otimes b^*) &= \phi_{kj}^{n-1}(v \otimes w) \phi_{ik}(b \otimes b^*) \\
		&= x_k \phi_{kj}^{n-1}(v \otimes w).
		\end{align*}
	We then also have
		\begin{align*}
		\phi_{ji}^n(w \otimes b^* \otimes b \otimes v) &= \phi_{jk}^{n-1}(w \otimes v) \phi_{ki}(b^* \otimes b) \\
		&= x_i \phi_{jk}^{n-1}(w \otimes v) \\
		&= \frac{x_ix_k}{xj} \phi_{kj}^{n-1}(v \otimes w) \\
		&= \frac{x_i}{xj} \phi_{ij}^n(b \otimes v \otimes w \otimes b^*).
		\end{align*}
	Therefore~\eqref{eq:first_pivotal_trick} is proved. To prove \eqref{eq:second_pivotal_trick} we proceed more directly,
		\begin{align*}
		\frac{x_i}{x_j} \ T_{ij}^n \circ \varphi_{ij}^n &= \frac{x_i}{x_j} \ T_{ij}^n \circ \Big( \sum\limits_k (\id_{\underline{n-1}} \otimes \ \varphi_{kj} \otimes \id_{\underline{n-1}}) \circ \varphi^{n-1}_{ik} \Big) \\
		&= \frac{x_i}{x_j} \sum\limits_k (\id_{\underline{1}} \otimes \ ( T^{n-1}_{ik} \circ \varphi^{n-1}_{ik} ) \otimes \id_{\underline{1}}) \circ (T^1_{kj} \circ \varphi_{kj}) \\
		&= \sum\limits_k (\id_{\underline{1}} \otimes \ \varphi^{n-1}_{ki} \otimes \id_{\underline{1}}) \circ \varphi_{jk} \\
		&= \varphi^n_{ji}
		\end{align*}
	as desired. \endproof

	\end{LEMMA}

	\begin{PROP} \label{prop:induce_pivotal_structure}

	Let $ \M $ be a module category over $ \tl(\beta) $ given by \eqref{eq:define_functor}. Then~\eqref{eq:new_condition} commutes. In other words, $ \M $ induces a pivotal structure on its full image.

	\proof

	As the pivotal structure on $ \tl(\beta) $ is given by the identity~\eqref{eq:new_condition} reduces to $ ^\vee \alpha = \alpha^\vee $ for all $ \alpha \in \Hom_{\D}(\M(\underline{m}),\M(\underline{n})) $. For $ a \in \Bmt{i}{j}{n} $, by~\eqref{eq:first_pivotal_trick}, we have
		\begin{align*}
		\alpha^\vee(a) &=(\id \otimes \phi_{ji}^m) \circ (\id \otimes \alpha_{ji} \otimes \id )\circ (\varphi_{ij}^n \otimes \id)(a) \\
		&= \sum\limits_{IJ} \lambda_{IJ}^n \phi^m_{ji}(\alpha_{ji}(b_J) \otimes a) \ b_I \\
		&= \frac{x_i}{x_j} \sum\limits_{IJ} \lambda_{IJ}^n \phi^m_{ij}(a \otimes \alpha_{ji}(b_J)) \ b_I
		\end{align*}
	where the $ \lambda_{IJ}^n $, the $ b_I $ and the $ b_J $ are such that
		\begin{align*}
		\varphi_{ij}^n (1) = \sum\limits_{IJ} \lambda_{IJ}^n \ b_I \otimes b_J \in \Bmt{i}{j}{n} \otimes \Bmt{j}{i}{n}.
		\end{align*}
	However, by~\eqref{eq:second_pivotal_trick}, we also have
		\begin{align*}
		^\vee\alpha(a) &=(\phi_{ij}^m \otimes \id ) \circ (\id \otimes \alpha_{ji} \otimes \id )\circ (\id \otimes \varphi_{ji}^n)(a) \\
		&= \frac{x_i}{x_j} \sum\limits_{IJ} \lambda_{IJ}^n \phi^m_{ij}(a \otimes \alpha_{ji}(b_J)) \ b_I \\
		&= \alpha^\vee(a)
		\end{align*}
	and so we are done.	\endproof

	\end{PROP}

Once again, let $ \Q $ be a symmetric quiver with non-degenerate eigenvalue $ \beta $ and let $ \M $ be given by~\eqref{eq:define_functor}. Under certain additional conditions on $ \Q $, $ \M $ will vanish on $ \N $ and give a module category over $ \C $. Such modules categories turn out to classify all modules categories over $ \C $, as described by the following theorem.

	\begin{THM}[\cite{MR2046203}, Theorem 3.12] \label{thm:ost_class_mod_over_tl}

	Indecomposable semisimple module categories over $ \C $ with finitely many simple objects are classified by the double Dynkin quivers of type A, D, E.

	\end{THM}

	\begin{COR} \label{cor:all_mod_cat_over_tl_are_pivotal}

	Every semisimple module category over $ \C $ with finitely many simple objects induces a pivotal structure on its full image.

	\proof

	This follows immediately from Theorem~\ref{thm:ost_class_mod_over_tl} and Proposition~\ref{prop:induce_pivotal_structure}. \endproof

	\end{COR}

A module category $ \M\colon \C \to \BMod{A} $ over an arbitrary monoidal category $ \C $ comes equipped with a natural action of $ \K_{\fld}(\C) $ on $ \K_{\fld}(\LMod{A}) $ given by $ [X] \cdot [V] = [\M(X) \otimes_A V] $. However, when $ \C $ is a spherical fusion category and $ \M $ induces a pivotal structure on its full image, we may consider the $ \TM $ construction. As $ \End_{\TC}(\tid) = \K_{\fld}(\C) $ this defines \emph{another} action of $ \K_{\fld}(\C) $ on $ \TM(\tid) = \End_D(A) $ (where $ \D $ is the full image of $ \M $). Exploiting graphical calculus in $ \D $, this action is given by
   \begin{align*}
   [X] \cdot \alpha =
	   \begin{array}{c}
		   \begin{tikzpicture}[scale=0.15,every node/.style={inner sep=0,outer sep=-1}]
		   \node [blue] at (5.5,0) {$ X $};
		   \draw [line width= 0.15cm,white] (0,0) ellipse (4 and 4);
		   \draw [thick, blue] (0,0) ellipse (4 and 4);
		   \node [draw,outer sep=0,inner sep=1,minimum height=12,minimum width=24,fill=white] (v5) at (0,0) {$  \alpha $};
		   \end{tikzpicture}
	   \end{array}.
   \end{align*}
In the case when $ \C $ is the semisimple quotient category constructed from $ \tl(\beta) $, these two actions coincide.

   \begin{PROP} \label{prop:iso_of_modules_over_kc}

    Let $ \M $ be a module category over $ \C $. W.l.o.g.\ we suppose that $ \M $ is given by \eqref{eq:define_functor}. For $ j \in \Q_0 $, let $ \tid_j $ be the corresponding idempotent on $ A $ and let $ V_i $ be the corresponding simple $ A $-module. Then the map
	   \begin{align*}
	   \Phi \colon \End_{\D}(A) &\to \K_{\fld}(\LMod{A}) \\
	   \tid_{j} &\mapsto x_j [V_j]
	   \end{align*}
   is an isomorphism of $ \K_{\fld}(\C) $-modules.

   \proof

   As $ \{[V_j]\} $ and $ \{\tid_j\} $ form a basis of $ \K_{\fld}(\LMod{A}) $ and $ \End_{\D}(A) $ respectively, $ \Phi $ is an isomorphism of vector spaces. As, for $ X $ in $ \C $,
	   \begin{align*}
	   [X] \cdot \tid_j = \sum\limits_i \phi_{ij}^X \circ \varphi_{ij}^X \tid_i = \sum\limits_i \frac{x_j}{x_i} \dim \MX{i}{j}{X}\tid_i
	   \end{align*}
   where $ \phi_{ij}^X $ and $ \varphi_{ij}^X $ is $ \M(\ann_X) $ and $ \M(\cre_X) $ respectively, we have
	   \begin{center}
		   \begin{tikzcd}[row sep=2cm,column sep=2cm,inner sep=1ex]
		   \tid_j \arrow[mapsto,swap]{d}[name=D]{\Phi}   \arrow[mapsto]{r}{[X]} &  \sum\limits_i \frac{x_j}{x_i} \dim \MX{i}{j}{X}\tid_i
		   \\
		   x_j [V_j] \arrow[mapsto]{r}{[X]} & x_j \sum\limits_{i} \dim \MX{i}{j}{X} \ [V_i]     \arrow[mapsto,swap]{u}[name=U]{\Phi^{-1}}
		   \arrow[to path={(U) node[scale=2] at (-0.6,0) {$\circlearrowleft$}  (D)}]{}
		   \end{tikzcd}
	   \end{center}
   as desired.   \endproof

   \end{PROP}

    \begin{REM}
    Let $ \M $ be a module category over $ \C $ arising from some quiver $ \Q $ via~\eqref{eq:define_functor} and let $ n $ be in $ \N $. Then $ \TM(\underline{n}) $ is the space of \emph{cycles} of length $ n $ in $ \Q $.
    \end{REM}

We are now ready to exploit the $ \TM $ construction to explain a well known pattern in the classification of modular invariants over $ \C $.

	\begin{THM}[C.I.Z.\ classification~\cite{MR918402}] \label{thm:ciz}

	The complete list of modular invariants over $ \C $ is as follows. To aid legibility, we present these modular invariants as partition functions, cf.~\eqref{eq:partition_function}.
		\begin{align*}
		{\A}_{h-1}=&\,\sum_{a=1}^{h-1}\,|\chi_a|^2\ ,\qquad\qquad
		\qquad & &\forall h\ge 3\\
		{\D}_{{h\over 2}+1}=&\,\sum_{a=1}^{h-1}\,\chi_a\,\chi_{J^{a-1}a}^*\ ,\qquad
		\qquad\qquad & &\text{ whenever}\ {h\over 2}\ \text{ is even}\\
		{\D}_{{h\over 2}+1}=&\,|\chi_1+\chi_{J1}|^2+|\chi_3+\chi_{J3}|^2+\cdots
		+2|\chi_{{h\over 2}}|^2\ ,\qquad & &\text{ whenever}\ {h\over 2}\ \text{ is\ odd}\\
		{\E}_6=&\,|\chi_1+\chi_7|^2+|\chi_4+\chi_8|^2+|\chi_5+\chi_{11}|^2\ ,
		\qquad\qquad & &\text{ for}\ h=12\\
		{\E}_7=&\,|\chi_1+\chi_{17}|^2+|\chi_5+\chi_{13}|^2+|\chi_7+\chi_{11}|^2
		&\\&\,+\chi_9\,(\chi_3+\chi_{15})^*+(\chi_3+\chi_{15})\,\chi^*_9+|\chi_9|^2
		\ ,\qquad\qquad & &\text{ for}\ h=18\\
		{\E}_8=&\,|\chi_1+\chi_{11}+\chi_{19}+\chi_{29}|^2+|\chi_7+\chi_{13}+
		\chi_{17}+\chi_{23}|^2\ ,\qquad & &\text{ for}\ h=30.
		\end{align*}
	where $ h $ is the Coxeter number of $ \C $ and $ J : \{ 1 , 2 , ... , h-1 \} \to \{ 1 , 2 , ... , h-1 \} $ maps $ a $ to $ h-a $.

	\end{THM}

As alluded to in the introduction, the classification of modular invariants over $ \C $ admits the following A-D-E pattern. Let $ \X $ be a double Dynkin quiver of type A,D or E. The eigenvalues of $ \X $ form a subset of $ \{ -[2]_q \mid q = e^{\frac{\pi i l}{h}}, 1 \leq l \leq h-1 \} $ for some $ h \in \mathbb{N} $. Then, for $ l \in \{ 1, 2, \dots h-1 \} $, the $ l $th diagonal entry in the modular invariant associated to $ \X $ gives the dimension of the corresponding eigenspace of $ \X $.

Let $ \X $ be an A-D-E double Dynkin quiver and let $ \M\colon \C \to \BMod{A} $ be the corresponding module category over $ \C $. It is known that applying $ \alpha $-induction, as described in Section~\ref{sec:mod_cat_and_alpha_ind}, to $ \M $ yields the modular invariant associated to $ \X $ by the list appearing in Theorem~\ref{thm:ciz}~\cite[Section 5]{MR1867545}. However, even once this connection has been established it is non-trivial to explain the A-D-E pattern described above. The $ \TM $ construction explains this pattern in the following way. Let $ Z $ be the modular invariant obtained by applying $ \alpha $-induction to $ \M $. By Theorem~\ref{thm:alpha_induction_equivalence} the entries of $ Z $ may be thought of as the dimensions of the simple multiplicity spaces in $ \TM $, in other words
 	\begin{align*}
 	Z = Z(\TM)
 	\end{align*}
where $ Z(\TM) $ is given by Definition~\ref{def:zf}.

 We recall that $ \End_{\TC}(\tid) $ is a semisimple commutative algebra generated by the orthogonal primitive idempotents $ \{ \tid_I \}_{I \in \I} $ where $ (\tid, \tid_I)\Sharp = (II^\vee, e_I^{I\vee})\Sharp $, see Remark~\ref{rem:simples_in_kc}.
The diagonal terms in $ Z $ therefore correspond to the dimensions of the weight spaces of the action of $ \End_{\TC}(\tid) = \K_{\fld}(\C) $ on $ TM(\tid) $. However by Proposition~\ref{prop:iso_of_modules_over_kc} this action coincides with the natural action of $ \K_{\fld}(\C) $ on $ \K_{\fld}(\LMod{A}) $. As the weight spaces of this action are given by the eigenspaces of $ \X $ this explains the pattern.

\appendix

\section{Lemmas on Frobenius Algebras} \label{sec:appendix}

The purpose of this section is to provide the necessary definitions on Frobenius algebras followed by certain technical results required in Section~\ref{sec:tm_as_a_frob_alg}.

    \begin{DEF} \label{def:frob_alg}

    Let $ \B $ be a monoidal category. A \emph{Frobenius algebra} $ A $ in $ \B $ is an algebra and a coalgebra in $ \B $ such that
        \begin{align}\label{eq:frob_condition}
        (\id_A \otimes \nabla) \circ (\Delta \otimes \id_A) = \Delta \circ \nabla  = (\nabla \otimes \id_A) \circ (\id_A \otimes \Delta)
        \end{align}
    where $ \nabla $ is the product and $ \Delta $ is the coproduct. Using the graphical notation
        \begin{align*}
        \nabla \ = \hspace{-1em}
            \begin{array}{c}
                \begin{tikzpicture}[scale=0.3]
        		\coordinate (v1) at (-0.5,0.3) {};
        		\coordinate (v2) at (0.5,0.3) {};
        		\coordinate (v3) at (0,-0.45) {};
        		\draw[thick] (-1,2) to[out=-90] ($ (v1) + (0.11,-0.1) $);
        		\draw[thick] (1,2) to[out=-90,in=45] ($ (v2) + (-0.11,-0.1) $);
        		\draw[thick] (0,0) -- (0,-2);
        		\fill [blue] (v1) -- (v2) -- (v3) -- cycle;
                \node at (-1,2.8) {$A$};
                \node at (1,2.8) {$A$};
                \node at (0,-3) {$A$};
        		\end{tikzpicture}
            \end{array}
        \quad \text{and} \ \ \quad \Delta  \ = \hspace{-1em}
            \begin{array}{c}
                \begin{tikzpicture}[scale=0.3, yscale = -1]
                \coordinate (v1) at (-0.5,0.3) {};
                \coordinate (v2) at (0.5,0.3) {};
                \coordinate (v3) at (0,-0.45) {};
                \draw[thick] (-1,2) to[out=-90] ($ (v1) + (0.11,-0.1) $);
                \draw[thick] (1,2) to[out=-90,in=45] ($ (v2) + (-0.11,-0.1) $);
                \draw[thick] (0,0) -- (0,-2);
                \fill [blue] (v1) -- (v2) -- (v3) -- cycle;
                \node at (-1,2.8) {$A$};
                \node at (1,2.8) {$A$};
                \node at (0,-3) {$A$};
                \end{tikzpicture}
            \end{array}
        \end{align*}
    we can rewrite Condition \eqref{eq:frob_condition} as
        \begin{align} \label{eq:grapical_frob_condition}
            \begin{array}{c}
                \begin{tikzpicture}[scale=0.2]
                \coordinate (v1) at (-0.5,0.3) {};
                \coordinate (v2) at (0.5,0.3) {};
                \coordinate (v3) at (0,-0.45) {};
                \draw[thick] (-3,2) -- (-3,-2);
                \draw[thick] (-1,2) to[out=-90] ($ (v1) + (0.11,-0.1) $);
                \draw[thick] (1,2) to[out=-90,in=45] ($ (v2) + (-0.11,-0.1) $);
                \draw[thick] (0,0) -- (0,-2);
                \fill [blue] (v1) -- (v2) -- (v3) -- cycle;
                \node at (0,-3.2) {$A$};
                \node at (-3,-3.2) {$A$};
                \coordinate (v4) at (-2.5,3.7) {};
                \coordinate (v5) at (-1.5,3.7) {};
                \coordinate (v6) at (-2,4.45) {};
                \draw[thick] (1,2) -- (1,6);
                \draw[thick] (-3,2) to[out=90,in=-135] ($ (v4) + (0.11,+0.1) $);
                \draw[thick] (-1,2) to[out=90,in=-45] ($ (v5) + (-0.11,+0.1) $);
                \draw[thick] (-2,3.7) -- (-2,6);
                \fill [blue] (v4) -- (v5) -- (v6) -- cycle;
                \node at (1,7.2) {$A$};
                \node at (-2,7.2) {$A$};
                \end{tikzpicture}
            \end{array}
        =
            \begin{array}{c}
        		\begin{tikzpicture}[scale=0.2]
        		\coordinate (v1) at (-0.5,0.3) {};
        		\coordinate (v2) at (0.5,0.3) {};
        		\coordinate (v3) at (0,-0.45) {};
        		\coordinate (v4) at (0,-2) {};
        		\draw[thick] (-1,2.75) to[out=-90] ($ (v1) + (0.11,-0.1) $);
        		\draw[thick] (1,2.75) to[out=-90,in=45] ($ (v2) + (-0.11,-0.1) $);
        		\draw[thick] (0,0) -- ($(v4) + (0,-0.75) $);
        		\fill [blue] (v1) -- (v2) -- (v3) -- cycle;
                \node at (-1,3.95) {$A$};
                \node at (1,3.95) {$A$};
        		\coordinate (v5) at (-0.5,-2.75) {};
        		\coordinate (v6) at (0.5,-2.75) {};
            	\draw[thick] (-1,-5.25) to[out=90, in = -135] ($ (v5) + (0.11,0.1) $);
            	\draw[thick] (1,-5.25) to[out=90,in=-45] ($ (v6) + (-0.11,0.1) $);
            	\fill [blue] (v4) -- (v5) -- (v6) -- cycle;
                \node at (-1,-6.45) {$A$};
                \node at (1,-6.45) {$A$};
        		\end{tikzpicture}
            \end{array}
        =
            \begin{array}{c}
                \begin{tikzpicture}[scale=0.2, xscale= -1]
                \coordinate (v1) at (-0.5,0.3) {};
                \coordinate (v2) at (0.5,0.3) {};
                \coordinate (v3) at (0,-0.45) {};
                \draw[thick] (-3,2) -- (-3,-2);
                \draw[thick] (-1,2) to[out=-90] ($ (v1) + (0.11,-0.1) $);
                \draw[thick] (1,2) to[out=-90,in=45] ($ (v2) + (-0.11,-0.1) $);
                \draw[thick] (0,0) -- (0,-2);
                \fill [blue] (v1) -- (v2) -- (v3) -- cycle;
                \node at (0,-3.2) {$A$};
                \node at (-3,-3.2) {$A$};
                \coordinate (v4) at (-2.5,3.7) {};
                \coordinate (v5) at (-1.5,3.7) {};
                \coordinate (v6) at (-2,4.45) {};
                \draw[thick] (1,2) -- (1,6);
                \draw[thick] (-3,2) to[out=90,in=-135] ($ (v4) + (0.11,+0.1) $);
                \draw[thick] (-1,2) to[out=90,in=-45] ($ (v5) + (-0.11,+0.1) $);
                \draw[thick] (-2,3.7) -- (-2,6);
                \fill [blue] (v4) -- (v5) -- (v6) -- cycle;
                \node at (1,7.2) {$A$};
                \node at (-2,7.2) {$A$};
                \end{tikzpicture}
            \end{array}.
        \end{align}
    We also use
        \begin{tikzpicture}[scale=0.3]
        \draw[thick] (0,1) -- (0,0);
        \draw[thick, fill=white] (0,1) circle (0.2);
        \end{tikzpicture}
    to denote the \emph{unit} and
        \begin{tikzpicture}[scale=0.3, yscale=-1]
        \draw[thick] (0,1) -- (0,0);
        \draw[thick, fill=white] (0,1) circle (0.2);
        \end{tikzpicture}
    to denote the \emph{counit}. $ A $ is called \emph{haploid} if is satisfies $ \Hom_{\B}(\tid, A) = \fld  $. If $ \B $ is braided then $ A $ is called \emph{commutative} if the underlying algebra structure is commutative i.e.\
        \begin{align} \label{eq:commutative_condition}
            \begin{array}{c}
                \begin{tikzpicture}[scale=0.3,every node/.style={inner sep=0,outer sep=-1}]
                \coordinate (v1) at (-0.5,0.3) {};
                \coordinate (v2) at (0.5,0.3) {};
                \coordinate (v3) at (0,-0.45) {};
                \draw[thick] (-0.75,1) node (v4) {} to[out=-90] ($ (v1) + (0.11,-0.1) $);
                \draw[thick] (0.75,1) node (v6) {} to[out=-90,in=45] ($ (v2) + (-0.11,-0.1) $);
                \draw[thick] (0,0) -- (0,-2);
                \fill [blue] (v1) -- (v2) -- (v3) -- cycle;
                \node at (-1,4) {$A$};
                \node at (1,4) {$A$};
                \node at (0,-3) {$A$};
                \node (v7) at (-1,3) {};
                \node (v5) at (1,3) {};
                \draw [thick] (v4) to[out=90,in=-90] (v5);
                \draw [line width= 0.2cm,white] (v6) to[out=90,in=-90] (v7);
                \draw [thick] (v6) to[out=90,in=-90] (v7);
                \end{tikzpicture}
            \end{array}
        = \nabla.
        \end{align}
    If $ \B $ is pivotal then $ A $ is called symmetric if is satisfies
        \begin{align} \label{eq:symmetric_condition}
            \begin{array}{c}
                \begin{tikzpicture}[scale=0.3,every node/.style={inner sep=0,outer sep=-1}]
                \coordinate (v1) at (-0.5,0.3) {};
                \coordinate (v2) at (0.5,0.3) {};
                \coordinate (v3) at (0,-0.45) {};
                \draw[thick] (-1,2) to[out=-90] ($ (v1) + (0.11,-0.1) $);
                \draw[thick] (1,1.5) node (v4) {} to[out=-90,in=45] ($ (v2) + (-0.11,-0.1) $);
                \draw[thick] (0,0) -- (0,-2);
                \fill [blue] (v1) -- (v2) -- (v3) -- cycle;
                \node at (-1,2.8) {$A$};
                \node at (3,-3) {$A^\vee$};
                \node (v5) at (3,1.5) {};
                \draw [thick] (v4) to[out=90,in=90] (v5);
                \node (v6) at (3,-2) {};
                \draw [thick] (v5) edge (v6);
                \draw[thick, fill=white] (0,-2) circle (0.2);
                \end{tikzpicture}
            \end{array}
        =
            \begin{array}{c}
                \begin{tikzpicture}[scale=0.3,every node/.style={inner sep=0,outer sep=-1},xscale=-1]
                \coordinate (v1) at (-0.5,0.3) {};
                \coordinate (v2) at (0.5,0.3) {};
                \coordinate (v3) at (0,-0.45) {};
                \draw[thick] (-1,2) to[out=-90] ($ (v1) + (0.11,-0.1) $);
                \draw[thick] (1,1.5) node (v4) {} to[out=-90,in=45] ($ (v2) + (-0.11,-0.1) $);
                \draw[thick] (0,0) -- (0,-2);
                \fill [blue] (v1) -- (v2) -- (v3) -- cycle;
                \node at (-1,2.8) {$A$};
                \node at (3,-3) {$A^\vee$};
                \node (v5) at (3,1.5) {};
                \draw [thick] (v4) to[out=90,in=90] (v5);
                \node (v6) at (3,-2) {};
                \draw [thick] (v5) edge (v6);
                \draw[thick, fill=white] (0,-2) circle (0.2);
                \end{tikzpicture}
            \end{array}.
        \end{align}
    \end{DEF}

    \begin{REM} \label{rem:decompose_morphisms}

    Let $ \B $ be a fusion category with complete set of simples $ \Irr(\B) $ and let $ A $ be an object in $ \B $. Any morphism $ \nabla $ from $ A \otimes A $ to $ A $ gives rise to the following morphisms,
        \begin{align*}
        \nabla_{X}^{Y,Z} \colon \Hom_{\B}(X,YZ) &\to \Hom(A_Y \otimes A_Z, A_X)\\
        \alpha &\mapsto \left(g \otimes h \mapsto
            \begin{array}{c}
    			\begin{tikzpicture}[scale=0.25,every node/.style={inner sep=0,outer sep=-1}]
    			\draw[thick] (-0.8,0) -- (-0.8,-2);
    			\draw[thick] (0.8,0) -- (0.8,-2);
    			\draw[thick] (0,0) -- (0,2);
    			\node [draw,outer sep=0,inner sep=2,minimum width=15,minimum height= 8,fill=white] at (0,0.4) {$\scriptstyle \alpha $};
    			\node at (0.5,1.7) {$\scriptscriptstyle X $};
    			\coordinate (v1) at (-0.5,-3) {} {} {} {} {};
    			\coordinate (v2) at (0.5,-3) {} {} {} {} {};
    			\coordinate (v3) at (0,-3.75) {} {} {};
    			\draw[thick] (-0.8,-2) to[out=-90] ($ (v1) + (0.11,-0.1) $);
    			\draw[thick] (0.8,-2) to[out=-90,in=45] ($ (v2) + (-0.11,-0.1) $);
    			\draw[thick] (0,-3.3) -- (0,-5);
    			\fill [blue] (v1) -- (v2) -- (v3) -- cycle;
    			\node at (0.5,-4.7) {$ \scriptscriptstyle A$};
    			\node [draw,outer sep=0,inner sep=1,minimum height = 10,fill=white] at (-0.8,-1.5) {$\scriptstyle g $};
    			\node [draw,outer sep=0,inner sep=1,minimum width = 7,minimum height = 10,fill=white] at (0.8,-1.5) {$\scriptstyle h $};
    			\end{tikzpicture}
            \end{array}\right)
        \end{align*}
    where $ X,Y,Z $ are in $ \B $ and $ A_X \defeq \Hom_{\B}(X,A) $. The full map $ \nabla $ is determined by $ \nabla_R^{S,T} $ for $ R,S,T \in \Irr(\B) $. Indeed we can recover it via
	\begin{align*}
	\bigoplus\limits_{RST} \sum\limits_{g,h,\alpha} \nabla_R^{S,T} (\alpha)(g \otimes h) \circ \alpha^* \circ (g^* \otimes h^*)
	= \nabla
	\end{align*}
where $ g $ ranges over a basis of $ A_S $, $ h $ ranges over a basis of $ A_T $ and $ \alpha $ ranges over a basis of $ \Hom_{\B}(R,ST) $. Similarly any morphism from $ A $ to $ A \otimes A $ can also be decomposed in the following way
    \begin{align*}
    \Delta_{ST}^R \colon \Hom_{\B}(ST,R) &\to \Hom(A_R, A_S \otimes A_T) \\
    \beta &\mapsto \left( f \mapsto \sum\limits_{g,h}
        \begin{array}{c}
			\begin{tikzpicture}[scale=0.25,every node/.style={inner sep=0,outer sep=-1}, xscale=1.1,yscale=-1]
			\draw[thick] (-0.8,0) -- (-0.8,-2);
			\draw[thick] (0.8,0) -- (0.8,-2);
			\draw[thick] (0,0) -- (0,2);
			\node [draw,outer sep=0,inner sep=2,minimum width=15,minimum height= 8,fill=white] at (0,0.4) {$\scriptstyle \beta $};
			\node at (0.5,1.7) {$\scriptscriptstyle R $};
			\coordinate (v1) at (-0.5,-2.75) {} {} {};
			\coordinate (v2) at (0.5,-2.75) {} {} {};
			\coordinate (v3) at (0,-3.625) {} {};
			\draw[thick] (-0.8,-2) to[out=-90] ($ (v1) + (0.11,-0.1) $);
			\draw[thick] (0.8,-2) to[out=-90,in=45] ($ (v2) + (-0.11,-0.1) $);
			\draw[thick] (0,-3) -- (0,-6.5);
			\fill [blue] (v1) -- (v2) -- (v3) -- cycle;
			\node at (0.5,-6.2) {$ \scriptscriptstyle R$};
			\node [draw,outer sep=0,inner sep=1,minimum height = 10,fill=white] at (-0.8,-1.5) {$\scriptstyle g^* $};
			\node [draw,outer sep=0,inner sep=1,minimum width = 7,minimum height = 10,fill=white] at (0.8,-1.5) {$\scriptstyle h^* $};
			\node [draw,outer sep=0,inner sep=1,minimum width = 7,minimum height = 9,fill=white] at (0,-5) {$\scriptstyle f $};
			\end{tikzpicture}
        \end{array}
    g \otimes h
    \right)
    \end{align*}
and then recovered via
	\begin{align*}
	\sum\limits_{\substack{RST \\ \beta, f}}  \Delta^R_{ST}(\beta)(f) \circ \beta^* \circ f^*
	= \Delta.
	\end{align*}
	\end{REM}

	\begin{LEMMA} \label{lem:build_an_algebra}

	Let $ A $ be an object in $ \B $ and let $ \nabla $ be in $ \Hom_{\B}(A \otimes A,A) $. Then $ \nabla $ is associative if
		\begin{align}\label{eq:condition_assoc}
		\nabla^{RS,T}_{RST} (\id) \big( \nabla^{R,S}_{RS} (\id) (f \otimes g) \otimes h \big) = \nabla^{R,ST}_{RST} (\id) \big(  f \otimes \nabla^{S,T}_{ST}(\id)(g \otimes h) \big)
		\end{align}
	for all $ R,S,T \in \Irr(\B) $, $ \alpha \in \Hom_{\B}(R,ST) $, $ f \in A_R $, $ g \in A_S $ and $ h \in A_T $. An element $ u \in A_{\tid} $ is a unit for $ \nabla $ if
		\begin{align} \label{eq:condition_unit}
		\nabla^{\tid,S}_S (\id) (u \otimes g) = g \quad \text{and} \quad \nabla^{S,\tid}_S (\id) (g \otimes u) = g
		\end{align}
	Furthermore, if $ \B $ is braided then $ \nabla $ is commutative if
		\begin{align} \label{eq:condition_commutative}
		\nabla^{T,S}_{ST} \big( \hspace{-.3em}
			\begin{array}{c}
				\begin{tikzpicture}[scale=0.25,every node/.style={inner sep=0,outer sep=-1}]
				\node (v4) at (-0.5,-1) {};
				\node (v6) at (1,-1) {};
				\node (v7) at (-0.5,0.5) {};
				\node (v5) at (1,0.5) {};
				\draw [thick] (v4) to[out=90,in=-90] (v5);
				\draw [line width= 0.15cm,white] (v6) to[out=90,in=-90] (v7);
				\draw [thick] (v6) to[out=90,in=-90] (v7);
				\node at (0,1.5) {};
				\end{tikzpicture}
			\end{array}
		\hspace{-.3em} \big) (h \otimes g) = \nabla_{ST}^{S,T}(\id)(g \otimes h).
		\end{align}

	\proof

	The first claim follows from the fact that, by decomposing the top of each strand (as in, for example, \cite[Lemma 3.3]{hardiman_king}), we have
		\begin{align*}
			\begin{array}{c}
				\begin{tikzpicture}[scale=0.3]
				\coordinate (v1) at (-0.5,0.3) {};
				\coordinate (v2) at (0.5,0.3) {};
				\coordinate (v3) at (0,-0.45) {};
				\coordinate (v10) at (0.5,-1.5) {} {} {};
				\coordinate (v20) at (1.5,-1.5) {} {} {};
				\coordinate (v30) at (1,-2.25) {} {} {};
				\draw[thick] (-1,2) to[out=-90] ($ (v1) + (0.11,-0.1) $);
				\draw[thick] (1,2) to[out=-90,in=45] ($ (v2) + (-0.11,-0.1) $);
				\draw[thick] (0,0) to[out=-90] ($ (v10) + (0.11,-0.1) $);
				\node at (-1,2.75) {$A$};
				\node at (1,2.75) {$A$};
				\node at (1,-3.75) {$A$};
				\node at (3,2.75) {$A$};
				\draw [thick] (3,2) to[out=-90,in=45] ($ (v20) + (-0.11,-0.1) $);
				\draw [thick] ($ (v30) + (0,0.3) $) edge (1,-3);
				\fill [blue] (v1) -- (v2) -- (v3) -- cycle;
				\fill [blue] (v10) -- (v20) -- (v30) -- cycle;
				\end{tikzpicture}
			\end{array}
		= \sum\limits_{\substack{R,S,T \\ f,g,h}} \nabla^{RS,T}_{RST} (\id) \big( \nabla^{R,S}_{RS} (\id) (f \otimes g) \otimes h \big) \circ (f^* \otimes g^* \otimes h^*)
		\end{align*}
	and
		\begin{align*}
			\begin{array}{c}
				\begin{tikzpicture}[scale=0.3,xscale=-1]
				\coordinate (v1) at (-0.5,0.3) {};
				\coordinate (v2) at (0.5,0.3) {};
				\coordinate (v3) at (0,-0.45) {};
				\coordinate (v10) at (0.5,-1.5) {} {} {};
				\coordinate (v20) at (1.5,-1.5) {} {} {};
				\coordinate (v30) at (1,-2.25) {} {} {};
				\draw[thick] (-1,2) to[out=-90] ($ (v1) + (0.11,-0.1) $);
				\draw[thick] (1,2) to[out=-90,in=45] ($ (v2) + (-0.11,-0.1) $);
				\draw[thick] (0,0) to[out=-90] ($ (v10) + (0.11,-0.1) $);
				\node at (-1,2.75) {$A$};
				\node at (1,2.75) {$A$};
				\node at (1,-3.75) {$A$};
				\node at (3,2.75) {$A$};
				\draw [thick] (3,2) to[out=-90,in=45] ($ (v20) + (-0.11,-0.1) $);
				\draw [thick] ($ (v30) + (0,0.3) $) edge (1,-3);
				\fill [blue] (v1) -- (v2) -- (v3) -- cycle;
				\fill [blue] (v10) -- (v20) -- (v30) -- cycle;
				\end{tikzpicture}
			\end{array}
		= \sum\limits_{\substack{R,S,T \\ f,g,h}}  \nabla^{R,ST}_{RST} (\id) \big(  f \otimes \nabla^{S,T}_{ST}(\id)(g \otimes h) \big) \circ (f^* \otimes g^* \otimes h^*).
		\end{align*}
	Similarly the second claim follows from
		\begin{align*}
			\begin{array}{c}
				\begin{tikzpicture}[scale=0.3]
				\coordinate (v1) at (-0.5,0.3) {};
				\coordinate (v2) at (0.5,0.3) {};
				\coordinate (v3) at (0,-0.45) {};
				\draw[thick] (-1,1.5) to[out=-90] ($ (v1) + (0.11,-0.1) $);
				\draw[thick] (1,2) to[out=-90,in=45] ($ (v2) + (-0.11,-0.1) $);
				\draw[thick] (0,0) -- (0,-2);
				\fill [blue] (v1) -- (v2) -- (v3) -- cycle;
				\draw[thick, fill=white] (-1,1.5) circle (0.2);
				\node at (1,2.8) {$A$};
				\node at (0,-3) {$A$};
				\end{tikzpicture}
			\end{array}
		= \sum\limits_{S,g} \nabla^{\tid,S}_S (\id) (u \otimes g) \circ g^* \quad \text{and} \quad
			\begin{array}{c}
				\begin{tikzpicture}[scale=0.3,xscale=-1]
				\coordinate (v1) at (-0.5,0.3) {};
				\coordinate (v2) at (0.5,0.3) {};
				\coordinate (v3) at (0,-0.45) {};
				\draw[thick] (-1,1.5) to[out=-90] ($ (v1) + (0.11,-0.1) $);
				\draw[thick] (1,2) to[out=-90,in=45] ($ (v2) + (-0.11,-0.1) $);
				\draw[thick] (0,0) -- (0,-2);
				\fill [blue] (v1) -- (v2) -- (v3) -- cycle;
				\draw[thick, fill=white] (-1,1.5) circle (0.2);
				\node at (1,2.8) {$A$};
				\node at (0,-3) {$A$};
				\end{tikzpicture}
			\end{array}
		= \sum\limits_{S,g} \nabla^{S,\tid}_S (\id) (g \otimes u) \circ g^*.
		\end{align*}
	and the third claim from
		\begin{align*}
		 	\begin{array}{c}
				\begin{tikzpicture}[scale=0.3,every node/.style={inner sep=0,outer sep=-1}]
				\coordinate (v1) at (-0.5,0.3) {};
				\coordinate (v2) at (0.5,0.3) {};
				\coordinate (v3) at (0,-0.45) {};
				\draw[thick] (-0.75,1) node (v4) {} to[out=-90] ($ (v1) + (0.11,-0.1) $);
				\draw[thick] (0.75,1) node (v6) {} to[out=-90,in=45] ($ (v2) + (-0.11,-0.1) $);
				\draw[thick] (0,0) -- (0,-1);
				\fill [blue] (v1) -- (v2) -- (v3) -- cycle;
				\node (v7) at (-1,3) {};
				\node (v5) at (1,3) {};
				\draw [thick] (v4) to[out=90,in=-90] (v5);
				\draw [line width= 0.2cm,white] (v6) to[out=90,in=-90] (v7);
				\draw [thick] (v6) to[out=90,in=-90] (v7);
				\node at (-1,4) {$ A $};
				\node at (1,4) {$ A $};
				\node at (0,-2) {$ A $};
				\end{tikzpicture}
			\end{array}
 		= \sum\limits_{\substack{S,T\\g,h}}
			\begin{array}{c}
				\begin{tikzpicture}[scale=0.25,every node/.style={inner sep=0,outer sep=-1}]
				\coordinate (v1) at (-0.425,-3) {} {} {} {} {} {};
				\coordinate (v2) at (0.575,-3) {} {} {} {} {} {};
				\coordinate (v3) at (0.075,-3.75) {} {} {} {};
				\draw[thick] (1.25,-1) to[out=-90] ($ (v1) + (0.11,-0.1) $);
				\draw[line width= 0.15cm,white] (-1.25,-1) to[out=-90,in=45] ($ (v2) + (-0.11,-0.1) $);
				\draw[thick] (-1.25,-1) to[out=-90,in=45] ($ (v2) + (-0.11,-0.1) $);
				\draw[thick] (0.075,-3.3) -- (0.075,-5);
				\fill [blue] (v1) -- (v2) -- (v3) -- cycle;
				\node at (2.5,1.25) {$  T$};
				\node at (-0.25,1.25) {$  S$};
				\node (v4) at (1.25,0) {};
				\node (v6) at (-1.25,-0.25) {};
				\node (v7) at (-1.25,4.25) {};
				\node (v5) at (1.25,4.25) {};
				\draw [thick] (v4) to[out=90,in=-90] (v5);
				\draw [line width= 0.15cm,white] (v6) to[out=90,in=-90] (v7);
				\draw [thick] (v6) to[out=90,in=-90] (v7);
				\node [draw,outer sep=0,inner sep=1,minimum width = 13,minimum height = 12,fill=white] at (-1.25,-0.5) {$ g $};
				\node [draw,outer sep=0,inner sep=1,minimum width = 13,minimum height = 12,fill=white] at (1.25,-0.5) {$ h $};
				\node [draw,outer sep=0,inner sep=1,minimum width = 13,minimum height = 12,fill=white] at (-1.25,3) {$ g^* $};
				\node [draw,outer sep=0,inner sep=1,minimum width = 13,minimum height = 12,fill=white] at (1.25,3) {$ h^* $};
				\end{tikzpicture}
			\end{array}
		= \sum\limits_{\substack{S,T\\g,h}}
			\begin{array}{c}
				\begin{tikzpicture}[scale=0.25,every node/.style={inner sep=0,outer sep=-1}]
				\coordinate (v1) at (-0.425,-3) {} {} {} {} {} {};
				\coordinate (v2) at (0.575,-3) {} {} {} {} {} {};
				\coordinate (v3) at (0.075,-3.75) {} {} {} {};
				\draw[thick] (-1.25,-1.5) to[out=-90] ($ (v1) + (0.11,-0.1) $);
				\draw[thick] (1.25,-1.5) to[out=-90,in=45] ($ (v2) + (-0.11,-0.1) $);
				\draw[thick] (0.075,-3.3) -- (0.075,-5);
				\fill [blue] (v1) -- (v2) -- (v3) -- cycle;
				\node (v4) at (-1.25,-0.5) {};
				\node (v6) at (1.25,-0.5) {};
				\node (v7) at (-1.25,2) {};
				\node (v5) at (1.25,2) {};
				\draw [thick] (v4) to[out=90,in=-90] (v5);
				\draw [line width= 0.15cm,white] (v6) to[out=90,in=-90] (v7);
				\draw [thick] (v6) to[out=90,in=-90] (v7);
				\node (v8) at (-1.25,4.25) {};
				\node (v9) at (1.25,4.25) {};
				\draw [thick] (v8) edge (v7);
				\draw [thick] (v9) edge (v5);
				\node [draw,outer sep=0,inner sep=1,minimum width = 13,minimum height = 12,fill=white] at (-1.25,-1) {$ h $};
				\node [draw,outer sep=0,inner sep=1,minimum width = 13,minimum height = 12,fill=white] at (1.25,-1) {$ g $};
				\node [draw,outer sep=0,inner sep=1,minimum width = 13,minimum height = 12,fill=white] at (-1.25,2.5) {$ g^* $};
				\node [draw,outer sep=0,inner sep=1,minimum width = 13,minimum height = 12,fill=white] at (1.25,2.5) {$ h^* $};
				\end{tikzpicture}
			\end{array}
		= \sum\limits_{\substack{S,T\\g,h}} \nabla^{T,S}_{ST}\big( \hspace{-.3em}
			\begin{array}{c}
				\begin{tikzpicture}[scale=0.25,every node/.style={inner sep=0,outer sep=-1}]
				\node (v4) at (-0.5,-1) {};
				\node (v6) at (1,-1) {};
				\node (v7) at (-0.5,0.5) {};
				\node (v5) at (1,0.5) {};
				\draw [thick] (v4) to[out=90,in=-90] (v5);
				\draw [line width= 0.15cm,white] (v6) to[out=90,in=-90] (v7);
				\draw [thick] (v6) to[out=90,in=-90] (v7);
				\node at (0,1.5) {};
				\end{tikzpicture}
			\end{array}
		\hspace{-.3em} \big) (h \otimes g) \circ (g^* \otimes h^*).
		\end{align*}
    \endproof

	\end{LEMMA}

    \begin{LEMMA} \label{lem:dualizing_maps}

	Let $ \B $ be a spherical fusion category and let $ A $ be an object in $ \B $ together with a collection of perfect pairings
    	\begin{align*}
    	\langle \Blank, \Blank \rangle_S \colon A_S \otimes A_{S^\vee} \to \fld
    	\end{align*}
    for all $ S \in \Irr(\B) $. Let $ c $ be a map from $ \Irr(\B) $ to $ \fld \setminus \{ 0 \} $. We consider the morphisms
        \begin{align*}
       		\begin{array}{c}
				\begin{tikzpicture}[scale=0.25,every node/.style={inner sep=0,outer sep=-1},yscale=-1]
				\node (v7) at (-1,2) {};
				\node (v9) at (1,2) {};
				\node (v1) at (0,1) {};
				\draw [thick] (v7) to[out=-90,in=180] (v1);
				\draw [thick] (v9) to[out=-90,in=0] (v1);
				\node at (-1,2.8) {$ \scriptstyle A $};
				\node at (1,2.8) {$ \scriptstyle A $};
				\end{tikzpicture}
       		\end{array}
       	= \sum_{S,b} c(S)
        	\begin{array}{c}
        	    \begin{tikzpicture}[scale=0.25,every node/.style={inner sep=0,outer sep=-1}]
	        	\node [draw,outer sep=0,inner sep=1,minimum width = 10, minimum height = 11, fill=white] (v7) at (-1,2) {$\scriptstyle b' $};
	        	\node [draw,outer sep=0,inner sep=1,minimum width = 10, minimum height = 11,fill=white] (v9) at (1,2) {$\scriptstyle b $};
	        	\node (v1) at (0,3.5) {};
	        	\draw [thick] (v7) to[out=90,in=180] (v1);
	        	\draw [thick] (v9) to[out=90,in=0] (v1);
	        	\draw [thick] (v9) -- (1,0);
	        	\draw [thick] (v7) -- (-1,0);
	        	\node at (1.5,3.5) {$ \scriptstyle S $};
	        	\node at (1.6,0.4) {$ \scriptstyle A $};
	        	\node at (-1.7,0.4) {$ \scriptstyle A $};
	        	\node at (-1.3,3.5) {$ \scriptstyle S^\vee $};
	        	\end{tikzpicture}
        	\end{array}
        \end{align*}
    and
        \begin{align*}
       		\begin{array}{c}
				\begin{tikzpicture}[scale=0.25,every node/.style={inner sep=0,outer sep=-1}]
				\node (v7) at (-1,2) {};
				\node (v9) at (1,2) {};
				\node (v1) at (0,1) {};
				\draw [thick] (v7) to[out=-90,in=180] (v1);
				\draw [thick] (v9) to[out=-90,in=0] (v1);
				\node at (-1,2.8) {$ \scriptstyle A $};
				\node at (1,2.8) {$ \scriptstyle A $};
				\end{tikzpicture}
			\end{array}
		= \sum_{S,b} \frac{1}{c(S)}
			\begin{array}{c}
				\begin{tikzpicture}[scale=0.25,every node/.style={inner sep=0,outer sep=-1}]
				\node [draw,outer sep=0,inner sep=1,minimum width = 12, minimum height = 11, fill=white] (v7) at (-1,2) {$\scriptstyle b^* $};
				\node [draw,outer sep=0,inner sep=1,minimum width = 12, minimum height = 11, fill=white] (v9) at (1,2) {$\scriptstyle b'^* $};
				\node (v1) at (0,0.5) {};
				\draw [thick] (v7) to[out=-90,in=180] (v1);
				\draw [thick] (v9) to[out=-90,in=0] (v1);
				\draw [thick] (v9) -- (1,4);
				\draw [thick] (v7) -- (-1,4);
				\node at (1.7,0.5) {$ \scriptstyle S^\vee $};
				\node at (1.8,3.5) {$ \scriptstyle A $};
				\node at (-1.9,3.5) {$ \scriptstyle A $};
				\node at (-1.5,0.43) {$ \scriptstyle S $};
				\end{tikzpicture}
			\end{array}
        \end{align*}
	where $ \{b\} $ is a basis of $ A_S $ and $ \{b'\} $ is the corresponding dual basis of $ A_{S^\vee} $ with respect to $ \langle \Blank, \Blank \rangle_S $. Then $ \left(A,
    	\begin{array}{c}
			\begin{tikzpicture}[scale=0.25,every node/.style={inner sep=0,outer sep=-1},yscale=-1]
			\node (v7) at (-1,2) {};
			\node (v9) at (1,2) {};
			\node (v1) at (0,1) {};
			\draw [thick] (v7) to[out=-90,in=180] (v1);
			\draw [thick] (v9) to[out=-90,in=0] (v1);
			\node at (-1,2.8) {$ \scriptstyle A $};
			\node at (1,2.8) {$ \scriptstyle A $};
			\end{tikzpicture}
		\end{array}
	,
		\begin{array}{c}
			\begin{tikzpicture}[scale=0.25,every node/.style={inner sep=0,outer sep=-1},yscale=1]
			\node (v7) at (-1,2) {};
			\node (v9) at (1,2) {};
			\node (v1) at (0,1) {};
			\draw [thick] (v7) to[out=-90,in=180] (v1);
			\draw [thick] (v9) to[out=-90,in=0] (v1);
			\node at (-1,2.8) {$ \scriptstyle A $};
			\node at (1,2.8) {$ \scriptstyle A $};
			\end{tikzpicture}
		\end{array}
	\right) $ is a dual object to $ A $. Furthermore, with respect to this duality, we have
		\begin{align} \label{eq:relationship_between_dualities}
		\langle f , (g^*)^\vee \rangle = c(S) g^*(f)
		\end{align}
	for all $ f,g \in A_S $.
    \proof

    We have
   		\begin{align*}
   			\begin{array}{c}
   				\begin{tikzpicture}[scale=0.25,every node/.style={inner sep=0,outer sep=-1}]
   				\node (v7) at (-1,2) {};
   				\node (v9) at (1,2) {};
   				\node (v1) at (0,1) {};
   				\draw [thick] (v7) to[out=-90,in=180] (v1);
   				\draw [thick] (v9) to[out=-90,in=0] (v1);
   				\node at (-3.6,0.8) {$ \scriptstyle A $};
   				\node at (0.4,3.2) {$ \scriptstyle A $};
   				\node (v2) at (-2,3) {};
   				\node (v3) at (-3,2) {};
   				\node (v5) at (1,3.5) {};
   				\node (v4) at (-3,0.5) {};
   				\draw [thick] (v7) to[out=90,in=0] (v2);
   				\draw [thick] (v2) to[out=180,in=90] (v3);
   				\draw [thick] (v3) edge (v4);
   				\draw [thick] (v9) edge (v5);
   				\end{tikzpicture}
   			\end{array}
   		= \sum\limits_{S,b} \frac{c(S)}{c(S)}
   			\begin{array}{c}
   				\begin{tikzpicture}[scale=0.25,every node/.style={inner sep=0,outer sep=-1}]
   				\node [draw,outer sep=0,inner sep=1,minimum width = 12, minimum height = 11, fill=white] (v7) at (-1,2) {$\scriptstyle b^* $};
   				\node [draw,outer sep=0,inner sep=1,minimum width = 12, minimum height = 11, fill=white] (v9) at (1,2) {$\scriptstyle b'^* $};
   				\node [draw,outer sep=0,inner sep=1,minimum width = 12, minimum height = 11, fill=white] (v70) at (-1,4.5) {$\scriptstyle b $};
   				\node [draw,outer sep=0,inner sep=1,minimum width = 12, minimum height = 11, fill=white] (v90) at (-3,4.5) {$\scriptstyle b' $};
   				\draw [thick] (v7) to[out=-90,in=-90] (v9);
   				\draw [thick] (v9) -- (1,6.5);
   				\draw [thick] (v7) -- (v70);
   				\node at (1.7,0.5) {$ \scriptstyle S^\vee $};
   				\node at (-1.5,0.43) {$ \scriptstyle S $};
   				\node at (-0.5,6.13) {$ \scriptstyle S $};
   				\node at (-3,6.2) {$ \scriptstyle S^\vee $};
   				\node at (1.75,6) {$ \scriptstyle A $};
   				\node at (-3.7,0.3) {$ \scriptstyle A $};
   				\draw [thick] (v90) to[out=90,in=90] (v70);
   				\node (v1) at (-3,0) {};
   				\draw [thick] (v90) edge (v1);
   				\end{tikzpicture}
   			\end{array}
   		= \sum\limits_{S,b}
   			\begin{array}{c}
   				\begin{tikzpicture}[scale=0.25,every node/.style={inner sep=0,outer sep=-1}]
   				\node (v7) at (-1,1.25) {};
   				\node [draw,outer sep=0,inner sep=1,minimum width = 12, minimum height = 11, fill=white] (v9) at (1,2) {$\scriptstyle b'^* $};
   				\node  (v70) at (-1,5.25) {};
   				\node [draw,outer sep=0,inner sep=1,minimum width = 12, minimum height = 11, fill=white] (v90) at (-3,4.5) {$\scriptstyle b' $};
   				\draw [thick] (v7) to[out=-90,in=-90] (v9);
   				\draw [thick] (v9) -- (1,6.5);
   				\draw [thick] (v7) -- (v70);
   				\node at (1.7,0.5) {$ \scriptstyle S^\vee $};
   				\node at (-3,6.2) {$ \scriptstyle S^\vee $};
   				\node at (1.75,6) {$ \scriptstyle A $};
   				\node at (-3.7,0.3) {$ \scriptstyle A $};
   				\draw [thick] (v90) to[out=90,in=90] (v70);
   				\node (v1) at (-3,0) {};
   				\draw [thick] (v90) edge (v1);
   				\end{tikzpicture}
   			\end{array}
   		= \sum\limits_{S,b}
   			\begin{array}{c}
				\begin{tikzpicture}[scale=0.25,every node/.style={inner sep=0,outer sep=-1}]
				\node [draw,outer sep=0,inner sep=1,minimum width = 12, minimum height = 11, fill=white] (v9) at (1,4.5) {$\scriptstyle b'^* $};
				\node [draw,outer sep=0,inner sep=1,minimum width = 12, minimum height = 11, fill=white] (v90) at (1,1.5) {$\scriptstyle b' $};
				\draw [thick] (v9) -- (1,6.5);
				\node at (2,3) {$ \scriptstyle S^\vee $};
				\node at (1.75,6) {$ \scriptstyle A $};
				\node at (1.75,0) {$ \scriptstyle A $};
				\node (v1) at (1,-0.5) {};
				\draw [thick] (v90) edge (v1);
				\draw [thick] (v9) edge (v90);
				\end{tikzpicture}
			\end{array}
		=
			\begin{array}{c}
				\begin{tikzpicture}[scale=0.25,every node/.style={inner sep=0,outer sep=-1}]
				\node at (1.75,6) {$ \scriptstyle A $};
				\node at (1.75,0) {$ \scriptstyle A $};
				\node (v1) at (1,-0.5) {};
				\draw [thick] (1,6.5) edge (v1);
				\end{tikzpicture}
			\end{array}
   		\end{align*}
	and, in the same way, we also have
		\begin{align*}
		   	\begin{array}{c}
				\begin{tikzpicture}[scale=0.25,every node/.style={inner sep=0,outer sep=-1},xscale=-1]
				\node (v7) at (-1,2) {};
				\node (v9) at (1,2) {};
				\node (v1) at (0,1) {};
				\draw [thick] (v7) to[out=-90,in=180] (v1);
				\draw [thick] (v9) to[out=-90,in=0] (v1);
				\node at (-3.6,0.8) {$ \scriptstyle A $};
				\node at (0.4,3.2) {$ \scriptstyle A $};
				\node (v2) at (-2,3) {};
				\node (v3) at (-3,2) {};
				\node (v5) at (1,3.5) {};
				\node (v4) at (-3,0.5) {};
				\draw [thick] (v7) to[out=90,in=0] (v2);
				\draw [thick] (v2) to[out=180,in=90] (v3);
				\draw [thick] (v3) edge (v4);
				\draw [thick] (v9) edge (v5);
				\end{tikzpicture}
			\end{array}
		=
		   	\begin{array}{c}
				\begin{tikzpicture}[scale=0.25,every node/.style={inner sep=0,outer sep=-1},yscale = 0.4]
				\node at (1.75,6) {$ \scriptstyle A $};
				\node at (1.75,0) {$ \scriptstyle A $};
				\node (v1) at (1,-0.5) {};
				\draw [thick] (1,6.5) edge (v1);
				\end{tikzpicture}
			\end{array}.
		\end{align*}
	To prove the second claim we simply compute
		\begin{align*}
		\langle f , (g^*)^\vee \rangle &= \left\langle f ,
			\begin{array}{c}
				\begin{tikzpicture}[scale=0.25,every node/.style={inner sep=0,outer sep=-1}]
				\node [draw,outer sep=0,inner sep=1,minimum width = 12, minimum height = 11, fill=white] (v7) at (-1,2) {$\scriptstyle g^* $};
				\node (v9) at (0.5,1.25) {};
				\node (v90) at (-2.5,2.75) {};
				\draw [thick] (v7) to[out=-90,in=-90] (v9);
				\draw [thick] (v9) -- (0.5,4.25);
				\node at (1.5,4) {$ \scriptstyle S^\vee $};
				\node at (-3.2,0.3) {$ \scriptstyle A $};
				\draw [thick] (v90) to[out=90,in=90] (v7);
				\node (v1) at (-2.5,0) {};
				\draw [thick] (v90) edge (v1);
				\end{tikzpicture}
			\end{array}
		\right\rangle \\
		&= c(S) \sum\limits_{b} \left\langle f ,
			\begin{array}{c}
				\begin{tikzpicture}[scale=0.25,every node/.style={inner sep=0,outer sep=-1}]
				\node [draw,outer sep=0,inner sep=1,minimum width = 12, minimum height = 11, fill=white] (v7) at (-1,2) {$\scriptstyle g^* $};
				\node (v9) at (1,1.25) {};
				\node [draw,outer sep=0,inner sep=1,minimum width = 12, minimum height = 11, fill=white] (v70) at (-1,4) {$\scriptstyle b $};
				\node [draw,outer sep=0,inner sep=1,minimum width = 12, minimum height = 11, fill=white] (v90) at (-3,4) {$\scriptstyle b' $};
				\draw [thick] (v7) to[out=-90,in=-90] (v9);
				\draw [thick] (v9) -- (1,6);
				\draw [thick] (v7) -- (v70);
				\node at (-0.5,5.63) {$ \scriptstyle S $};
				\node at (-3,5.7) {$ \scriptstyle S^\vee $};
				\node at (2,5.5) {$ \scriptstyle S^\vee $};
				\node at (-3.7,0.8) {$ \scriptstyle A $};
				\draw [thick] (v90) to[out=90,in=90] (v70);
				\node (v1) at (-3,0.5) {};
				\draw [thick] (v90) edge (v1);
				\end{tikzpicture}
			\end{array}
		\right\rangle \\
		&= c(S) \sum\limits_{b} g^*(b) \langle f, b' \rangle = c(S) g^*(f).
		\end{align*}
    \endproof

    \end{LEMMA}

	\begin{LEMMA} \label{lem:checking_balanced_condition}

	Let $ \B $ be a spherical fusion category and let $ A $ be an algebra object in $ \B $ (with product $ \nabla $) together with structure maps that make $ A $ self-dual. Then $ A $ satisfies \eqref{eq:unique_coproduct} if and only if
		\begin{align}
			\begin{split} \label{eq_frob_cond}
			d(S)\left( g^* \circ \nabla_S^{R,T^\vee}\left(
					\begin{array}{c}
							\begin{tikzpicture}[scale=0.25,every node/.style={inner sep=0,outer sep=-1}]
							\draw[thick] (0,0) -- (0,-2);
							\draw[thick] (-0.6,0) -- (-0.6,2);
							\draw[thick] (2,0.7)-- (2,-2);
							\draw[thick] (0.6,0.7) to[out=90,in=90] (2,0.7);
							\node [draw,outer sep=0,inner sep=2,minimum width=20,fill=white] at (0,0) {$\scriptstyle \beta $};
							\node at (-0.6,-1.65) {$\scriptscriptstyle R $};
							\node at (1.3,-1.55) {$\scriptscriptstyle T^\vee $};
							\node at (-1.1,1.6) {$\scriptscriptstyle S $};
							\end{tikzpicture}
					\end{array}
			\right) \right)(f \otimes (h^*)^\vee) \\
		 = d(T) \left( h^* \circ \nabla_T^{S^\vee,R}\left(
					\begin{array}{c}
							\begin{tikzpicture}[scale=0.25,every node/.style={inner sep=0,outer sep=-1}, xscale=-1]
							\draw[thick] (0,0) -- (0,-2);
							\draw[thick] (-0.6,0) -- (-0.6,2);
							\draw[thick] (2,0.7)-- (2,-2);
							\draw[thick] (0.6,0.7) to[out=90,in=90] (2,0.7);
							\node [draw,outer sep=0,inner sep=2,minimum width=20,fill=white] at (0,0) {$\scriptstyle \beta $};
							\node at (-0.5,-1.65) {$\scriptscriptstyle R $};
							\node at (1.2,-1.5) {$\scriptscriptstyle S^\vee $};
							\node at (-1.1,1.6) {$\scriptscriptstyle T $};
							\end{tikzpicture}
					\end{array}
			\right) \right)((g^*)^\vee \otimes f)
			\end{split}
		\end{align}
	for all $ R,S,T \in \Irr(\B) $, $ \beta \in \Hom_{\B}(ST,R) $, $ h,f \in A_R $, $ g,h \in A_T $ and $ f,g \in A_S $.

	\proof

	Decomposing the coproduct given by the left-hand side of \eqref{eq:unique_coproduct} gives
		\begin{align*}
			\left(\begin{array}{c}
				\begin{tikzpicture}[scale=0.2]
				\coordinate (v1) at (-0.5,0.3) {};
				\coordinate (v2) at (0.5,0.3) {};
				\coordinate (v3) at (0,-0.45) {};
				\draw[thick] (-1,2) to[out=-90] ($ (v1) + (0.11,-0.1) $);
				\draw[thick] (1,1) to[out=-90,in=45] ($ (v2) + (-0.11,-0.1) $);
				\draw[thick] (1,1) to[out=90,in=90] (3,1);
				\draw[thick] (0,0) -- (0,-2);
				\draw[thick] (3,1) -- (3,-2);
				\fill [blue] (v1) -- (v2) -- (v3) -- cycle;
				\end{tikzpicture}
			\end{array}\right)^R_{S,T} (\beta) (f) &= \sum\limits_{g,h}
			\begin{array}{c}
				\begin{tikzpicture}[scale=0.25,every node/.style={inner sep=0,outer sep=-1}, xscale=1.1,yscale=-1]
				\draw[thick] (-0.8,0) -- (-0.75,-2.5);
				\draw[thick] (0,0) -- (0,2);
				\node (v6) at (-1,-7.25) {};
				\node (v7) at (-1,-5.75) {};
				\node at (0.5,1.7) {$\scriptscriptstyle R $};
				\draw[thick] (0.75,-0.25) to[out=-90,in=90] (1.75,-1.5) node (v5) {};
				\coordinate (v1) at (-0.5,-4.25) {} {} {};
				\coordinate (v2) at (0.5,-4.25) {} {} {};
				\coordinate (v3) at (0,-3.5) {} {} {};
				\draw[thick] (-0.75,-2.5) to[out=-90,in=90] ($ (v3) + (0,-0.1) $);
				\draw[thick] (1.25,-4.75) node (v4) {} to[out=-180,in=-45] ($ (v2) + (-0.11,0.1) $);
				\draw [thick] (v7) to[out=90,in=-135] ($ (v1) + (0.15,0.1) $);
				\fill [blue] (v1) -- (v2) -- (v3) -- cycle;
				\draw [thick] (v4) to[out=0,in=-90] (v5);
				\node at (-0.5,-7) {$ \scriptscriptstyle R$};
				\draw [thick] (v6) edge (v7);
				\node [draw,outer sep=0,inner sep=2,minimum width=15,minimum height= 11,fill=white] at (0,0.4) {$\scriptstyle \beta $};
				\node [draw,outer sep=0,inner sep=1,minimum width = 11,minimum height = 11,fill=white] at (-0.75,-2) {$\scriptstyle g^* $};
				\node [draw,outer sep=0,inner sep=1,minimum width = 11,minimum height = 11,fill=white] at (1.75,-2) {$\scriptstyle h^* $};
				\node [draw,outer sep=0,inner sep=1,minimum width = 11,minimum height = 11,fill=white] at (v7)  {$\scriptstyle f $};
				\end{tikzpicture}
			\end{array}
		g \otimes h
		= \frac{1}{d(R)} \sum\limits_{g,h}
			\begin{array}{c}
				\begin{tikzpicture}[scale=0.25,every node/.style={inner sep=0,outer sep=-1}, xscale=1.1,yscale=-1]
				\draw[thick] (-0.8,0) -- (-0.75,-2.5);
				\draw[thick] (0,0) -- (0,2) node (v10) {};
				\node (v6) at (-1,-7.5) {};
				\node (v7) at (-1,-5.75) {};
				\node at (0.5,1.7) {$\scriptscriptstyle R $};
				\draw[thick] (0.75,-0.25) to[out=-90,in=90] (1.75,-1.5) node (v5) {};
				\coordinate (v1) at (-0.5,-4.25) {} {} {};
				\coordinate (v2) at (0.5,-4.25) {} {} {};
				\coordinate (v3) at (0,-3.5) {} {} {};
				\draw[thick] (-0.75,-2.5) to[out=-90,in=90] ($ (v3) + (0,-0.1) $);
				\draw[thick] (1.25,-4.75) node (v4) {} to[out=-180,in=-45] ($ (v2) + (-0.11,0.1) $);
				\draw [thick] (v7) to[out=90,in=-135] ($ (v1) + (0.15,0.1) $);
				\fill [blue] (v1) -- (v2) -- (v3) -- cycle;
				\draw [thick] (v4) to[out=0,in=-90] (v5);
				\node at (-0.5,-7) {$ \scriptscriptstyle R$};
				\draw [thick] (v6) edge (v7);
				\node [draw,outer sep=0,inner sep=2,minimum width=15,minimum height= 11,fill=white] at (0,0.4) {$\scriptstyle \beta $};
				\node [draw,outer sep=0,inner sep=1,minimum width = 11,minimum height = 11,fill=white] at (-0.75,-2) {$\scriptstyle g^* $};
				\node [draw,outer sep=0,inner sep=1,minimum width = 11,minimum height = 11,fill=white] at (1.75,-2) {$\scriptstyle h^* $};
				\node [draw,outer sep=0,inner sep=1,minimum width = 11,minimum height = 11,fill=white] at (v7)  {$\scriptstyle f $};
				\node (v9) at (-2.5,2) {};
				\node (v8) at (-2.5,-7.5) {};
				\draw [thick] (v8) edge (v9);
				\draw [thick] (v9) to[out=90,in=90] (0,2);
				\draw [thick] (v8) to[out=-90,in=-90] (v6);
				\end{tikzpicture}
			\end{array}
		g \otimes h
	\end{align*}
	\begin{align*}
		&= \frac{1}{d(R)}  \sum\limits_{g,h}
			\begin{array}{c}
				\begin{tikzpicture}[scale=0.25,every node/.style={inner sep=0,outer sep=-1}, xscale=1.1]
				\draw[thick] (-0.75,0) -- (-0.75,-2.75);
				\draw[thick] (-1,0) -- (-1,2) node (v10) {};
				\node (v6) at (-1,-7.5) {};
				\node (v7) at (-1,-5.75) {};
				\node at (-0.25,-1) {$\scriptscriptstyle R $};
				\draw[thick] (2.5,1.5) node (v11) {} to[out=-90,in=90] (2.5,-7) node (v5) {};
				\coordinate (v1) at (0,-4.25) {} {} {} {};
				\coordinate (v2) at (0.5,-3.5) {} {} {} {} {};
				\coordinate (v3) at (-0.5,-3.5) {} {} {} {};
				\draw[thick] (-0.75,-2.75) to[out=-90,in=135] ($ (v3) + (0.15,-0.1) $);
				\draw[thick] (0.875,-3.25) node (v4) {} to[out=180,in=45] ($ (v2) + (-0.15,-0.1) $);
				\draw [thick] (-1,-5) to[out=90,in=-90] ($ (v1) + (0,0.1) $);
				\fill [blue] (v1) -- (v2) -- (v3) -- cycle;
				\draw [thick] (v4) to[out=0,in=90] (1.25,-3.5) node (v14) {};
				\node at (-0.5,-7) {$ \scriptscriptstyle S$};
				\node at (0.75,-7) {$ \scriptscriptstyle T$};
				\draw [thick] (v6) edge (v7);
				\node (v9) at (-2.5,2) {};
				\node (v8) at (-2.5,-7.5) {};
				\draw [thick] (v8) edge (v9);
				\draw [thick] (v9) to[out=90,in=90] (-1,2);
				\draw [thick] (v8) to[out=-90,in=-90] (v6);
				\node (v12) at (-0.25,1.5) {};
				\draw [thick] (v11) to[out=90,in=90] (v12);
				\node (v13) at (-0.25,1) {};
				\draw [thick] (v12) edge (v13);
				\node [draw,outer sep=0,inner sep=2,minimum width=15,minimum height= 11,fill=white] at (-0.75,0.4) {$\scriptstyle \beta $};
				\node [draw,outer sep=0,inner sep=1,minimum width = 11,minimum height = 11, fill=white] at (-0.75,-2.25) {$\scriptstyle f $};
				\node [draw,outer sep=0,inner sep=1,minimum width = 11,minimum height = 11, fill=white] (v15) at (1.25,-5.75) {$\scriptstyle h^* $};
				\node [draw,outer sep=0,inner sep=1,minimum width = 11, minimum height = 11,fill=white] at (v7)  {$\scriptstyle g^* $};
				\draw [thick] (v14) edge (v15);
				\node (v16) at (1.25,-7) {};
				\draw [thick] (v5) to[out=-90,in=-90] (v16);
				\draw [thick] (v15) edge (v16);
				\end{tikzpicture}
			\end{array}
		g \otimes h = \frac{1}{d(R)} \sum\limits_{g,h}
			\begin{array}{c}
				\begin{tikzpicture}[scale=0.25,every node/.style={inner sep=0,outer sep=-1}, xscale=1.1]
				\draw[thick] (-0.75,0.25) to[out=-90,in=90] (-1.25,-1.75);
				\draw[thick] (-1,0.5) -- (-1,2) node (v10) {};
				\node (v6) at (-1,-7.5) {};
				\node (v7) at (-1,-5.75) {};
				\node at (-0.25,-0.5) {$\scriptscriptstyle R $};
				\node (v16) at (1.25,0.25) {};
				\draw[thick] (1.25,2) node (v11) {} to[out=-90,in=90] (v16) {};
				\coordinate (v1) at (0,-4.25) {} {} {} {};
				\coordinate (v2) at (0.5,-3.5) {} {} {} {} {};
				\coordinate (v3) at (-0.5,-3.5) {} {} {} {};
				\draw[thick] (-1.25,-2.75) to[out=-90,in=135] ($ (v3) + (0.15,-0.1) $);
				\draw [thick] (-1,-5) to[out=90,in=-90] ($ (v1) + (0,0.1) $);
				\node (v4) at (1.25,-2.75) {};
				\draw [thick] (v4) to[out=-90,in=45] ($ (v2) + (-0.15,-0.1) $);
				\fill [blue] (v1) -- (v2) -- (v3) -- cycle;
				\node at (-0.5,-7) {$ \scriptscriptstyle S$};
				\draw [thick] (v6) edge (v7);
				\node (v9) at (-2.5,2) {};
				\node (v8) at (-2.5,-7.5) {};
				\draw [thick] (v8) edge (v9);
				\draw [thick] (v9) to[out=90,in=90] (-1,2);
				\draw [thick] (v8) to[out=-90,in=-90] (v6);
				\node (v12) at (-0.25,2) {};
				\draw [thick] (v11) to[out=90,in=90] (v12);
				\node (v13) at (-0.25,1.5) {};
				\draw [thick] (v12) edge (v13);
				\node [draw,outer sep=0,inner sep=2,minimum width=15,minimum height= 11,fill=white] at (-0.75,0.9) {$\scriptstyle \beta $};
				\node [draw,outer sep=0,inner sep=1,minimum width = 11,minimum height = 11, fill=white] at (-1.25,-2.25) {$\scriptstyle f $};
				\node [draw,outer sep=0,inner sep=1,minimum width = 11,minimum height = 11, fill=white] (v15) at (1.25,-2.25) {$\scriptstyle (h^*)^\vee $};
				\node [draw,outer sep=0,inner sep=1,minimum width = 11, minimum height = 11,fill=white] at (v7)  {$\scriptstyle g^* $};
				\draw [thick] (v15) edge (v16);
				\end{tikzpicture}
			\end{array}
		g \otimes h \\
		&= \frac{d(S)}{d(R)} \sum\limits_{g,h}
			\begin{array}{c}
				\begin{tikzpicture}[scale=0.25,every node/.style={inner sep=0,outer sep=-1}, xscale=1.1]
				\draw[thick] (-0.75,0.25) to[out=-90,in=90] (-1.25,-1.75);
				\draw[thick] (-1,0.5) -- (-1,2.75) node (v10) {};
				\node (v6) at (-1,-7.5) {};
				\node (v7) at (-1,-5.75) {};
				\node at (-0.25,-0.5) {$\scriptscriptstyle R $};
				\node (v16) at (1.25,0.25) {};
				\draw[thick] (1.25,2) node (v11) {} to[out=-90,in=90] (v16) {};
				\coordinate (v1) at (0,-4.25) {} {} {} {};
				\coordinate (v2) at (0.5,-3.5) {} {} {} {} {};
				\coordinate (v3) at (-0.5,-3.5) {} {} {} {};
				\draw[thick] (-1.25,-2.75) to[out=-90,in=135] ($ (v3) + (0.15,-0.1) $);
				\draw [thick] (-1,-5) to[out=90,in=-90] ($ (v1) + (0,0.1) $);
				\node (v4) at (1.25,-2.75) {};
				\draw [thick] (v4) to[out=-90,in=45] ($ (v2) + (-0.15,-0.1) $);
				\fill [blue] (v1) -- (v2) -- (v3) -- cycle;
				\node at (-1.5,-7.25) {$ \scriptscriptstyle S$};
				\node at (-1.5,2.5) {$ \scriptscriptstyle S$};
				\draw [thick] (v6) edge (v7);
				\node (v12) at (-0.25,2) {};
				\draw [thick] (v11) to[out=90,in=90] (v12);
				\node (v13) at (-0.25,1.5) {};
				\draw [thick] (v12) edge (v13);
				\node [draw,outer sep=0,inner sep=2,minimum width=15,minimum height= 11,fill=white] at (-0.75,0.9) {$\scriptstyle \beta $};
				\node [draw,outer sep=0,inner sep=1,minimum width = 11,minimum height = 11, fill=white] at (-1.25,-2.25) {$\scriptstyle f $};
				\node [draw,outer sep=0,inner sep=1,minimum width = 11,minimum height = 11, fill=white] (v15) at (1.25,-2.25) {$\scriptstyle (h^*)^\vee $};
				\node [draw,outer sep=0,inner sep=1,minimum width = 11, minimum height = 11,fill=white] at (v7)  {$\scriptstyle g^* $};
				\draw [thick] (v15) edge (v16);
				\end{tikzpicture}
			\end{array}
		g \otimes h\\
		&= \frac{d(S)}{d(R)} \sum\limits_{g,h} \left( g^* \circ \nabla_S^{R,T^\vee}\left(
			\begin{array}{c}
				\begin{tikzpicture}[scale=0.25,every node/.style={inner sep=0,outer sep=-1}]
				\draw[thick] (0,0) -- (0,-2);
				\draw[thick] (-0.6,0) -- (-0.6,2);
				\draw[thick] (2,0.7)-- (2,-2);
				\draw[thick] (0.6,0.7) to[out=90,in=90] (2,0.7);
				\node [draw,outer sep=0,inner sep=2,minimum width=20,fill=white] at (0,0) {$\scriptstyle \beta $};
				\node at (-0.6,-1.65) {$\scriptscriptstyle R $};
				\node at (1.3,-1.55) {$\scriptscriptstyle T^\vee $};
				\node at (-1.1,1.6) {$\scriptscriptstyle S $};
				\end{tikzpicture}
			\end{array}
		\right) \right)(f \otimes (h^*)^\vee) \ g \otimes h
		\end{align*}
	In an analogous way, we also have
		\begin{align*}
		\left(
			\begin{array}{c}
				\begin{tikzpicture}[scale=0.2,xscale=-1]
				\coordinate (v1) at (-0.5,0.3) {};
				\coordinate (v2) at (0.5,0.3) {};
				\coordinate (v3) at (0,-0.45) {};
				\draw[thick] (-1,2) to[out=-90] ($ (v1) + (0.11,-0.1) $);
				\draw[thick] (1,1) to[out=-90,in=45] ($ (v2) + (-0.11,-0.1) $);
				\draw[thick] (1,1) to[out=90,in=90] (3,1);
				\draw[thick] (0,0) -- (0,-2);
				\draw[thick] (3,1) -- (3,-2);
				\fill [blue] (v1) -- (v2) -- (v3) -- cycle;
				\end{tikzpicture}
			\end{array}\right)^R_{S,T} (\beta) (f) = \frac{d(T)}{d(R)} \sum\limits_{g,h} \left( h^* \circ \nabla_T^{S^\vee,R}\left(
			\begin{array}{c}
				\begin{tikzpicture}[scale=0.25,every node/.style={inner sep=0,outer sep=-1}, xscale=-1]
				\draw[thick] (0,0) -- (0,-2);
				\draw[thick] (-0.6,0) -- (-0.6,2);
				\draw[thick] (2,0.7)-- (2,-2);
				\draw[thick] (0.6,0.7) to[out=90,in=90] (2,0.7);
				\node [draw,outer sep=0,inner sep=2,minimum width=20,fill=white] at (0,0) {$\scriptstyle \beta $};
				\node at (-0.5,-1.65) {$\scriptscriptstyle R $};
				\node at (1.2,-1.5) {$\scriptscriptstyle S^\vee $};
				\node at (-1.1,1.6) {$\scriptscriptstyle T $};
				\end{tikzpicture}
			\end{array}
			\right) \right)((g^*)^\vee \otimes f) \ g \otimes h
		\end{align*}
	which proves the proposition. \endproof

	\end{LEMMA}

\bibliographystyle{alpha}
\bibliography{extending_the_trace}

\end{document}